\documentclass[10pt]{preprint}

\usepackage[margin=3.5cm]{geometry}
\usepackage[utf8x]{inputenc}
\usepackage[english]{babel}
\usepackage[full]{textcomp}

\usepackage[osf]{newtxtext}

\usepackage{amssymb}
\usepackage{amsmath}
\usepackage{amsthm}

\usepackage{mhequ}
\usepackage{booktabs}
\usepackage{tikz}
\usepackage{mathrsfs}
\usepackage[noadjust]{cite}
\usepackage{microtype}
\usepackage{comment}
\usepackage{slashed}
\usepackage{mathtools}
\usepackage{centernot}
\usepackage{footnote}
\usepackage{enumerate}
\usepackage[shortlabels]{enumitem}
\usepackage{stackrel}
\usepackage{longtable}
\usepackage{cprotect}
\usepackage{xstring}
\usepackage{calc}
\usepackage{bbm}
\usepackage{dsfont}

\usepackage[colorlinks=true, pdfstartview=FitV, linkcolor=colorLink, citecolor=colorCite, urlcolor=colorLink, linktocpage=true]{hyperref}

\usepackage[capitalize]{cleveref}

\usepackage{upgreek}
\usepackage{parskip}
\usepackage{graphicx}
\usepackage{orcidlink}

\makeatletter
\renewcommand{\paragraph}{%
\@startsection{paragraph}{4}%
{\z@}{1.5ex \@plus 1.5ex \@minus .2ex}{-0.7em}%
{\normalfont\normalsize\bfseries}%
}
\makeatother

\makeatletter
\def\thm@space@setup{%
  \thm@preskip=\parskip \thm@postskip=0pt
}
\makeatother

\setlist[itemize]{leftmargin=5mm}

\DeclareSymbolFont{timesoperators}{T1}{ptm}{m}{n}
\SetSymbolFont{timesoperators}{bold}{T1}{ptm}{b}{n}
\DeclareMathAlphabet{\mathbb}{U}{jkpsyb}{m}{n}
\SetMathAlphabet{\mathbb}{bold}{U}{jkpsyb}{bx}{n}

\allowdisplaybreaks
\definecolor{colorLink}{RGB}{0,100,162}
\definecolor{colorCite}{RGB}{8,124,100}

\def\R{\mathbb{R}}

\def\RR{\mathbb{R}}
\def\R{\mathbb{R}}
\def\PP{\mathbb{P}}
\def\cP{\mathcal{P}}
\def\EE{\mathbb{E}}
\def\uE{\mathsf{E}}
\def\ue{\mathsf{e}}
\def\uf{\mathsf{f}}
\def\ud{\mathsf{d}}
\def\to{\rightarrow}
\def\a{\alpha}

\def\d{\text{d}}
\def\d{{\rm{d}}}
\def\eps{\varepsilon}
\newcommand{\bb}[1]{\mathbb{#1}}
\newcommand{\ca}[1]{\mathcal{#1}}
\newcommand{\scr}[1]{\mathscr{#1}}
\newcommand{\worknote}[1]{}

\let\eps\upvarepsilon
\let\alpha\upalpha
\let\beta\upbeta
\let\delta\updelta
\let\gamma\upgamma
\let\mu\upmu
\let\eta\upeta
\let\nu\upnu
\let\rho\uprho
\let\chi\upchi
\let\xi\upxi
\let\zeta\upzeta
\let\tau\uptau
\let\varphi\upvarphi
\let\lambda\uplambda
\let\theta\uptheta
\let\pi\uppi
\let\Upsilon\Upupsilon
\let\Theta\Uptheta
\let\Psi\Uppsi
\let\Xi\Upxi

\def\dash{\leavevmode\unskip\kern0.18em--\penalty\exhyphenpenalty\kern0.18em}
\def\slash{\leavevmode\unskip\kern0.15em/\penalty\exhyphenpenalty\kern0.15em}

\makeatletter
\renewcommand{\operator@font}{\mathgroup\symtimesoperators}
\makeatother

\makeatletter
\DeclareRobustCommand{\TitleEquation}[2]{\texorpdfstring{\StrLeft{\f@series}{1}[\@firstchar]$\if%
b\@firstchar\boldsymbol{#1}\else#1\fi$}{#2}}
\makeatother

\makeatletter
\newcommand{\pushright}[1]{\ifmeasuring@#1\else\omit\hfill$\displaystyle#1$\fi\ignorespaces}
\newcommand{\pushleft}[1]{\ifmeasuring@#1\else\omit$\displaystyle#1$\hfill\fi\ignorespaces}
\makeatother

\renewcommand{\bar}{\overline}

\renewcommand{\tilde}{\widetilde}

\newcommand{\bfX}{\mathbf{X}} 

\makeatletter
\newcommand{\oset}[3][0ex]{%
  \mathrel{\mathop{#3}\limits^{
    \vbox to#1{\kern-2\ex@
    \hbox{$\scriptstyle#2$}\vss}}}}
\makeatother

\usetikzlibrary{shapes.misc}
\usetikzlibrary{shapes.symbols}
\usetikzlibrary{shapes.geometric}
\usetikzlibrary{decorations}
\usetikzlibrary{decorations.markings}
\usetikzlibrary{calc}
\usetikzlibrary{external}
\usetikzlibrary{arrows}
\usetikzlibrary{patterns}

\theoremstyle{plain}

\newtheorem{theorem}{Theorem}[section]
\crefname{theorem}{Theorem}{Theorems}

\newtheorem{corollary}[theorem]{Corollary}
\crefname{corollary}{Corollary}{Corollaries}

\newtheorem{lemma}[theorem]{Lemma}
\crefname{lemma}{Lemma}{Lemmas}

\newtheorem{proposition}[theorem]{Proposition}
\crefname{proposition}{Proposition}{Propositions}

\crefname{claim}{Claim}{Claims}

\crefname{theoremA}{Theorem}{Theorems}

\crefname{propositionA}{Proposition}{Propositions}

\theoremstyle{definition}

\newtheorem{definition}[theorem]{Definition}
\crefname{definition}{Definition}{Definitions}

\crefname{notation}{Notation}{Notations}

\crefname{acknowledgements}{Acknowledgements}{Acknowledgements}

\newtheorem{assumption}[theorem]{Assumption}
\crefname{assumption}{Assumption}{Assumptions}

\newtheorem{remark}[theorem]{Remark}
\crefname{remark}{Remark}{Remarks}

\crefname{observation}{Observation}{Observations}

\crefname{manualassumptioninner}{Assumption}{Assumptions}

\crefname{customthm}{Theorem}{Theorems}

\numberwithin{equation}{section}

\usepackage[textsize=scriptsize,linecolor=white,bordercolor=red,backgroundcolor=white]{todonotes}
\usepackage{autonum}

\allowdisplaybreaks

\begin{document}
\title{Random dynamical systems for McKean--Vlasov SDEs via rough path theory}

\author{Benjamin Gess$^{1}$, Rishabh S.\ Gvalani$^{2}$, Shanshan Hu$^{3}$}

\institute{
    Institut f\"{u}r Mathematik, Technische Universit\"{a}t Berlin \& Max–Planck–Institute for Mathematics in the Sciences, Leipzig, Email: \href{mailto:benjamin.gess@tu-berlin.de}{\color{black} \texttt{benjamin.gess@tu-berlin.de}} 
    \and
    D-MATH, ETH Z\"{u}rich, Email: \href{mailto:rgvalani@ethz.ch}{\color{black}\texttt{rgvalani@ethz.ch}}
    \and
    Institut f\"{u}r Mathematik, Technische Universit\"{a}t Berlin, Email: \href{mailto:shanshan.hu@tu-berlin.de}{\color{black} \texttt{shanshan.hu@tu-berlin.de}} 
    }

\maketitle

\begin{abstract}
    The existence of random dynamical systems for McKean--Vlasov SDEs is established. This is approached by considering the joint dynamics of the corresponding nonlinear Fokker-Planck equation governing the law of the system and the underlying stochastic differential equation (SDE) as a dynamical system on the product space $\RR^d \times \mathcal{P}(\RR^d)$. 

    The proof relies on two main ingredients: At the level of the SDE, a pathwise rough path-based solution theory for SDEs with time-dependent coefficients is implemented, while at the level of the PDE a well-posedness theory is developed, for measurable solutions and allowing for degenerate diffusion coefficients. 

    The results apply in particular to the so-called ensemble Kalman sampler (EKS), proving the existence of an associated RDS under some assumptions on the posterior, as well as to the Lagrangian formulation of the Landau equation with Maxwell molecules. As a by-product of the main results, the uniqueness of solutions non-linear Fokker--Planck equations associated to the EKS is shown.
\end{abstract}
\setcounter{tocdepth}{2}
\tableofcontents
\section{Introduction}
In this paper, we construct a random dynamical system (RDS) for the following coupled system of a distribution dependent SDE, and its nonlinear Fokker-Planck equation
\begin{subequations}\label{system0}
\begin{align}
    \d Y_t &=  b(Y_t,\mu_t)\, \d t+\sigma(Y_t, \mu_t)\,\d W_t,\;\;\;\;\;Y_0=y\,,\label{SDE0}\\
    \partial_t \mu_t &=\frac{1}{2}D^2:\left((\sigma\sigma^*)(\cdot,\mu_t)\mu_t\right) -\nabla\cdot(b(\cdot,\mu_t)\mu_t)\, ,\;\;\;\;\;\mu_0=\mu \, ,\label{PDE0}
\end{align}
\end{subequations}
where $b:\RR^d \times \cP(\RR^d)\to \R^d, \sigma:\RR^d \times \cP(\RR^d)\to \R^{d\times d}$ are the drift and diffusion coefficients, $\cP(\R^d)$ is the space of Borel probability measures, and $W_t$ is a $d$-dimensional Brownian motion. Setting $\mu_0=\delta_y$, recovers a system of McKean--Vlasov-type with $\mu_t$ being the law of $Y_t$. 

The distribution dependent SDE~\eqref{SDE0} describes the mean-field limit of weakly interacting diffusions (see \cite{Sznitman.1991.165}). This is a fundamental class of interacting particle systems appearing in several applications, ranging from the theory of random matrices~\cite{LLX20} and plasma physics~\cite{Wang.2018.SPA595} to machine learning and global optimisation algorithms~\cite{Garbuno-InigoHoffmannLiStuart.2020.SJoADS412,CFV22,CJLZ21} and biological models of chemotaxis~\cite{S00}. These systems have been studied from a wide variety of perspectives, for instance, well-posedness, propagation of chaos, ergodicity, phase transitions, etc. We refer the reader to~\cite{CarrilloGvalaniPavliotisEtAl.2020.ARMA635,CD22} and the references therein for more details. 

Despite their fundamental nature, the dynamical behaviour of the trajectories of~\eqref{SDE0} has received only limited attention. This article takes the first step in remedying this by establishing the existence of a  RDS associated to the system comprised of~\eqref{SDE0} and~\eqref{PDE0}. 

The theory of random dynamical systems allows the application of methods and concepts from the theory of dynamical systems in the context of stochastic processes, such as solutions to SDEs. Notably, this includes the use of the multiplicative ergodic theorem,  leading to the existence of Lyapunov spectra and to the construction of local, invariant manifolds, notions of bifurcation, and synchronisation/chaos~\cite{Arnold.1998.586}. However, in order to make this dynamics toolbox accessible, it has to be shown that the solutions to an SDE form a random dynamical system. Since the usual constructions in stochastic analysis do not control the occurrence of $\mathbb{P}$-zero sets, the standard strong solution theory which would only provide a version of this property that holds true $\mathbb{P}$-a.s.~with the null set depending on $s,t$ and the initial datum. The usual approach for SDEs as delineated in~\cite[Chapter 2]{Arnold.1998.586} involves first proving a `crude' version of the property~\eqref{eq:perfectcocyle} and then arguing that one can apply perfection theorems (see~\cite[Theorems 3 and 4]{KagerScheutzow.1997.EJP8} or~\cite[Theorem 28]{AS95}) to improve it to~\eqref{eq:perfectcocyle}. This approach fundamentally relies on the selection of a version of the solution that is  $\mathbb{P}$-a.s.~continuous in $s,t$ and the initial datum, which itself fundamentally relies on Kolmogorov's continuity theorem. No analog of this argument is known in infinite dimensional stochastic systems. However, due to the distribution dependence of~\eqref{SDE0}, the construction of an RDS for a distribution dependent SDE is inherently infinite dimensional and, thus, cannot proceed via the established route. This challenge is resolved in the present work by instead developing a path-by-path rough path approach. 

A second obstacle met in the construction is the (ir)regularity of solutions to the nonlinear Fokker--Planck equation~\eqref{PDE0}. Key examples of distribution dependent SDEs, like the ensemble Kalman sampler (see~\cref{ssub:the_ensemble_kalman_sampler}) lead to SDEs~\eqref{SDE0} with potentially degenerate diffusivities $(\sigma \sigma^*)(\cdot,\cdot)$. Therefore, we do not rely on any smoothing effect from the second-order term in ~\eqref{PDE0}, and have to work with merely measure-valued solutions. This obstacle is resolved in this work by developing a duality argument which results in a stability estimate for the PDE~\eqref{PDE0} in a  one-parameter family of metrics   on the space of Borel probability measures. 

The application to the ensemble Kalman sampler (see~\eqref{eks0}), for example, in the setting sampling from a Gaussian, leads to linearly growing drift in $x$. As a consequence, solutions to the dual Kolmogorov equation exhibit polynomial growth. This  necessitates working with probability measures that possess at least finite second moments and requires the introduction of a metric (see~\eqref{metric-measure}) that accommodates test functions exhibiting at least quadratic growth.

The choice of these metrics requires case. They are chosen to be exactly compatible with the growth and regularity assumptions on the coefficients and thus allow to close a Gr\"{o}nwall-type estimate and obtain the required stability. This would not be possible using other more well-known metrics on probability measures, like the Wasserstein metric or the total variation metric.  Nevertheless, it turns out that each metric in this family generates the same topology as an associated $p$-Wasserstein metric, thus leading to a solution map which is continuous in the $p$-Wasserstein metric.

A related challenge arises at the level of the RDE. Since we aim to treat examples coming from sampling in which the drift can grow linearly  at infinity, this requires the development of a rough path theory for RDEs with unbounded drifts and time-inhomogeneous coefficients. We do this by using a time-dependent version of the so-called Doss--Sussmann transformation (see~\cite{D22} for a similar approach).

We state an informal version of the main result about the existence of random dynamical system below. The mathematically precise version of the above result can be found in \cref{meas,exisrds}.
\begin{theorem}\label{gener-RDS}
Given a probability space $(\Omega,\mathcal{F},\mathbb{P})$ and family $(\theta_t)_{t\geq 0}:\Omega \to \Omega$ of measure-preserving shifts, under appropriate assumptions on coefficients $b,\sigma$ (see~ \cref{assumPDE}, \cref{ass:RDE}, 
\cref{assume3}), there exists an $p\geq2$  and a continuous-in-time random dynamical system $\varphi:[0,\infty)\times\Omega \times \R^d \times \cP_p(\R^d) \to \R^d \times \cP_p(\R^d)$ associated with the coupled equations~\eqref{SDE0} and~\eqref{PDE0}.
That is, $\varphi\left(t, \omega, (y,\mu)\right)$ is a solution of the system~\eqref{SDE0} and~\eqref{PDE0} with initial datum $(y,\mu)$ and, more importantly, it is a ``perfect cocycle'', i.e.\ for any $s, t\geq 0$, $\omega\in\Omega$, $(y,\mu)\in\RR^d\times \ca{P}_{p}(\RR^d)$,
\begin{equation}
\varphi\left(t+s, \omega, (y,\mu)\right)
=\varphi(t, \theta_s\omega,\varphi(s, \omega,(y,\mu))) \, .
\label{eq:perfectcocyle}
\end{equation}
\end{theorem}
The proof proceeds by first developing a well-posedness and stability theory for the nonlinear Fokker--Planck equation~\eqref{PDE0}. Subsequently, its solution is fed into the~\eqref{SDE0} leaving us with an SDE with multiplicative noise and time-dependent coefficients for which we develop a pathwise solution theory based on the theory of controlled rough paths~\cite{Gub04}. This allows us to obtain a solution map which is continuous with respect to time, the initial datum, and the (suitably enhanced) noise. Thus, instead of using a perfection theorem, we can prove the perfect cocycle property~\eqref{eq:perfectcocyle} by hand.

\begin{remark}\label{rem:fromrdstomv}
 Note that given an RDS $\varphi$ for the system described by~\eqref{SDE0} and~\eqref{PDE0}, one can recover the dynamics and the cocycle property at the level of the McKean--Vlasov system by setting $\mu=\delta_{y}$. More importantly, we are required to treat the joint system~\eqref{system0} in which the solution of~\eqref{PDE0} is not necessarily the law of~\eqref{SDE0}. Indeed, if we treated~\eqref{SDE0} by itself it would be impossible to construct an RDS since the evolution of $X_t$ depends on $\mu_t$. On the other hand, if we insisted on choosing initial data of the form $(y,\delta_{y}) \in \R^d \times \mathcal{P}_p(\R^d)$, then we would enforce the constraint $\mu_t :=\ca{L}_{Y_t}=\mathrm{Law}(Y_t)$. However, there would be no reasonable state space on which one can define the RDS since even if the system~\eqref{system0} starts on the set $(y ,\delta_{y})_{y \in \R^d}$ it would leave this set immediately. 
\end{remark}

\subsection{Applications}
In this section, we present particular systems to which this work's main theorem applies. 
\subsubsection{The ensemble Kalman sampler} 
\label{ssub:the_ensemble_kalman_sampler}
The results of~\cref{gener-RDS} apply in particular to the so-called ensemble Kalman sampler which takes the following form:
\begin{align}\label{eks0}
    \d Y_t=-{\rm{Cov}}(Y_t)\nabla V(Y_t)\,\d t+\sqrt{2{\rm{Cov}}(Y_t)}\,\d W_t \, ,
\end{align}
where ${\rm{Cov}}(Y_t)$ is the covariance matrix of the random variable $Y_t$ and $W_t$ is a $d$-dimensional Brownian motion. This system was introduced in~\cite{Garbuno-InigoHoffmannLiStuart.2020.SJoADS412} (see also~\cite{NR19}) as a potential method of sampling from the posterior $\rho \sim e^{-V}$. The system~\eqref{eks0} is expected to exhibit better convergence properties as $t\to \infty$ as compared to the standard overdamped Langevin dynamics (obtained by replacing ${\rm Cov}(Y_t)$ with $\mathrm{Id}$ in~\eqref{eks0}). This is due to the fact that the presence of the covariance term ${\rm Cov}(Y_t)$ forces the system to `sample faster' in directions along which the variance of the target distribution is larger. This is most easily demonstrated in the setting $\rho=\mathcal{N}(0,\Sigma)$. It was shown in~\cite{BEHMS25} that the law $\rho_t$ of~\eqref{eks0} converges exponentially fast to $\rho$ with rate independent of $\Sigma$, while the law of the overdamped Langevin dynamics converges exponentially but with rate dependent on the condition number of $\Sigma$. Clearly, the first scenario is preferable, especially for $\Sigma$ which are strongly anisotropic. 
Additionally,~\eqref{eks0} can also be efficiently approximated by an interacting particle system. Consider the following system of $N$ interacting particles:
\begin{equation}
 \d Y_t^i=-{\rm Cov}(\mu^N_t)\nabla V(Y_t^i)\, \d t+\sqrt{{\rm Cov}(\mu^N_t)}\,\d W_t^i,\;\;\;Y_0^i\sim \rho_0 \, ,
\end{equation}
where the $Y_t^i, i=1,...N,$ represent the positions of the $N$ particles, the $W_t^i, i=1,...N$ represent $N$ independent $d$-dimensional Brownian motions, and 
\begin{align}
    \mu^N_t:=\frac{1}{N}\sum_{i=1}^N \delta_{Y_t^i}\, .
\end{align}
Then, it was conjectured in~\cite{Garbuno-InigoHoffmannLiStuart.2020.SJoADS412} and shown in~\cite{GHV24} that, under appropriate assumptions on $V,\rho_0$, $Y_t^i$ converges in law to the solution of~\eqref{eks0}. For our purposes, we treat~\eqref{eks0} by decoupling the evolution of the law from $Y_t$ and thus bringing it into the form of~\eqref{system0}, i.e.\
\begin{subequations}\label{eksdecoupled}
\begin{align}
    \d Y_t= & \, -{\rm{Cov}}(\mu_t)\nabla V(Y_t)\,\d t+\sqrt{2{\rm{Cov}}(\mu_t)}\,\d W_t\, , \quad Y_0=y \, , \label{eksdecoupledSDE}\\
    \partial_t \mu_t =&\,D^2:\left({\rm Cov}(\mu_t)\mu_t\right) -\nabla\cdot({\rm Cov}(\mu_t) \nabla V \mu_t)\, ,\;\;\;\;\;\mu_0=\mu \, \label{eksdecoupledPDE}.
\end{align}
\end{subequations}
Here~\eqref{eksdecoupledSDE} is the same as~\eqref{eks0} except that it is driven by the covariance of some measure $\mu_t$ which itself is a solution of the PDE~\eqref{eksdecoupledPDE}. Setting $\mu=\delta_{y}$, one recovers~\eqref{eks0}. Applying~\cref{gener-RDS} to~\eqref{eksdecoupled}, we obtain the following result.
\begin{corollary} [\cref{eks}]
Assume $V$ satisfies certain regularity and growth conditions and that $(y,\mu)\in \R^d\times \cP_p(\R^d)$ for some $p\geq 2$ depending on the choice of $V$. Then, there exists an RDS $\varphi:[0,\infty)\times\Omega \times \R^d \times \cP_p(\R^d) \to \R^d \times \cP_p(\R^d)$ associated to~\eqref{eksdecoupled}.
\end{corollary}
The above result opens the door to studying dynamical aspects of the ensemble Kalman sampling algorithm. A natural next question is whether the algorithm exhibits the phenomenon of synchronisation-by-noise, which is a desirable property if one wishes to construct higher-order sampling schemes (see~\cite{LPP15}).

\subsubsection{Lagrangian formulation of the Landau equation} 
 \label{sub:lagrangian_formulation_of_the_landau_equation}
Consider the spatially homogeneous Landau equation 
\begin{equation}\label{eq:Landau}
    \partial_t f_t=\frac{1}{2}\nabla\cdot\left(\int_{\RR^d}a(\cdot-z)\left(f_t(z)\nabla f_t-f_t\nabla f_t(z)\right)\d z\right),
\end{equation}
for some collision kernel $a:\RR^d\to \RR^{d\times d}$, which describes the distribution of charged particles in a collisional plasma in the grazing collision limit. In particular, if $d=3$, the collision kernel is chosen to be
\begin{align}
    a(x)=|x|^{2+\gamma}\left(I-\frac{x\otimes x}{|x|^2}\right),\,\,\,\,\,x\in\RR^d \,
\end{align} 
for some constant $\gamma \in [-3,1]$. The choice $\gamma=-3$ is the most physically interesting choice as it corresponds to Coulomb interactions between the charged particles. The case $\gamma=0$ corresponds to Maxwell molecules and $\gamma\in (0,1]$ corresponds to hard potentials. The Landau equation can be interpreted as the nonlinear Fokker-Planck equation corresponding to the following distribution-dependent SDE 
\begin{align}
    \d Y_t=(b_0* \mathrm{Law}(Y_t))(Y_t)\,\d t+(\sigma_0* \mathrm{Law}(Y_t))(Y_t)\,\d W_t\, .
\end{align}
for the choice $b_0=\nabla \cdot a$ and $\sigma_0$ such that $\sigma_0\sigma_0^*=a$, where for a measurable function $g$ and a probability measure $\mu$, $(g\ast \mu)(x)\coloneqq \int_{\RR^d}g(x-z)\,\d \mu(z)$. It was shown in~\cite[Theorem 6.1]{Wang.2018.SPA595} that the above SDE has a unique strong solution  for $\gamma \in [0,1]$. Since $\mathrm{Law}(Y_t)=f_t$, the above SDE provides a Lagrangian interpretation of the Landau equation~\eqref{eq:Landau}. Note that if $\gamma=0$, then
\begin{align}
    b_0(y)=-2y\, , \quad \sigma_0(y)= & \
        \begin{pmatrix}
            y_2 & 0 &y_3\\
            -y_1 & y_3 & 0 \\
            0 & -y_2 & -y_1
        \end{pmatrix}\, .
\end{align}
We then apply our general result~\cref{exisrds} to the following system which we can obtain by setting $\gamma=0$ for the Landau system.
\begin{subequations}\label{eq:landaudecoupled}
        \begin{align}
        \d Y_t=&\,(-2Y_t+2m(\mu_t))\,\d t+(\sigma_0(Y_t)-\sigma_0(\mu_t))\,\d W_t,\,Y_0=y\label{eq:Landau1}\, ,\\
    \partial_t \mu_t=&\,  -\nabla \cdot \big((-2\cdot+2m(\mu_t))\mu_t\big)+\frac{1}{2}D^2: \big[(\sigma_0(\cdot)-\sigma_0(\mu_t))(\sigma_0(\cdot)-\sigma_0(\mu_t))^*\mu_t\big],\,\mu_0=\mu\label{eq:Landau2}\,, 
        \end{align}
    \end{subequations}
    where
    \begin{align}
\sigma_0(\mu)\coloneqq  \int_{\R^3}\sigma_0(y)\, \ \mathrm{d}\mu(y)\,, \;\;\;\;\;
        m(\mu)\coloneqq \int_{\RR^3}y\,\d\mu \label{eq:mean}\, .
    \end{align}
Since the well-posedness theory presented in this article is based on controlled rough paths which requires relatively strong growth and regularity assumptions on $b,\sigma$, we  are forced to restrict ourselves to the setting $\gamma=0$. We have the following result.
\begin{corollary} [\cref{landauprop}]
Assume $\gamma=0$, i.e.\ Maxwell molecules,  $(y,\mu)\in \R^d\times \cP_p(\R^d)$ for some $p\geq 2$. Then, there exists an RDS $\varphi:[0,\infty)\times\Omega \times \R^d \times \cP_p(\R^d) \to \R^d \times \cP_p(\R^d)$ associated to~\eqref{eq:Landau1} and~\eqref{eq:Landau2}.
\end{corollary}

\subsection{Overview on the existing literature} 
The generation of RDS by finite-dimensional SDEs has been thoroughly understood, and we refer to \cite[Theorem 2.3.26, Theorem 2.3.32]{Arnold.1998.586} and the references therein for the details.

The theory of rough paths factorizes the solution of a stochastic differential equation into a step of lifting the driving stochastic process into a rough paths, neglecting $\mathbb{P}$-zero sets as needed, and a second, purely analytic step arguing path-by-path to construct solutions to the equation. As a result, $\mathbb{P}$-zero sets appearing in this approach are agnostic to the coefficients and initial data of the equation, and, thus no further perfection arguments are required. This approach has been put to use in \cite{Bailleul2017} in order to construct RDS by means of rough paths theory, a key remaining step being the introduction of the concept of rough cocycles, and the validation that (certain) rough paths lifts are rough cocylcles. 

McKean-Vlasov SDEs driven by rough paths have received considerable interest in recent years. Initiated by \cite{CassLyons.2015.PLMS383} for the case of drift depending linearly on the interaction term, subsequently extended in \cite{Bailleul.2015.RMI901} to drifts allowing nonlinear interaction, and by \cite{Bailleul2020,Bailleul2021,Coghi2021} to allowing interaction dependency also in the diffusion coefficients. A partial extension to include also unbounded coefficients has been given in  \cite{Coghi2023} in the context of Ensemble Kalman filtering. 
    
\subsection{Organisation of the paper}

In \cref{setup}, we introduce the notation and preliminary concepts needed for the rest of paper. In particular, this includes the required prerequisites from the theory of controlled rough paths along with the relevant notions of solution for~\eqref{PDE0} and~\eqref{SDE0}.  In \cref{sec-well}, we introduce the assumptions we work with and present the precise versions of the main results of this paper, including the well-posedness theory for~\eqref{PDE0} and~\eqref{SDE0} along with the theorem on the existence of an RDS. \cref{sec:proofs} contains the proof of all the main results. Finally, in~\cref{sec-appli}, we verify all the assumptions on coefficients for the ensemble Kalman system~\eqref{eksdecoupled} and Landau system~\eqref{eq:landaudecoupled}. We collect auxilliary results that are more technical in ~\cref{sec:proofconvergenceequiva,Moment estimates,Differentiability of backward PDE,sec:calculations}.

\section{Set up and preliminaries}\label{setup}
Let $\bb{N}$ be the set of natural numbers including 0. We denote by $\alpha=(\alpha_1,\ldots,\alpha_d)\in \bb{\bb{N}}^{d}$  a multi-index with size  $|\bar{\alpha}|\coloneqq\alpha_1+\ldots+\alpha_d$ and for a function $f$ defined on $[0,T]\times \RR^d$ or $[0,\infty) \times \R^d$, we define
$$\partial^{\alpha} f\coloneqq\left(\frac{\partial}{\partial y_1}\right)^{\alpha_1}\left(\frac{\partial}{\partial y_2}\right)^{\alpha_2}\ldots\left(\frac{\partial}{\partial y_d}\right)^{\alpha_d} f \, ,$$
where the derivatives in $t$ at $t=0,T$ are understood in a one-sided sense. For any continuously differentiable functions $f:[0,T]\times \RR^d\to \R$, $m,n\in \mathbb{N}$, we define the seminorms
\begin{equation}
     \label{eq:seminorms}
     [f]_{m,n} =  \sum_{\substack{\bar{\alpha}\in \bb{N}^{d},\\
\alpha_0=m, |\bar{\alpha}|= n}} \sup\limits_{y}\sup\limits_{ t} |\partial^\alpha f|(t, y) \, ,
 \end{equation} 
 and the associated norm
 \begin{equation}
     \| f \|_{m,n}=\sum_{0 \leq k \leq m, 0\leq \ell \leq n} [f]_{k,\ell} \, .
 \end{equation}
In particular, if $f: \RR^d \to \RR$  we define the norm
 $$\|f\|_n:=\sum_{\substack{\bar{\alpha}\in \bb{N}^{d}\,,\\
 0\leq |\bar{\alpha}|\leq  n}} \sup\limits_{y} |\partial^\alpha f|(y)\,.$$
 We will also need the more general seminorms $[f]_{\beta,n}$, for  $\beta \in (0,1), n \in \mathbb{N}$, which we define as follows
\begin{equation}
    \label{eq:Holderseminorm}
    [f]_{\beta,n} := \sum_{\substack{\bar{\alpha}\in \bb{N}^{d}\,,\\
\alpha_0=0, |\bar{\alpha}|= n}} \sup\limits_{y}\sup\limits_{ t\neq s} \frac{|\partial^\alpha f(t, y) - \partial^\alpha f(s, y)|}{|t-s|^\beta}\, .
\end{equation}
We denote by $[f]_{\mathrm{Lip},n}$, the above norm when $\beta=1$ and 
\begin{align}
	[f]_{0, \rm{Lip}}:=\sup_{t\in[0,T]}\sup_{x\neq y\in\RR^d}\frac{|f(t, x)-f(t, y)|}{|x-y|}\,. 
	\end{align}
Note that whenever we consider an $f$ valued in $\R^m$, then we will use the same notations for the seminorms/norms with the understanding that we take the maximum over all components. 

We now introduce the required prerequisites from the theory of random dynamical systems. Fix a metric space $(E, d)$  with Borel $\sigma$-algebra $\ca{E}$ and consider a metric dynamical system $(\Omega, \ca{F},\bb{P}, (\theta_t)_{t \in \RR})$, i.e.\ $(\Omega, \ca{F},\bb{P})$ is a probability space and $\theta_t:\Omega\to \Omega, t\in \RR$ is a group of mappings satisfying: 
\begin{enumerate}[(i)]
\item $(\omega,t)\mapsto \theta_t\omega$ is $(\ca{B}(\RR)\otimes\ca{F}, \ca{F})$-measurable;
\item (Flow property) $\theta_0={\rm{Id}}_{\Omega}$ and $\theta_{t+s}=\theta_t
\circ\theta_s$;
\item (Measure preserving) For each $t \in\RR$, $\theta_t\PP=\PP$.
\end{enumerate}
\begin{definition}\label{def:rds}
A random dynamical system $\varphi$ (see \cite[Definition 1.1.1]{Arnold.1998.586}) on $(E,\ca{E})$ over a metric dynamical system $(\Omega,\ca{F},\bb{P},(\theta_t)_{t \geq 0})$ is a mapping 
\begin{align}
\varphi:[0, \infty)\times \Omega\times E\to E\,,\;\;\;\;\;\;\;\;\;(t,\omega,e)\mapsto \varphi(t,\omega)e\,,
\end{align}
such that
\begin{enumerate}[(i)]
\item 
$\varphi$ is $(\mathcal{B}([0, \infty))\otimes \mathcal{F}\otimes \mathcal{E} ; \mathcal{E})$-measurable; 

\item The mappings $\varphi(t,\omega):=\varphi(t,\omega,\cdot)\colon E\to E$ form a cocycle over $\theta$, i.e. for all $\omega\in\Omega$, for all $0 \leq s ,t <\infty $, we have  $\varphi(0,\omega)={\rm{Id}}_{E}$ and 
\begin{align}\label{cocy}
\varphi(t+s,\omega)=\varphi(t,\theta_s\omega)\circ \varphi(s,\omega) \, . \;\;\;\;\;\;\;\;\;\;\;\;\;\;\;\;\;\;\text{(Cocycle \;property) }
\end{align} 
\end{enumerate}
\end{definition}

We will now describe the choice of metric space $E$ for our setting. To this end, we denote by $\mathcal{P}(\bb{R}^d)$ the space of all Borel probability measures on $\RR^d$ and, for any $p\geq 1$, by $\mathcal{P}_{p}(\bb{R}^d)$ the subspace of $\mathcal{P}(\R^d)$ of measures with finite $p$-th moment, i.e.\
$$\mathcal{P}_{p}(\bb{R}^d)\coloneqq\Big\{\mu\in \ca{P}(\bb{R}^d), \;M_{p}(\mu):=\int_{\bb{R}^d}|y|^{p}\,\,\d \mu(y) <\infty\Big\} \, .$$
We define the following class of functions
$$\ca{K}_p=\big\{\varphi\in C^2(\RR^d): |\nabla \varphi|(y)\leq 1+|y|^{p-1},\,
|D^2\varphi|(y)\leq 1+|y|^{p-1}\,\big\}\,.$$ 
\begin{lemma}\label{convergenceequiva}
For any $\mu,\nu\in\ca{P}_p(\RR^d)$, $p\geq 1$, define 
\begin{align}
\mathcal{D}_p(\mu,\nu)&\coloneqq\sup\limits_{\varphi\in \ca{K}_p}
\int_{\RR^d}\varphi\, \d (\mu-\nu)\,,\;\;\;\;\;\;\;\;\mu,\nu\in\ca{P}_p(\RR^d)\,,\label{metric-measure}
\end{align}
and consider the $p$-Wasserstein distance
\begin{align} 
d_p(\mu,\nu)\coloneqq\inf\limits_{\pi\in\mathcal{C}(\mu,\nu)}\Big(\int_{\RR^d\times \RR^d}|x-y|^p\;\pi(\d x \d y)\Big)^{\frac{1}{p}},
\end{align}
with $\mathcal{C}(\mu,\nu)$ being the set of all couplings of $\mu$ and $\nu$. 
Then, the following statements hold true:
\begin{enumerate}[(i)]

\item $\ca{D}_p$ is a  complete, separable metric on $\cP_p(\R^d)$;

\item $\ca{D}_p$ metrises the same topology as $d_p$ on $\ca{P}_p(\R^d)$.

\end{enumerate}

\end{lemma} 
 We relegate the proof of the above lemma to~\cref{sec:proofconvergenceequiva}. We now equip the space $\ca{P}_p(\R^d)$ with the metric $\ca{D}_p$ and  define $\uE_p:=\R^d \times \ca{P}_{p}(\R^d)$. We will denote a typical element of $\uE_p$ as $\ue= (y, \mu)$ and will use $\ue_1 = y \in \RR^d$ and $\ue_2 =\mu \in \mathcal{P}_p(\RR^d)$ to refer to its components.
 We equip $\uE_p$  with the product topology, i.e.\ the one  generated by the following metric:
\begin{equation}
       \ud_p(\ue,\uf)=|\ue_1-\uf_1|+\ca{D}_p(\ue_2,\uf_2) \, ,
\end{equation}
and  its Borel $\sigma$-algebra which is exactly the product $\sigma$-algebra $\ca{E}_p=\ca{B}(\RR^d)\otimes \ca{B}(\ca{P}_p(\RR^d))$.
Before we can describe the  choice of the metric dynamical system  $(\Omega,\mathcal{F},\mathbb{P},(\theta_{t})_{t\in\RR})$ in this paper, we need to introduce some concepts from the theory of rough paths which we borrow from~\cite{FrizHairer.2014.251}.


A rough path on an interval $[0,T]$ with values in $\R^d$ consists of a continuous function $X:[0,T]\to \R^d$ and  an enhanced continuous ``second order process'' $\mathbb{X}:[0,T]^2\to \R^{d\times d}$, subject to certain algebraic and analytical conditions.  Precisely speaking,  for $\alpha\in(\frac{1}{3},\frac{1}{2}]$, we define the space $\scr{C}^{\alpha}([0,T];\R^d)$ of $\alpha$-H${\rm{\ddot{o}}}$lder rough paths over $\R^d$ as those pairs $\mathbf{X}\coloneqq(X,\mathbb{X})$ such that \begin{align}
&\|X\|_{\alpha}\coloneqq\sup_{s\neq t\in[0,T]}\frac{| X_{s,t}|}{|t-s|^{\alpha}}<\infty\,,\,\,\,\,\,\,\,\,\,\,\,\,\,\,\,\|\mathbb{X}\|_{2\alpha}\coloneqq\sup_{s\neq t\in[0,T]}\frac{|\mathbb{X}_{s,t}|}{|t-s|^{2\alpha}}<\infty\,,
\end{align}
where $ X_{s,t}=X_t-X_s$, and such that Chen's relation
$$\mathbb{X}_{s,t}-\mathbb{X}_{s,u}-\mathbb{X}_{u,t}=X_{s,u}\otimes X_{u,t}$$
is satisfied (see \cite[Definition 2.1]{FrizHairer.2014.251}). We also introduce the $\alpha$-H$\rm{\ddot{o}}$lder homogeneous rough path norm 
$$\|\bfX\|_{\alpha}\coloneqq\|X\|_{\alpha}+\sqrt{\|\bb{X}\|_{2\alpha}}\,.$$   The most fundamental example and the one that is relevant for this paper is that of $d$-dimensional standard Brownian motion $W$ enhanced with
$$\bb{W}_{s,t}\coloneqq\int_s^tW_{s,r}\otimes \d W_r\, .$$ 
We will use $\mathbf{W}$ to denote the pair $(W,\mathbb{W})$. We now introduce the concept of controlled rough paths. Given a path $X\in C^{\alpha}([0,T];\RR^d)$, a path $Y \in C^\alpha([0,T]; \RR^d)$ is said to be controlled by $X$ if there exists $Y'\in C^{\alpha}([0,T];\RR^{d \times d})$, referred to as the Gubinelli derivative, so that the remainder term $R^Y$ given by 
\begin{align}\label{re}
Y_{s,t}=Y_s'X_{s,t}+R_{s,t}^Y\,,
\end{align}
satisfies $\|R^Y\|_{2\alpha}<\infty$. This defines the space of controlled rough paths $(Y, Y')\in \scr{D}_X^{2\alpha}([0,T]; \R^d)$ and we endow the space $\scr{D}_X^{2\alpha}$ with the seminorm 
\begin{align}
\|Y,Y'\|_{X,2\alpha}\coloneqq\|Y'\|_{\alpha}+\|R^Y\|_{2\alpha}\,.
\end{align}
We could also realise $\scr{D}_X^{2\alpha}$ as a Banach space if we equip it with the norm $(Y, Y')\mapsto |Y_0|+|Y_0'|+\|Y,Y'\|_{X,2\alpha}$. With these objects at hand, we are able to define the rough integral of controlled rough path $Y$ against a rough path $\bfX$ as
$$\int_0^T Y_s \,\d\bfX_s\coloneqq\lim\limits_{|\ca{P}|\to 0}\sum_{[u,v]\in\ca{P}}(Y_uX_{u,v}+Y_u'{\bb{X}}_{u,v})\,,$$
with $\ca{P}$ being a partition over $[0,T]$. Moreover, it follows from \cite[Theorem 4.10]{FrizHairer.2014.251} that one has the bound
\begin{align}\label{es:roughintegral-remainder}
	\Big|\int_{s}^t Y_r \d \bfX_r-Y_sX_{s,t}-Y_s'{\bb{X}}_{s,t} \Big|\leq\, C\big(\|X\|_{
\alpha}\|R^Y\|_{2\alpha}+\|\bb{X}\|_{2\alpha}\|Y'\|_{\alpha} \big)|t-s|^{3\alpha}\,.
\end{align}

Given two rough paths $\bfX,\tilde{\bfX}\in \scr{C}^{\alpha}([0,T];\RR^d)$, the $\alpha$-H$\rm{\ddot{o}}$lder rough path metric is defined as follows
$$\varrho_{\alpha}(\bfX,\tilde{\bfX})\coloneqq\sup\limits_{s\neq t\in[0,T]}\frac{|X_{s,t}-\tilde{X}_{s,t}|}{|t-s|^{\alpha}}+\sup\limits_{s\neq t\in[0,T]}\frac{|\bb{X}_{s,t}-\bb{\tilde{X}}_{s,t}|}{|t-s|^{2\alpha}}\,.$$
Suppose $(Y,Y')\in \scr{D}_{X}^{2\alpha}$, $(\tilde{Y},\tilde{Y}')\in \scr{D}_{\tilde{X}}^{2\alpha}$, we also define the ``distance''
$$\|Y,Y'; \tilde{Y}, \tilde{Y}'\|_{X,\tilde{X},2\alpha}\coloneqq\|Y'-\tilde{Y}'\|_{\alpha}+\|R^Y-R^{\tilde{Y}}\|_{2\alpha}\,.$$

Given all these notions, we are finally in a position to define the metric dynamical system $(\Omega,\mathcal{F},\mathbb{P},(\theta_{t})_{t\in\RR})$. We start by setting $\bar{\Omega}\coloneqq C_0(\RR;\R^d)$ the space of continuous $\R^d$-valued paths $t \mapsto \omega(t)$ which are $0$ at $t=0$ and equip it with compact-open topology. We set $\bar{\ca{F}}$ to be the associated Borel $\sigma$-algebra and $\bar{\mathbb{P}}$ to be the Wiener measure. We now define the $(\theta_t)_{t\in\RR}:\bar{\Omega}\to \bar{\Omega}$ as the natural left shift
\begin{equation}
(\theta_t \omega)(\cdot):=\omega(t+\cdot)-\omega(t) \, .
\end{equation}
It follows from~\cite[Appendix A.3]{Arnold.1998.586} that $\theta$ is $(\ca{B}(\RR\otimes\bar{\ca{F}},\bar{\ca{F}})$-measurable and leaves $\bar{\mathbb{P}}$ invariant and so the tuple $(\bar{\Omega},\bar{\mathcal{F}},\bar{\mathbb{P}},(\theta_{t})_{t\in\RR})$ is a metric dynamical system. This is not enough for our purposes, we would like to enhance elements of $\bar{\Omega}$ with additional second-order information since we will be dealing with rough paths. We denote by 
\begin{equation}
W_t(\omega)\coloneqq\omega(t) \, ,
\end{equation}
the Brownian motion associated to each $\omega \in \bar{\Omega}$. We also introduce the following dyadic approximation $W^n_t$ which is a piecewise linear and continuous function with value exactly $W_t$ for all $t= k2^{-n}$, $k=\dots,-n,\dots,-1,0, 1,\dots, \infty.$ We note that, by~\cite[Propositions 3.5 and 3.6]{FrizHairer.2014.251}, for all $m=1,2,\dots$ the following set
\begin{align}
&\Omega_m:= \bigg\{\omega\in \bar{\Omega}:(W(\omega),\mathbb{W}^{\,\rm Strat.}(\omega))\in \scr{C}^\alpha([-m,m];\R^d) \text{ is well-defined}, \quad  \\ 
&\qquad\qquad\quad\,\Big(W^n(\omega),\int_{-m}^{\cdot} W^n(\omega) \otimes \d W^n(\omega)\Big) \xrightarrow{n \to \infty } (W,\mathbb{W}^{\,\rm Strat.})\in \scr{C}^\alpha([-m,m];\R^d),\\
&\qquad\qquad\quad\text{ for all }\alpha \in (1/3,1/2)  \bigg\} \, ,
\end{align}
has measure $1$, where $\mathbb{W}^{\, \rm Strat.}$ is the Stratonovich lift of $W$ and the integral is understood as Riemann--Stieltjes integral. We then define
\begin{equation}
\Omega= \bigcap_{m=1}^\infty \Omega_m \, ,
\end{equation}
and note that for all $\omega \in \Omega$, $(W,\mathbb{W}^{\,\rm Strat.}) \in \scr{C}^\alpha([-T,T];\R^d)$ is well-defined  and the dyadic approximation converges in $\scr{C}^\alpha([-T,T];\R^d)$ for all $T<\infty$ and $\alpha \in (1/3,1/2)$. Clearly, the lift $(W,\mathbb{W}^{\,\text{It\^{o}}})$ is also well-defined for all $T<\infty$. We now check that $\theta_t$ leaves $\Omega$ invariant  for all $t \in\RR$. To do this, we need to understand what action the shift has on the second-order information of the path $W$. To this end, we note that for all $r \in\RR$, $W(\theta_r \omega)$ can be equipped with the  second-order process $\mathbb{W}_{s+r,t+r}^{\, \rm Strat.} (\omega)$ and we have the identity (see~\cite[Lemma 3.2]{NeamctuKuehn.2021.SJMA3912})
\begin{equation}\label{eq:roughshift}
 \mathbb{W}_{s,t}^{\, \rm Strat.} (\theta_r \omega)=\mathbb{W}_{s+r,t+r}^{\, \rm Strat.} (\omega) \, ,
\end{equation} 
for all $s,t,r\in\RR$. The analogous construction and conclusion also holds for the It\^o lift $\mathbb{W}^{\,\text{It\^{o}}}$. Furthermore, one can readily check that $(W^n(\theta_r\omega),\int_{-T}^\cdot W^n(\theta_r\omega) \otimes \d W^n(\theta_r\omega))$ converges to $(W(\theta_r \omega),\mathbb{W}(\theta_r\omega))$ in $\scr{C}^\alpha([-T,T];\R^d)$ for all $T < \infty,r\in\RR$ and $\alpha \in (1/3,1/2)$.  Thus, for all $r\in\RR$, $\theta_r \Omega \subseteq \Omega$. We now define $\ca{F}=\bar{\ca{F}}\cap \Omega$ and note that $\theta$ must be $(\ca{B}(\RR\otimes\ca{F},\ca{F})$-measurable. If we denote by $\mathbb{P}$ the restriction of $\bar{\mathbb{P}}$  to $\Omega$, clearly $\theta_t$ leaves $\mathbb{P}$ invariant for all $t\in\RR$. This completes the construction of the metric dynamical system $(\Omega,\mathcal{F},\mathbb{P},(\theta_{t})_{t\in\RR})$.
\section{Main Results}\label{sec-well}
In this section, we collect the assumptions and background, along with the statements of the main results of this paper. More precisely, we prove the well-posedness of McKean--Vlasov PDEs in \cref{ssec:well-PDE}, the well-posedness of time-inhomogeneous   RDEs with drift in \cref{ssec:wellposednessRDE}, and, finally, the existence of an RDS is established in \cref{ssec:RDS}.
\subsection{Well-posedness of McKean-Vlasov PDEs}\label{ssec:well-PDE}
We consider the well-posedness of the  nonlinear Fokker--Planck equation
\begin{align}\label{MVPDE}
\partial_t \mu_t &=-\nabla\cdot(b(\cdot,\mu_t)\mu_t)+\frac{1}{2}D^2:\left((\sigma\sigma^*)(\cdot,\mu_t)\mu_t\right)\,.
\end{align}
We first introduce the assumptions we need to establish a well-posedness theory for McKean--Vlasov PDE~\eqref{MVPDE}.
\begin{assumption}\label{assumPDE}
Assume there exists $\kappa\geq 2$ and a  function $F:\RR_+^3\to \RR_+$, satisfying $F(x,y,z)\leq F(x,y_1,z_1), y\leq y_1,z\leq z_1$, such that the following conditions hold.

\begin{enumerate}[(i)]
    \item 
    %
    For any  $x,y\in \RR^d$, $\mu, \nu\in \ca{P}_{\kappa}(\RR^d)$,
    \begin{align}\label{localinstate}
        &| b(x,\mu)-b(y,\mu)|+|\sigma(x,\mu)-\sigma(y,\mu)|
        \leq F(\kappa, M_{\kappa}(\mu),0)|x-y|\,,
    \end{align}
    \worknote{This condition is to make sure the well-posedness of classical SDE~\eqref{claSDE}.}
    and
    \begin{align}\label{LocalLipschitzcontinuity}
        |b(y,\mu)-b(y,\nu)|+|\sigma\sigma^*(y,
        \mu)-\sigma\sigma^*(y,
        \nu)|\leq  F(\kappa, M_{\kappa}(\mu),M_{\kappa}(\nu))\big(1+|y|\big)\mathcal{D}_{\kappa}(\mu,\nu)\,,
    \end{align}
    where  $|\cdot|$ is understood as the Schatten-2 norm.  
    \item There exists a constant $C$ such that for any $y\in\RR^d$ and $\mu\in \ca{P}_{\kappa}(\RR^d)$,
    \begin{align}\label{Weakcoercivity}
        \langle b(y,\mu), y\rangle &\leq C\Big(1+|y|^2+M_{\kappa}^{\frac{2}{\kappa}}(\mu)\Big)\,,\\
        |\sigma(y,\mu)|&\leq C\Big(1+|y|+M_{\kappa}^{\frac{1}{\kappa}}(\mu)\Big)\,.\label{Growthsigma}
    \end{align}
    
    \item 
    For any $\mu\in\ca{P}_{\kappa}(\RR^d)$, $b(\cdot,\mu)$ and $\sigma(\cdot,\mu)$ are twice continuously differentiable with respect to $y$. Moreover, for any  $y\in\RR^d$ and $\mu\in \ca{P}_{\kappa}(\RR^d)$, 
    \begin{align}
        |b(0,\mu)|+|\nabla b(\cdot,\mu)|(y)+|\nabla \sigma(\cdot,\mu)|(y)&\leq  F(\kappa, M_{\kappa}(\mu),0)\,,\label{Growth1}\\
        |D^2 b(\cdot,\mu)|(y)+|D^2 \sigma(\cdot,\mu)|(y)&\leq  F(\kappa, M_{\kappa}(\mu),0)\,.\label{Growthseconderiva}
    \end{align}
\end{enumerate}

\end{assumption}
Given the assumptions on the coefficients, we now introduce the relevant notion of solution below.
\begin{definition}[Solution to the nonlinear Fokker--Planck equation]\label{defofPDE}
We say that a curve $[0,T]\ni t\mapsto \mu_t \in \ca{P}_p(\R^d)$, $p\geq 1$, is a measure-valued solution to~\eqref{MVPDE}  with initial condition $\mu_0 \in \cP_p(\R^d)$ if 
\begin{enumerate}[(i)]
    \item $(t, y)\mapsto b(y,\mu_t)$ and $(t, y)\mapsto \sigma(y,\mu_t)$ are continuous;
    \item for every $\phi\in \ca{K}_p$, and $t\in [0,T]$, the following identity holds true
\begin{align}
\begin{split}
&\int_{\RR^d}\phi(y)\, \d \mu_t(y)-\int_{\RR^d}\phi(y)\,\d\mu_0( y)\\
=\,\,&\int_0^t\int_{\RR^d}\Big(\langle b(y,\mu_s),\nabla\phi(y)\rangle+\frac{1}{2}(\sigma\sigma^{*})(y,\mu_s): D^2\phi(y)\Big)\, \d\mu_s( y)\,\d s \, .\label{eq:PDEsolu}
\end{split}
\end{align}
\end{enumerate}
In this case, we call that the nonlinear Fokker--Planck equation ~\eqref{MVPDE} is well-posed in $\ca{P}_{p}(\RR^d)$.
\end{definition}

\begin{remark}
It follows from the condition~\eqref{Growth1} that   $b$ has locally linear growth, that is, for any $y\in\R^d$, $\mu\in\ca{P}_{\kappa}(\R^d)$, 
\begin{equation}\label{Growthb}
	|b(y,\mu)|\leq |b(y,\mu)-b(0,\mu)|+|b(0,\mu)|\leq  F(\kappa, M_{\kappa}(\mu),0)(1+|y|)\,.
\end{equation}
Moreover, by condition~\eqref{LocalLipschitzcontinuity}, we have that \begin{align}
	|\sigma(y,\mu)-\sigma(y,\nu)|\leq \sqrt{F(\kappa, M_{\kappa}(\mu),M_{\kappa}(\nu))\big(1+|y|\big)\mathcal{D}_{\kappa}(\mu,\nu)}\,,
\end{align}
due to the fact that $|A^{\frac{1}{2}}-B^{\frac{1}{2}}|^2\leq \|A-B\|_{\rm{Tr}}$ (cf. \cite[Lemma 4.1]{PowersSto.1970.CMP1}) for matrices $A$ and $B$.
\end{remark}

\begin{theorem}\label{ThmWell-posednessPDE}
Under \cref{assumPDE},  for any initial distribution $\mu_0\in\ca{P}_{p}(\RR^d)$ with $p\geq \kappa$, there exists a solution $[0,T]\ni t\mapsto \mu_t \in \ca{P}_p(\R^d)$, in the sense of ~\cref{defofPDE} to the McKean--Vlasov PDE~\eqref{MVPDE}. Furthermore, the solution satisfies the following  estimates
\begin{subequations}
\begin{align}
 M_{p,T}(\mu)\coloneqq  &\,\,\,\sup_{t\in [0,T]}M_p (\mu_t) < \infty\,, \label{eq:strongclass}\\
\ca{D}_p(\mu_t,\mu_s) \leq& \,\,\, C |t-s| \label{eq:timeregmu} \, ,
\end{align} 
\end{subequations}
for all $t ,s \in [0,T]$ and for some constant $C=C\big( T, M_{p,T}(\mu)\big)$. Additionally, for any two initial data $\mu_0,\rho_0\in\ca{P}_{p}(\RR^d)$ such that the associated solutions $(\mu_t)_{t\in [0,T]},(\rho_t)_{t\in[0,T]}$ satisfy~\eqref{eq:strongclass}, we have the following stability estimate
\begin{equation}\label{PDE:Stability}
\ca{D}_p(\mu_t,\rho_t) \leq \,\,C\mathcal{D}_p(\mu_0,\rho_0) \, ,
\end{equation}
where the constant $C$ depends on $d$, $T$, $M_{\kappa, T}(\mu)$, $M_{\kappa, T}(\rho)$. 
\end{theorem}

\subsection{Well-posedness of time-inhomogeneous RDEs with drift}\label{ssec:wellposednessRDE}
In this section, we discuss well-posedness and stability of solutions to the RDE
\begin{equation}
\d Y_t=b(Y_t,\mu_t)\,\d t+\sigma(Y_t,\mu_t)\,\d\bfX_t \, ,
\label{eq:RDE}
\end{equation}
where $(\mu_t)_{t\in [0,T]}$ is a solution of the McKean--Vlasov PDE ~\eqref{MVPDE} in the sense of~\cref{defofPDE}. It is possible to interpret the above RDE as an RDE with time-dependent coefficients
\begin{equation}
\d Y_t=\bar{b}(t,Y_t)\,\d t+\bar{\sigma}(t,Y_t)\,\d\bfX_t \, ,
\label{eq:RDE2}
\end{equation}
where  $\bar{b},\bar{\sigma}$ are defined as follows
\begin{align}
\begin{split}
\bar{b}_{\mu}(t, y):= & \,b(y,\mu_t)\,, \\
\bar{\sigma}_{\mu}(t, y):= & \, \sigma(y,\mu_t) \, . 
\label{eq:bardef}
\end{split}
\end{align}
We drop the subscript $\mu$ when it is clear from the context which solution we are referring to. We first introduce the notion of a solution.
\begin{definition}[Solution to RDEs]\label{def:RDEsol}
Given a rough path $\bfX\in \scr{C}^{\alpha}([0,T];\RR^d)$ with $\alpha\in (1/3,1/2)$. For any $\xi\in\RR^d$, we say that $Y\in C^{\alpha}([0,T];\RR^d)$ is a solution of the RDE~\eqref{eq:RDE2} if $(Y,Y')\in\scr{D}_{X}^{2\alpha}([0,T];\RR^d)$ with $Y'=\bar{\sigma}(\cdot,Y)$ such that the  integral equation is satisfied, namely, for any $t\in[0,T]$,
\begin{align}
Y_t=\xi+\int_0^t \bar{b}(s,Y_s)\,\d s+\int_0^t\bar{\sigma}(s,Y_s)\,\d \bfX_s \, .
\end{align}
\end{definition}
\begin{definition}[Solution to backward RDEs]\label{def:backRDEsol}
Given a rough path $\bfX\in \scr{C}^{\alpha}([0,T];\RR^d)$ with $\alpha\in (1/3,1/2)$. Given $\xi\in\RR^d$, we say that $Y\in C^{\alpha}([0,T];\RR^d)$ is a solution to the backward RDE 
\begin{equation}
\d Y_t=\bar{b}(t,Y_t)\,\d t+\bar{\sigma}(t,Y_t)\,\d \bfX_t \, ,\label{eq:backRDE}
\end{equation}
if $(Y,Y')\in\scr{D}_{X}^{2\alpha}([0,T];\RR^d)$ with $Y'=\bar{\sigma}(\cdot,Y)$ such that the  integral equation is satisfied, namely, for any $t\in[0,T]$,
\begin{align}
\xi=Y_t+\int_t^T \bar{b}(s,Y_s)\,\d s+\int_t^T\bar{\sigma}(s,Y_s)\,\d \bfX_s \, .
\end{align}

\end{definition}

We now introduce the assumptions needed to treat~\eqref{eq:RDE2}.
\begin{assumption}
\label{ass:RDE}
Assume that $b:\R^d \times \ca{P}(\R^d)\to \R^d,\,\sigma:\R^d \times \ca{P}(\R^d)\to \R^{d\times d}$ are such that for any solution $(\mu_t)_{t\in[0,T]}$ of the PDE ~\eqref{MVPDE} in the sense of~\cref{defofPDE} satisfying~\eqref{eq:strongclass} and~\eqref{eq:timeregmu},  the associated $\bar{b},$ defined as in~\eqref{eq:bardef} satisfies 
\begin{align}
&\|\bar{b}(\cdot, 0)\|_0 <\infty\,, \quad \quad  \text{and} \quad \quad \quad [\bar{b}]_{0,\mathrm{Lip}} <\infty\,, 
\end{align}
and the associated $\bar{\sigma} $  satisfies
\begin{align}
	\lVert\bar{\sigma}\rVert_{1,4} <\infty 
\end{align}
or  has a linear form $$\bar{\sigma}(t, y)=a_0y+a_1(t)\,,$$
where $a_0$ is a linear map from $\R^d$ to $\R^{d}$ and $a_1(\cdot)$ is a continuously differentiable function.
\end{assumption}

We now present the main result of this section.
\begin{theorem}[Existence, uniqueness, and stability for the RDE]
\label{thm:MVRDE}
Let $T>0$, $\alpha\in (1/3,1/2)$ and $\bfX\in\scr{C}^{\alpha}([0,T];\RR^d)$. Given $\xi\in\RR^d$ and a solution $(\mu_t)_{t\in[0,T]}$ of the McKean--Vlasov PDE ~\eqref{MVPDE} (in the sense of~\cref{defofPDE}) satisfying~\eqref{eq:strongclass} and~\eqref{eq:timeregmu}, assume $b,\sigma$ satisfy~\cref{ass:RDE}. Then, there exists a solution $(Y,Y')\in\scr{D}_{X}^{2\alpha}([0,T];\RR^d)$ to the RDE~\eqref{eq:RDE} in the sense of \cref{def:RDEsol}. 

	Additionally, the constructed solution to the RDE~\eqref{eq:RDE2} is unique. To be more precise, consider 
	two solutions $(\mu_t)_{t\in[0,T]},(\rho_t)_{t\in[0,T]}$ of the PDE~\eqref{MVPDE} (in the sense of~\cref{defofPDE}) satisfying~\eqref{eq:strongclass} and~\eqref{eq:timeregmu}, 
	$\xi,\tilde{\xi}\in\RR^d$, $\bfX, \tilde{\bfX}\in \scr{C}^{\alpha}([0,T];\RR^d)$.	Let $(Y,Y')\in \scr{D}_{X}^{2\alpha}([0,T];\RR^d)$ and $(\tilde{Y},\tilde{Y}')\in \scr{D}_{\tilde{X}}^{2\alpha}([0,T];\RR^d)$ be the associated  solutions of the two RDEs and assume that for some constant $M<\infty$, it holds that $\|\bfX\|, \,\|\tilde{\bfX}\|\leq M, $
 $|Y_0|+|Y_0'|+\|Y,Y'\|_{X, 2\alpha}\leq M,\,|\tilde{Y}_0|+|\tilde{Y}_0'|+\|\tilde{Y},\tilde{Y}'\|_{\tilde{X}, 2\alpha}\leq M$.
Then, we have the following stability estimate:
\begin{align}\label{sta-RDE}
&\|Y,Y';\tilde{Y},\tilde{Y}'\|_{X,\tilde{X},2\alpha}\\ 
\lesssim\, &\,\, |\xi-\tilde{\xi}|+\varrho_{\alpha}(\bfX,\tilde{\bfX})+[\bar{\sigma}_{\mu}-\bar{\sigma}_{\rho}]_{{\rm{Lip}},0}+[\nabla(\bar{\sigma}_\mu-\bar{\sigma}_\rho)]_{{\rm{Lip}},0}+[\nabla(\bar{\sigma}_{\mu}-\bar{\sigma}_{\rho})]_{0,0}\\
&+[D^2(\bar{\sigma}_{\mu}-\bar{\sigma}_{\rho})]_{0,0}+\sup\limits_{s\neq t\in[0,T]}\frac{\int_s^t|\bar{b}_{\mu}(u,\tilde{Y}_u)-\bar{b}_{\rho}(u,\tilde{Y}_u)|\d u}{|t-s|^{2\alpha}}\,,
\end{align}
where the constant depends on $M, \alpha, T, \sigma$.\end{theorem}
\vspace{1mm}   
\subsection{Generation of random dynamical systems}\label{ssec:RDS}
 Before proceeding, we introduce the following additional assumption to ensure the measurability of the RDS we will construct.
\begin{assumption}\label{assume3}
We assume that there exists an increasing continuous function $G: [0,\infty)\to [0,\infty)$ with $G(0)=0$, such that for any $x, y\in\RR^d$, $\mu,\nu\in\ca{P}_{\kappa}(\RR^d)$,
 \begin{align}\label{condiformeasurability}
&|\sigma(x,\mu)-\sigma(y,\nu)|+\big|\big((\nabla\sigma)\sigma\big)(x,\mu)-\big((\nabla\sigma)\sigma\big)(y,\nu)\big|\\
\leq \,\,\,& 
F(\kappa,M_{\kappa}(\mu),M_{\kappa}(\nu))\left[|x-y|+G(\ca{D}_{\kappa}(\mu,\nu))\right]\,,
 \end{align}
where  for each $i=1,...,d$,
$$\big((\nabla \sigma)\sigma\big)_i(y,\mu)=\sum_{j=1}^{d}\sum_{k=1}^{d}\sigma_{jk}(y,\mu)\frac{\partial \sigma_{ik}}{\partial y_{j}}(y,\mu)\,.$$    
 \end{assumption}
 
\vspace{1mm} 
\begin{theorem}\label{exisrds} 
Let \cref{assumPDE,ass:RDE,assume3} be satisfied, and $p\geq \kappa$. Then, there exists a random dynamical system $\varphi$ on the product space $\uE_p$  associated with the coupled equations
\begin{subequations}
	\begin{align}
    \d Y_t &=  b(Y_t,\mu_t)\,\d t+\sigma(Y_t, \mu_t)\,\d {\mathbf{W}}_t\,,\,\,\,\,\,\,\,\,\,\,Y_0=y\,,\label{drift-SDE}\\
    \partial_t \mu_t &=-\nabla\cdot (b(\cdot,\mu_t)\mu_t)+\frac{1}{2}D^2:\big((\sigma\sigma^*)(\cdot,\mu_t)\mu_t\big)\,,\,\,\,\,\,\,\,\,\,\,\mu_0=\mu\,.\label{drift-PDE}
\end{align}
\end{subequations}
 More precisely, let $S_t:\cP_p(\R^d)\to \cP_p(\R^d)$ be the solution map of~\eqref{drift-PDE} and let $\phi_t(\omega,\cdot,\cdot):\Omega \times \R^d\times \cP_p(\R^d)\to \R^d$ be the solution map of~\eqref{drift-SDE}. Define \begin{align} 
&\varphi: [0,\infty) \times \Omega \times \uE_p \longrightarrow \uE_p\,,\\
&\,\,\,\,\,\,\,\,\,\,\,\,\,\,\,(t,\omega,y,\mu)
\longmapsto  
       ( \phi_{t}(\omega,y,\mu), 
        S_{t}(\mu))\,.
\end{align}
Then, $\varphi$ is an RDS over the metric dynamical system $(\Omega,\mathcal{F},(\theta_t)_{t\in \R},\mathbb{P})$ in the sense of~\cref{def:rds}.
Moreover, $\varphi$ is continuous in time under the metric $\mathsf{d}_p$ on $\uE_p$. 
\end{theorem}
\vspace{3mm}
\section{Proofs of the main results}\label{sec:proofs}
\vspace{3mm}
\subsection{\texorpdfstring{Proof of~\cref{ThmWell-posednessPDE}}{Proof of Theorem~\ref{ThmWell-posednessPDE}}}\label{ssub:proof_of_thmwell-posednesspde}
In this section, we will present the proof of the well-posedness of the PDE~\eqref{MVPDE} as stated in~\cref{ThmWell-posednessPDE}.  For any $\kappa\geq 2$, we first start with initial distribution $\mu_0\in\ca{P}_{p}(\RR^d)$ for $p>\kappa$ and obtain a unique solution in $\ca{P}_{p}(\RR^d)$, and then extend the solution to $\ca{P}_{\kappa}(\RR^d)$.

\paragraph{Existence in $\ca{P}_{p}(\RR^d)$.} For any $T>0$ and $n\geq1$,  let $\delta_n\coloneqq\frac{T}{n}, \delta^n_t\coloneqq \left\lfloor\frac{t}{\delta_n} \right \rfloor \delta_n.$  We then introduce the following approximating equations 
\begin{align}
\begin{split}\label{approxiSDE}
\d Y^{n}_t=b(Y^{n}_t,\mu^{n}_{\delta^n_t})\,\d t+\sigma(Y^{n}_t,\mu^{n}_{\delta^n_t})\,\d W_t\,,\,\,\,\,\,\,\,\,\,\,\,\,\,\,\,\,\,\,\ca{L}_{Y^{n}_0}=\mu_0\,,
\end{split}
\end{align}
where  $\mu^{n}_{\delta^n_t}\coloneqq\ca{L}_{Y^{n}_{\delta^n_t}}.$ More precisely, let $t_k^n=k\delta_n, k=0, 1,\ldots, n-1$. When  $t\in[0,t_{1}^n]$, this reads as \begin{align}
\begin{split}\label{approximation}
\d Y^{n}_t=b(Y^{n}_t,\mu^{n}_0)\,\d t+\sigma(Y^{n}_t,\mu^{n}_0)\,\d W_t\,,\,\,\,\,\,\,\,\,\,\,\,\,\,\,\,\,\,\,\ca{L}_{Y^n_0}=\mu_0\,.
\end{split}
\end{align}
Then, for any $k=1,\ldots,n-1$ and $t\in(t_{k}^n,t_{k+1}^n]$, we consider
\begin{equation}\label{eq4}
\d Y^{n}_t=b(Y^{n}_t,\mu^{n}_{t_{k}^n})\,\d t+\sigma(Y^{n}_t,\mu^{n}_{t_{k}^n})\,\d W_t\,,
\end{equation}
with initial value $Y^{n}_{t_{k}^n}$ and $\mu^{n}_{t_{k}^n}=\mathcal{L}_{Y^{n}_{t_{k}^n}}$.

Notice that (\ref{approximation}) is a classical SDE (not distribution-dependent). In view of the  condition~\eqref{localinstate}, growth condition (ii) of \cref{assumPDE}, the coefficients of (\ref{approximation}) satisfy the conditions in \cite[Theorem 3.1.1]{LiuRockner.2015.266}, hence, (\ref{approximation}) admits a unique strong solution $(Y_t^n)_{t\in [0,t_1^n]}$. Analogously, (\ref{eq4}) admits a unique solution $(Y_t^n)_{t\in (t_k^n, t_{k+1}^n]}, k=1, 2,...,n-1$.  Now, we claim that
\begin{align}\label{uniformlymomentestimate1}
\sup_{n\geq 1}\mathbb{E}\Bigg[\sup_{t\in[0,T]}|Y^{n}_t|^p\Bigg]\leq C_{p,T}\,,
\end{align}
and  for any $s, t\in [0,T]$, there exist positive constants $\alpha, \beta$ such that 
\begin{align}
\sup_{n\geq 1}\EE|Y_t^{n}-Y_s^{n}|^{\alpha}&\leq C_{T, \alpha, \kappa} |t-s|^{1+\beta}\,.\label{uniformcontinu1}
\end{align}
The proofs of these two claims are presented in \cref{Moment estimates}. According to~\eqref{uniformlymomentestimate1},~\eqref{uniformcontinu1}, and the Arzel\'a--Ascoli theorem \cite[Problem 4.11]{KaratzasShreve.Brownian}, the sequence of probability measures$${\mathbf{P}}^n\coloneqq\PP\circ (Y^n)^{-1}, n\geq 1,$$ on the sample space $(C([0,T]), \ca{B}(C([0,T])))$ are tight  and thus, there exists a subsequence $n$, denoted by $n$ for notational simplicity, and a probability measure $\mathbf{P}$ on the sample space such that $\mathbf{P}^n\xrightarrow{w} \mathbf{P}$ as $n\to \infty$. We thus have for any $t\in [0,T]$, as $n\to \infty$,
$$\mu_t^n\eqqcolon\Pi_t\mathbf{P}^n\xrightarrow{w} \Pi_t\mathbf{P}\eqqcolon\mu_t\,.$$
It follows from~\eqref{uniformlymomentestimate1} and the weak lower semi-continuity of  $L^p$-norms that
\begin{align}
	\sup_{t\in [0,T]}M_{p}(\mu_t)<\infty\,.
\end{align}
Moreover,  using the moment bounds~\eqref{uniformlymomentestimate1} with $p>\kappa$, we have the stronger convergence that  for all $t\in [0,T]$, as $n\to \infty$, 
\begin{align}
	&\mu_t^n\xrightarrow{d_{\kappa}} \mu_t\,,\label{es:dkappaconvergence}
	\end{align}
and by \cref{convergenceequiva},
	\begin{align}
	\mu_t^n\xrightarrow{\ca{D}_{\kappa}} \mu_t\,.\label{es:Dconvergence}
\end{align}
Furthermore,  for any $s, t\in [0,T]$,
\begin{align}
	d_{\kappa}^{\kappa}(\mu_t,\mu_s)&\leq \,\,d_{\kappa}^{\kappa}(\mu_t,\mu_t^n)+d_{\kappa}^{\kappa}(\mu_t^n,\mu_s^n)+d_{\kappa}^{\kappa}(\mu_s^n,\mu_s)\\
	&\leq\,\, d_{\kappa}^{\kappa}(\mu_t,\mu_t^n)+\sup_{n\geq 1}\EE|Y_t^n-Y_s^n|^{\kappa}+d_{\kappa}^{\kappa}(\mu_s^n,\mu_s)\,,
\end{align}
 together with~\eqref{uniformcontinu1},~\eqref{es:dkappaconvergence} and \cref{convergenceequiva}, it follows that as $s\to t$,
\begin{align}\label{es: Dtimeconvergence}
	\ca{D}_{\kappa}(\mu_t,\mu_s)\to 0\,.
\end{align}
Since  $\mu^n_t=\ca{L}_{Y^n_t},\,t\in [0,T],$ It\^{o}'s formula  \cite[Theorem 3.6]{KaratzasShreve.Brownian} shows that $\mu_t^n$ solves the PDE 
\begin{align}
\partial_t\mu^{n}_t=-\nabla\cdot\big(b(\cdot,\mu^{n}_{\delta^n_t})\mu^{n}_t\big)+\frac{1}{2}D^2:\big(\sigma\sigma^*(\cdot,\mu^{n}_{\delta^n_t})\mu^{n}_t\big)\,,\,\,\,\,\,\,\,\,\,\,\,\,\,\,\,\,\mu^{n}_0=\mu_0\,,
\end{align}
in $\ca{P}_{\kappa}(\RR^d)$, in the sense that for every $\phi\in \ca{K}_{\kappa}$,
\begin{align}\label{approintegral}
\begin{split}
	&\int_{\RR^d}\phi(y)\,\d\mu^{n}_t( y)-\int_{\RR^d}\phi(y)\,\d\mu_0( y)\\
=&\int_0^t\int_{\RR^d}\Big(\langle b(y,\mu^{n}_{\delta^{n}_s}),\nabla\phi(y)\rangle+\frac{1}{2}(\sigma\sigma^{*})(y,\mu^{n}_{\delta^{n}_s}):D^2\phi(y)\Big)\,\d\mu^{n}_s(y)\d s\,.
\end{split}
\end{align}

Now, we aim at showing $(\mu_t)_{t\in [0,T]}$ also solves the PDE~\eqref{MVPDE} in the sense of \cref{defofPDE} by taking the limit $n\to \infty$ on both sides of~\eqref{approintegral}.
We first note that as $n\to\infty$, the left hand side of~\eqref{approintegral} converges to
\begin{equation}
	\int_{\RR^d}\phi(y)\,\d\mu_t( y)\, ,
\end{equation} 
because $\phi\in \ca{K}_{\kappa}$ and $\mu_t^n\xrightarrow{\ca{D}_{\kappa}}\mu$.
For the right hand side of~\eqref{approintegral}, we rewrite
\begin{align}
	&\int_0^t\int_{\RR^d}\langle b(y,\mu^{n}_{\delta^{n}_s}),\nabla\phi(y)\rangle\;\d\mu^{n}_s( y)\,\d s\\
    =&\,\,\int_0^t\int_{\RR^d}\langle b(y,\mu_s),\nabla\phi(y)\rangle\;\d\mu^{n}_s( y)\d s+\int_0^t\int_{\RR^d}\langle b(y,\mu^{n}_{\delta^{n}_s})-b(y,\mu_s),\nabla\phi(y)\rangle\;\d\mu^{n}_s( y)\d s\\
	\eqqcolon&\,\,\int_0^t ({\rm{I}}) \,\d s+\int_0^t ({\rm{II}})\, \d s\,.
\end{align}
Using a similar argument to the one used in the proof of  \cref{convergenceequiva}, we are able to show that for any $|g|(y)\leq C(1+|y|^{\kappa})$, as $n\to \infty$,
\begin{align}\label{convergence: C2function}
	\int_{\RR^d} g(y) \; \d \mu^n(y)\to \int_{\RR^d} g(y) \; \d \mu(y)\,,
\end{align}
provided that $\ca{D}_{\kappa}(\mu^n, \mu)\to 0$ and $\sup_{n}M_{p}(\mu_n, \mu)<\infty$. Due to  growth condition of \cref{assumPDE} (iii), we have
$$\langle b(y,\mu_{s}),\nabla\phi(y)\rangle\leq \,F(\kappa, M_{\kappa}(\mu_s),0)(1+|y|^{\kappa})\leq F\Big(\kappa, \sup_{s\in [0,T]}M_{\kappa}(\mu_s),0\Big)(1+|y|^{\kappa})\,\, .$$ 
Thus, it follows from~\eqref{convergence: C2function} and the dominated convergence theorem that as $n\to \infty$, 
\begin{align}\label{convergence:b1}
\int_0^t ({\rm{I}})\;\d s \to \int_0^t\int_{\RR^d}\langle b(y,\mu_s),\nabla\phi(y)\rangle\;\d\mu_s( y)\d s\,.
\end{align}
By local Lipschitz continuity of \cref{assumPDE} (i), we see that
\begin{align}
&\langle \,b(y,\mu^{n}_{\delta^{n}_s})-b(y,\mu_s),\nabla\phi(y)\,\rangle \\
\leq \,\,\,&F\Big(\kappa, \sup_{n\geq 1}\sup_{s\in[0,T]}M_{\kappa}(\mu^n_s),\sup_{s\in[0,T]}M_{\kappa}(\mu_s)\Big)\big(1+|y|\big)^{\kappa}\mathcal{D}_{\kappa}(\mu^{n}_{\delta^{n}_s},\mu_s)\\
\leq\,\,\, &F\Big(\kappa, \sup_{n\geq 1}\sup_{s\in[0,T]}M_{\kappa}(\mu^n_s),\sup_{s\in[0,T]}M_{\kappa}(\mu_s)\Big)\big(1+|y|\big)^{\kappa}\Big(\ca{D}_{\kappa}\big(\mu_{\delta_s^n}^n, \mu_s^n\big)+\ca{D}_{\kappa}\big( \mu_s^n, \mu_s\big) \Big)\,.
\end{align}
It follows from the convergences~\eqref{es:Dconvergence},~\eqref{es: Dtimeconvergence} and  the dominated convergence theorem again, as $n\to \infty$,
$$\int_0^t ({\rm{II}}) \;\d s \to 0\,,$$
together with~\eqref{convergence:b1}, we derive that as $n\to \infty$\,,
\begin{align}\label{convergence:b}
	\int_0^t\int_{\RR^d}\langle b(y,\mu^{n}_{\delta^{n}_s}),\nabla\phi(y)\rangle\;\d\mu^{n}_s( y)\d s\to \int_0^t\int_{\RR^d}\langle b(y,\mu_s),\nabla\phi(y)\rangle\;\d\mu_s( y)\d s\,.
\end{align}
Using a similar argument, as $n\to \infty$\,,
\begin{align}\label{convergence:sigma}
	\frac{1}{2}\int_0^t\int_{\RR^d}(\sigma\sigma^{*})(y,\mu^{n}_{\delta^{n}_s}):D^2\phi(y)\,\d\mu^{n}_s( y)\d s\to \frac{1}{2}\int_0^t\int_{\RR^d}(\sigma\sigma^{*})(y,\mu_s):D^2\phi(y)\;\d\mu_s( y)\d s\,.
\end{align}
In combination of the convergences~\eqref{convergence:b},~\eqref{convergence:sigma}, we conclude that $(\mu_t)_{t\in [0,T]}$ solves PDE ~\eqref{MVPDE} in $\ca{P}_{\kappa}(\RR^d)$ in the sense of \cref{defofPDE}. Furthermore, by growth condition of \cref{assumPDE} (iii), we have the following time-continuity, for any $s, t\in [0,T]$,
\begin{align}
	\ca{D}_{\kappa}(\mu_t, \mu_s)&\leq\,\,\,\sup_{\phi\in\ca{K}_{\kappa}}\int_{\RR^d}\phi \; \d (\mu_t-\mu_s)\label{es:continuityintimebefore}\\
	&=\,\,\,\sup_{\phi\in\ca{K}_{\kappa}}\bigg\{\int_s^t\int_{\RR^d}\langle b(y,\mu_r),\nabla\phi(y)\rangle\;\d\mu_r( y)\d s\\
	&\;\;\;\;\;\;\;\;\;\;\;+\frac{1}{2}\int_s^t\int_{\RR^d}(\sigma\sigma^{*})(y,\mu_r):D^2\phi(y)\;\d\mu_r( y)\d s\bigg\}\\
	&\leq \,\,\,C|t-s|\,,
	\end{align}
where the constant $C$ depending on  $M_{\kappa,T}(\mu)$.
\vspace{3mm}

\paragraph{Uniqueness.} Let $(\mu_t)_{t\in [0,T]}$ and $(\rho_t)_{t\in [0,T]}$ be two solutions satisfying~\eqref{eq:strongclass}. Consider the following backward PDE for $u:[0,T]\times \RR^d\to \RR$,
\begin{align}
\begin{split}\label{backPDE}
\partial_t u&=-\langle  b(y, \mu_t), \nabla u\rangle-\frac{1}{2}\sigma\sigma^*(y, \mu_t): D^2u\,,\,\,\,\,\,\,\,\,\,\,\,\,t\in[0,T]\,,\\
u(T)&=\varphi,
\end{split}
\end{align}
where $\varphi\in \ca{K}_{\kappa}$. Moreover, we consider at the 
SDE
\begin{align}
\begin{split}\label{FeynmanSDE}
\d\xi_r&=b(\xi_r, \mu_r)\,\d r+\sigma(\xi_r, \mu_r)\,\d W_r\,,\;\;\;\;\;r\in [t,T]\,,\\
\xi_t&=y,
\end{split}
\end{align}
and set $\xi_r^{t, y}$ to be the unique solution evaluated at time $r$ started from $\xi_t=y$. Define 
\begin{align}
u(t, y)\coloneqq\EE\big[\varphi(\xi_T^{t, y})\big],
\end{align}
where $\varphi$ is the terminal condition of~\eqref{backPDE}. 
Because $|\varphi(x)|+|\nabla \varphi(x)|+|D^2 \varphi(x)|\leq C(1+|x|^{\kappa})$, $b$ and $\sigma$ are continuously differentiable with respect to $y$, it follows from \cite[V.7, Lemma 1]{Krylov.1995.271} that $u$ is continuously differentiable with respect to $y$ and fulfills 
\begin{equation}
	|\nabla u| (t, y)\leq C(1+|y|^{\kappa-1})\,,
	\label{firstderif}
\end{equation}
where the constants $C$ depending on $d, T, M_{\kappa, T}(\mu)$.
The proof is presented  in  \cref{Differentiability of backward PDE}. Furthermore, $u$ is continuously differentiable with respect to $t$ and twice continuously differentiable with respect to $y$ with 
\begin{equation}
	|D^2 u| (t, y)\leq C(1+|y|^{\kappa-1})\,,\label{secderif}
\end{equation}
where the constants $C$ depending on $d, T, M_{\kappa, T}(\mu)$, under the assumption that  $b, \sigma$ are twice continuously differentiable with respect to $y$, satisfy the growth conditions~\eqref{Growthb},~\eqref{Growthsigma},~\eqref{Growth1},~\eqref{Growthseconderiva}, and functions $(t, y)\mapsto b(y, \mu_t)$, $(t, y)\mapsto \sigma(y, \mu_t)$,  are jointly continuous in $(t, y)$. Most importantly, $u$ is a classical solution of the backward PDE~\eqref{backPDE} (See \cite[V.7, Theorem 4]{Krylov.1995.271}).  Now, we aim at proving that 
 $$\ca{D}_{\kappa}(\mu_T,\rho_T)\leq C \ca{D}_{\kappa}(\mu_0, \rho_0)$$ for a constant $C$,  which would imply the uniqueness of the solution to McKean--Vlasov PDE ~\eqref{MVPDE}.  To this end, informally speaking, we calculate the  derivative in time of
$$\int_{\RR^d}u(t, y)\, \d (\mu_t-\rho_t)(y)\,,$$
for almost every $t\in [0,T]$.  For all $h> 0$,  we Taylor expand as follows
\begin{align}
	&\frac{1}{h}\bigg[\int_{\RR^d}u(t+h, y)\, \d (\mu_{t+h}-\rho_{t+h})(y)-\int_{\RR^d}u(t,y)\, \d (\mu_t-\rho_t)(y)\bigg]\\
	=\,\,\,&\frac{1}{h}\bigg[\int_{\RR^d}u(t+h, y)\, \d (\mu_{t+h}-\rho_{t+h})(y)-\int_{\RR^d}u(t+h, y)\, \d (\mu_{t}-\rho_{t})(y)\bigg]\\
	&+\frac{1}{h}\bigg[\int_{\RR^d}u(t+h, y)\, \d (\mu_{t}-\rho_{t})(y)-\int_{\RR^d}u(t,y)\, \d (\mu_t-\rho_t)(y)\bigg]\\
	=\,\,\,&\frac{1}{h}\bigg[\int_{\RR^d}u(t, y)\, \d (\mu_{t+h}-\rho_{t+h})(y)-\int_{\RR^d}u(t, y)\, \d (\mu_{t}-\rho_{t})(y)\bigg]\\
	&+\int_{\RR^d}\partial_t u(t^*_{t,y}, y) \,\d (\mu_{t+h}-\rho_{t+h})(y)-\int_{\RR^d}\partial_t u(t^*_{t,y}, y)\, \d (\mu_{t}-\rho_{t})(y)\\
	&+\int_{\RR^d}\frac{1}{h} \big(u(t+h, y)-u(t, y)\big)\, \d (\mu_t-\rho_t)(y)\\
	=\,\,\,&\bigg[\frac{1}{h}\int_t^{t+h}\int_{\RR^d}\Big(\langle b(y,\mu_s),\nabla u(t, y)\rangle+\frac{1}{2}(\sigma\sigma^{*})(y,\mu_s):D^2 u(t, y)\Big)\, \d \mu_{s}(y)\d s\bigg]\\
	&-\bigg[\frac{1}{h}\int_t^{t+h}\int_{\RR^d}\Big(\langle b(y,\rho_s),\nabla u(t, y)\rangle+\frac{1}{2}(\sigma\sigma^{*})(y,\rho_s):D^2 u(t, y)\Big)\, \d \rho_{s}(y)\d s\bigg]\\
	&+\bigg[\int_{\RR^d}\partial_t u(t^*_{t,y}, y) \,\d (\mu_{t+h}-\rho_{t+h})(y)-\int_{\RR^d}\partial_t u(t^*_{t,y}, y)\, \d (\mu_{t}-\rho_{t})(y)\bigg]\\
	&+\int_{\RR^d}\frac{1}{h} \big(u(t+h, y)-u(t, y)\big)\, \d (\mu_t-\rho_t)(y)\\
	\eqqcolon& \,\,\,({\rm{I}})-({\rm{II}})+({\rm{III}})+({\rm{IV}})\,,
\end{align}
where $t^*_{t,y}\in[t,t+h]$, and in the second to last equality, we used the fact that for any $t\in [0,T]$, $u_t\in\ca{K}_{\kappa}$ and $(\mu_t)_{t\in [0,T]}$, $(\rho_t)_{t\in [0,T]}$ solves the McKean--Vlasov PDE in the sense of \cref{defofPDE}. In view of $\sup_{t\in [0,T]}M_{p}(\mu_t)<\infty$ and growth condition of \cref{assumPDE} (iii), the Lebesgue differentiation theorem then implies that at almost every $t\in [0,T]$, as $h\to 0$,
 \begin{equation}\label{es:Iconvergence}
 	({\rm{I}})\to \int_{\RR^d}\Big(\langle b(y,\mu_t),\nabla u(t, y)\rangle+\frac{1}{2}(\sigma\sigma^{*})(y,\mu_t):D^2 u(t, y)\Big)\, \d \mu_{t}(y),
 \end{equation}
  and 
 \begin{equation}\label{es:IIconvergence}
 	({\rm{II}})\to \int_{\RR^d}\Big(\langle b(y,\rho_t),\nabla u(t, y)\rangle+\frac{1}{2}(\sigma\sigma^{*})(y,\rho_t):D^2 u(t, y)\Big)\, \d \rho_{t}(y).
 \end{equation}
Since $u$ is a solution to backward PDE~\eqref{backPDE} satisfying the regularity estimates~\eqref{firstderif},~\eqref{secderif}, we have
  \begin{align}
 	\sup_{t^*_{t,y}\in [0,T]}\big|\partial_t u(t^*_{t,y}, y)\big|&\leq \sup_{t^*_{t,y}\in [0,T]}\big|b(y, \mu_{t^*_{t,y}})\big|\cdot |\nabla u(t^*_{t,y},y)|+\big|\sigma\sigma^*(y, \mu_{t^*_{t,y}})\big|\big|D^2u(t^*_{t,y},y)\big|\\
 	&\leq F\Big(\kappa, \sup_{t\in [0,T]}M_{\kappa}(\mu_t),0\Big)(1+|y|^{\kappa})\,,
 \end{align}
 and then, by time continuity~\eqref{es: Dtimeconvergence} of $(\mu_t)_{t\in [0,T]}$ and~\eqref{convergence: C2function}, as $h\to 0$,
 \begin{equation}\label{es:IIIconvergence}
 	({\rm{III}})\to 0\,. 	
 	 \end{equation}
 It is clear that, by dominated convergence theorem and the formulation of backward PDE~\eqref{backPDE}, as $h\to 0$,
 \begin{align}\label{es:IVconvergence}
 \begin{split}
 		({\rm{IV}})\to& \int_{\RR^d} \partial_t u(t, y)\,\d (\mu_{t}-\rho_{t})(y)\\
 	&=\int_{\RR^d}\left(-\langle \nabla b(y, \mu_t), u(t, y)\rangle-\frac{1}{2} \sigma\sigma^*(y, \mu_t): D^2u(t, y)\right)\,\,\d (\mu_{t}-\rho_{t})(y)\,.
 \end{split}
 	 \end{align}
Therefore, taking the convergences~\eqref{es:Iconvergence},~\eqref{es:IIconvergence},~\eqref{es:IIIconvergence},~\eqref{es:IVconvergence} into account, we  conclude that
\begin{align}
&\int_{\RR^d}u(t, y)\, \d (\mu_t-\rho_t)(y)\\
=\,\,\,&\int_{\RR^d}u(0, y)\, \d (\mu_0-\rho_0)(y)+\int_0^t\int_{\RR^d}\langle b(y,\mu_s)-b(y,\rho_s), \nabla u(s, y)\rangle\,\d\rho_s(y)\d s\\
&+\frac{1}{2}\int_0^t \int_{\RR^d}  \big(\sigma\sigma^*(y,\mu_s)-\sigma\sigma^*(y,\rho_s): D^2 u(s, y)\big)\d\rho_s(y)\d s\\
\leq\,\,\, &\int_{\RR^d}u(0, y)\, \d (\mu_0-\rho_0)(y)+C_F\int_0^T \int_{\RR^d}(|\nabla u|+|D^2 u|)(s,y)(1+|y|)\ca{D}_{\kappa}(\mu_s,\rho_s)\,\d\rho_s(y)\d s\\
\leq\,\,\, &\int_{\RR^d}u(0, y)\, \d (\mu_0-\rho_0)(y)+ C_F\int_0^T \ca{D}_{\kappa}(\mu_s,\rho_s)\d s\,,\label{es:uniqunessstability}
\end{align}
where $C_F\coloneqq F\Big(\kappa, \sup_{s\in [0,T]}M_{\kappa}(\mu_s),\sup_{s\in [0,T]}M_{\kappa}(\rho_s)\Big)$, and we made use of the regularity estimates~\eqref{firstderif},~\eqref{secderif} in the last inequality. Hence,  taking  $t=T$ and $\varphi$ over $\ca{K}_{\kappa}$, we see by Gr\"{o}nwall's inequality that
\begin{align}\label{es:stabefore}
\ca{D}_{\kappa}(\mu_T,\rho_T)&\leq C\ca{D}_{\kappa}(\mu_0,\rho_0)\,,
\end{align}
with the constant $C$ depends on  $T,\kappa$, $M_{\kappa, T}(\mu)$, $M_{\kappa, T}(\rho)$. The uniqueness of the solution to the McKean--Vlasov PDE~\eqref{MVPDE} follows. 
\vspace{2mm}
\paragraph{Existence in $\ca{P}_{\kappa}(\RR^d)$.}
Now,  one is ready to extend the solution theory to $\ca{P}_{\kappa}(\RR^d)$, i.e.\  for any $\mu_0\in\ca{P}_{\kappa}(\RR^d)$, there exists a unique solution $(\mu_t)_{t\in [0,T]}$ to~\eqref{MVPDE} in the sense of \cref{defofPDE} satisfying~\eqref{eq:strongclass} and~\eqref{eq:timeregmu}. To do this,  let $\ca{P}_{\kappa, T}(\R^d)$ be the space of  curves of probability measures on $\R^d$, i.e.\ $[0,T]\ni t\mapsto\mu_t\in\ca{P}_{\kappa}(\RR^d)$, equipped with the metric $\ca{D}_{\kappa,T}\coloneqq \sup_{t\in[0,T]}\ca{D}_{\kappa}(\mu_t,\nu_t)$ for any two curves $(\mu_t)_{t\in [0,T]}$, $(\nu_t)_{t\in [0,T]}$, which is a complete metric space due to the completeness of $(\ca{P}_{\kappa}(\RR^d), \ca{D}_{\kappa})$ (see~\cref{convergenceequiva}). We notice that for any $\mu_0\in \ca{P}_{\kappa}(\RR^d)$, there exists a sequence of $(\mu_0^m)_{m\geq 1}\subset \ca{P}_{p}(\RR^d), p>\kappa$, such that $\mu_0^m\xrightarrow{d_{\kappa}}\mu_0$ as $m\to\infty$. For each initial distribution $\mu_0^m\in\ca{P}_{p}(\RR^d)$, $m\geq 1$, we can find a corresponding solution $(\mu_t^m)_{t\in [0,T]}$ satisfying~\eqref{eq:strongclass} and~\eqref{eq:timeregmu}. Moreover, it follows from \eqref{eq:PDEsolu} with the choice $\phi(x)=|x|^{\kappa}$  and using~\cref{assumPDE}(ii) that $$\sup_{m\geq 1}\sup_{t\in[0,T]}M_{\kappa}(\mu^m_t)<\infty$$ by  taking . Furthermore, for any $m_1, m_2\geq 1$, it follows from the stability estimate~\eqref{es:stabefore} that
\begin{align}
	\sup\nolimits_{t\in[0,T]}\ca{D}_{\kappa}(\mu_t^{m_1}, \mu_t^{m_2})\leq\, C \ca{D}_{\kappa}(\mu_0^{m_1}, \mu_0^{m_2})\,,
\end{align} 
where the right hand side converges to $0$ as $m_1, m_2\to \infty$.  Thus, $\{(\mu_t^m)_{t\in [0,T]}\}_{m\geq 1}$ is a Cauchy sequence in $(\ca{P}_{\kappa,T}(\RR^d),\ca{D}_{\kappa,T})$, and then there is an element $(\mu_t)_{t\in[0,T]}\in \ca{P}_{\kappa,T}(\RR^d)$ such that $\sup_{t\in[0,T]}\ca{D}_{\kappa}(\mu_t^m, \mu_t)\to 0$ as $m\to\infty$. Therefore, due to the fact that 
\begin{align}
	\partial_t \mu_t^m &=-\nabla\cdot(b(\cdot,\mu_t^m)\mu_t^m)+\frac{1}{2}D^2:\left((\sigma\sigma^*)(\cdot,\mu_t^m)\mu_t^m\right) \, ,
\end{align}
using a similar argument as~\eqref{approintegral} and ~\eqref{es:uniqunessstability}, it follows that $(\mu_t)_{t\in [0, T]}$ is the unique solution to~\eqref{MVPDE} in $\ca{P}_{\kappa}(\RR^d)$. In addition, for any initial data $\mu_0\in\ca{P}_p(\RR^d)$, since there exists a solution to~\eqref{MVPDE} in $\ca{P}_p(\RR^d)$, one is able to improve the time-continuity~\eqref{es:continuityintimebefore} to 
\begin{align}
	\ca{D}_{p}(\mu_t,\mu_s)\leq C|t-s|\,,
\end{align}
for a constant $C$ depending on $p, M_{p,T}(\mu)$, which completes the proof of \cref{ThmWell-posednessPDE}.

\vspace{5mm}
\subsection{Proof of \texorpdfstring{\cref{thm:MVRDE}}{Theorem~\ref{thm:MVRDE}}} To prove this theorem, it suffices to prove a  general result on the well-posedness of time-inhomogeneous RDE with drift and diffusion coefficients satisfying the bounds of~\cref{ass:RDE}.
\begin{proposition} 
Let $T>0$, $\alpha\in (1/3, 1/2)$, and $\bfX\in\scr{C}^{\alpha}([0,T];\RR^d)$. Given $\xi\in\RR^d$, assume that $b: [0, T]\times \RR^d\to \RR^d$, $\sigma: [0, T]\times \RR^d\to \RR^{d\times d}$ satisfy the bounds in~\cref{ass:RDE}. Then, there exists a solution  $(Y, Y')\in \scr{D}_{X}^{2\alpha}([0,T];\RR^d)$, in the sense of \cref{def:RDEsol}, to the RDE
\begin{align}
	Y_t=\xi+\int_0^t b(s, Y_s)\,\d s+\int_{0}^t \sigma(s, Y_s)\,\d \bfX_s\,.
\end{align} 
Now fix $\xi,\tilde{\xi}\in\RR^d$, $\bfX,\tilde{\bfX}\in \scr{C}^{\alpha}([0,T];\RR^d)$ and let $(Y,Y')\in \scr{D}_{X}^{2\alpha}([0,T];\RR^d)$ and $(\tilde{Y},\tilde{Y}')\in \scr{D}_{\tilde{X}}^{2\alpha}([0,T];\RR^d)$ be the solutions of the following two RDEs
\begin{align}
Y_t=\xi+\int_0^t b(s,Y_s)\,\d s+\int_0^t \sigma(s,Y_s)\,\d \bfX_s\,,\;\;\;\;\;\;\tilde{Y}_t=\tilde{\xi}+\int_0^t\tilde{b}(s,\tilde{Y}_s)\,\d s+\int_0^t\tilde{\sigma}(s,\tilde{Y}_s)\,\d \tilde{\bfX}_s\,.
\label{eq:twoRDEs}
\end{align}
Let us assume for some $M< \infty$, $$\|\bfX\|, \|\tilde{\bfX}\|\leq M<\infty, $$
and $$|Y_0|+|Y_0'|+\|Y,Y'\|_{X, 2\alpha}\leq M,\;\;\;\;\;\;\;\;|\tilde{Y}_0|+|\tilde{Y}_0'|+\|\tilde{Y},\tilde{Y}'\|_{\tilde{X}, 2\alpha}\leq M.$$
Then, we have the following stability estimate:
\begin{align}
\|Y,Y';\tilde{Y},\tilde{Y}'\|_{X,\tilde{X},2\alpha}
&\,\lesssim\,\, |\xi-\tilde{\xi}|+\varrho_{\alpha}(\bfX,\tilde{\bfX})+[\sigma-\tilde{\sigma}]_{{\rm{Lip}},0}+[\nabla(\sigma-\tilde{\sigma})]_{{\rm{Lip}},0}+[\nabla(\sigma-\tilde{\sigma})]_{0,0}\\
&\,\,\,\,\,\,\,\,\,+[D^2(\sigma-\tilde{\sigma})]_{0,0}+\sup\limits_{s\neq t\in[0,T]}\frac{\int_s^t|b(u,\tilde{Y}_u)-\tilde{b}(u,\tilde{Y}_u)|\d u}{|t-s|^{2\alpha}}\,,
\end{align}
where the constant depends on $M, \alpha, T,\sigma$.
\end{proposition}

\vspace{2mm}
\begin{proof} 
{\bf Existence.} We first consider $\sigma(t,y)=a_0 y +a_1(t)$ in which case the RDE has the form
\begin{align}\label{linearRDE}
	Y_t=\xi+\int_0^t b(s, Y_s)\,\d s+\int_0^t a_0Y_s+a_1(s)\,\d \bfX_s\,.
\end{align}
Let $1/3<\alpha''<\alpha'<\alpha<1/2$. Define the map $\ca{F}:\scr{D}_X^{2\alpha'}\to \scr{D}_X^{2\alpha'}, (Y, Y')\mapsto (F, F')$ as follows
\begin{align}
	\ca{F}(Y, Y')=(F, F'):=\Big(\xi+\int_0^t b(s, Y_s)\,\d s+\int_0^t a_0Y_s+a_1(s)\,\d \bfX_s,\, a_0Y+a_1(\cdot)\Big)
\end{align}
and calculate $\|\ca{F}(Y, Y')\|_{X, 2\alpha'}$.  By the estimate of \cite[Theorem 4.10]{FrizHairer.2014.251}, we obtain
\begin{align}
\big|R^{F}_{s,t}\big|&\leq \,|a_0||Y_s'||\bb{X}_{s,t}|+\bigg|\int_s^t b(u, Y_u)\,\d u \bigg|\\
&\,\,\,\,\,\,\,\,+C_{\alpha'}\big(\|X\|_{\alpha'}\big\|R^{a_0Y+a_1(\cdot)}\big\|_{2\alpha'}+\|\bb{X}\|_{2\alpha'}\|a_0Y'\|_{\alpha'}\big)|t-s|^{3\alpha'}\\
	&\leq\, |a_0|\,\|Y'\|_0\,|\bb{X}_{s, t}|+C_{\alpha'}|a_0|\,\|X\|_{\alpha'}\big\|R^{Y}\big\|_{2\alpha'}|t-s|^{3\alpha'}+\bigg|\int_s^t b(u, Y_u)\,\d u \bigg|\\
	&\,\,\,\,\,\,+C_{\alpha'}\|X\|_{\alpha'}|a_1(t)-a_1(s)|\,|t-s|^{3\alpha'}+|a_0|C_{\alpha'}\|\bb{X}\|_{2\alpha'}\|Y'\|_{\alpha'}|t-s|^{3\alpha'}
\end{align}
and then the $2\alpha'$-H$\rm{\ddot{o}}$lder norm 
\begin{align}
	\|R^{F}\|_{2\alpha'}\leq &\,\,|a_0|\|Y'\|_0\|\bb{X}\|_{2\alpha'}+|a_0|C_{\alpha'}\|X\|_{\alpha'}\big\|R^{Y}\big\|_{2\alpha'}T^{\alpha'}+C_{\alpha'}\|X\|_{\alpha'}\|a_1\|_{2\alpha'}T^{3\alpha'}\\
	&+C_{\alpha'}|a_0|\,\|\bb{X}\|_{2\alpha'}\|Y'\|_{\alpha'}T^{\alpha'}+\|b(\cdot, 0)\|_0T^{1-2\alpha'}+C\|Y\|_0T^{1-2\alpha'}\\
	\leq &\, \,|a_0|\,\|Y'\|_0\|\bb{X}\|_{2\alpha'}+C_{\alpha'}|a_0|\,\|X\|_{\alpha'}\big\|R^{Y}\big\|_{2\alpha'}T^{\alpha'}+C_{\alpha'}C_{a_1}\|X\|_{\alpha'}T^{1+\alpha'}\\
	&+C_{\alpha'}|a_0|\,\|\bb{X}\|_{2\alpha'}\|Y'\|_{\alpha'}T^{\alpha'}+\|b(\cdot, 0)\|_0T^{1-2\alpha'}\\
	&+C|\xi|T^{1-2\alpha'}+C(|a_0||\xi|+|a_1(0)|)\|X\|_{\alpha'}T^{1-\alpha'}+C\|R^Y\|_{2\alpha'}T\,,
\end{align}
where $C_{a_1}$ is the Lipschitz constant of the function $a_1(\cdot)$ and we use the estimate
$$\|Y\|_0\leq |\xi|+|Y_0'|\,\|X\|_{\alpha'}T^{\alpha'}+\|R^Y\|_{2\alpha'}T^{2\alpha'}$$
in the last inequality. Notice that
\begin{align}\label{es:alpha}
	\|a_0Y+a_1(\cdot)\|_{\alpha'}\leq |a_0|\,\|Y\|_{\alpha'}+\|a_1\|_{\alpha'}\,,
\end{align} 
together with 
\begin{align}
	\|Y\|_{\alpha'}&\leq \|Y'\|_0\|X\|_{\alpha'}+T^{\alpha'}\|R^{Y}\|_{2\alpha'},\\
	&\leq |Y'_0|\,\|X\|_{\alpha'}+T^{\alpha'}\|Y'\|_{\alpha'}\|X\|_{\alpha'}+T^{\alpha'}\|R^{Y}\|_{2\alpha'}\,,
\end{align}
we derive that 
\begin{align}\label{es:F}
	\|\ca{F}(Y, Y')\|_{X, 2\alpha'}\leq&\, |a_0|(|Y_0'|\,\|X\|_{\alpha'}+T^{\alpha'}\|X\|_{\alpha'}\|Y'\|_{\alpha'}+T^{\alpha'}\|R^Y\|_{2\alpha'})+C_{a_1}T^{1-\alpha'}\\
	&+|a_0|\,\|Y'\|_0\|\bb{X}\|_{2\alpha'}+|a_0|C_{\alpha'}\|X\|_{\alpha'}\big\|R^{Y}\big\|_{2\alpha'}T^{\alpha'}\\
	&+C_{\alpha'}C_{a_1}\|X\|_{\alpha'}T^{1+\alpha'}+|a_0|C_{\alpha'}\|\bb{X}\|_{2\alpha'}\|Y'\|_{\alpha'}T^{\alpha'}\\
	&+\|b(\cdot, 0)\|_0T^{1-2\alpha'}+C|\xi|T^{1-2\alpha'}\\
	&+C(|a_0||\xi|+|a_1(0)|)\|X\|_{\alpha'}T^{1-\alpha'}+C\|R^Y\|_{2\alpha'}T\\
	\leq&\,C_{a_0, a_1, \alpha', b, T}\big(T^{1-2\alpha'}+\|X\|_{\alpha}T^{\alpha-\alpha'}+\|\bb{X}\|_{2\alpha}T^{2(\alpha-\alpha')}\big)\\
	&\cdot \big(1+|\xi|+\|Y, Y'\|_{X,2\alpha'}\big)\,,
   	\end{align}
where we choose, for a fixed number $\lambda\in (0,1)$, sufficiently small $T_0$ such that 
$$C_{a_0, a_1, \alpha', b, T}\big(T^{1-2\alpha'}+\|X\|_{\alpha}T^{\alpha-\alpha'}+\|\bb{X}\|_{2\alpha}T^{2(\alpha-\alpha')}\big)\leq \lambda\,.$$
We restrict the controlled rough path to the following convex compact set in $\scr{D}^{2\alpha''}_{X}([0,T_0];\RR^d)$,
\begin{align}
B_{[0, T_0]}\coloneqq\Big\{(Y, Y')\in \scr{D}_{X}^{2\alpha'}([0,T_0];\RR^d): &Y_0=\xi, Y_0'=a_0\xi+a_1(0), \\
&\|Y, Y'\|_{X,2\alpha'}\leq \frac{\lambda}{1-\lambda}\big(1+|\xi|\big)\Big\},
\end{align}
where the compactness follows from the \cref{es:compactnessofcontrolledroughpath}, the convexity follows from the linearity of the controlled rough path, thus, by~\eqref{es:F},
\begin{align}
	\|\ca{F}(Y, Y')\|_{X, 2\alpha''}\leq \frac{\lambda}{1-\lambda}\big(1+|\xi|\big)\,,
\end{align}
which indicates that $\ca{F}: B_{[0, T_0]}\to B_{[0,T_0]}$. Together with the boundedness and continuity of $\ca{F}$, Schauder fixed point theorem yields that there exists a fixed point $(Y_t)_{t\in[0, T_0]}$ of the map $\ca{F}$, which is the solution of RDE~\eqref{linearRDE} on the interval $[0,T_0]$. 
Furthermore, since $\bfX\in \mathscr{C}^{\alpha}([0, T_0];\RR^d)$, it turns out that this solution $Y$ is automatically an element in $\scr{D}_X^{2\alpha}([0, T_0];\R^d)$ following similar arguments as in the proof of \cite[Theorem 8.3]{FrizHairer.2014.251}. Indeed, since $|Y_{s,t}|\leq \|Y'\|_{0}|X_{s,t}|+\|R^Y\|_{2\alpha''}|t-s|^{2\alpha''}$, $X\in C^{\alpha}$, and $2\alpha''>\alpha $,  it holds  that $Y\in C^{\alpha}$.  Additionally, $Y'=\sigma(\cdot,Y)\in C^{\alpha}$ since $\|\sigma\|_{1,3}<\infty$. Finally,  $\lVert R^Y\rVert_{2\alpha}<\infty $ using the fact that ${\bb{X}}\in C^{2\alpha}$ and
\begin{align}
|R^{Y}_{s,t}|&=|Y_{s,t}-Y_s'X_{s,t}|\\
&\leq \bigg|\int_s^tb(r,Y_r)\,\d r+\int_s^t \sigma(r,Y_r)\,\d\bfX_r-\sigma(s,Y_s)X_{s,t}\bigg|\\
&\leq (C\|Y\|_{0}+\|b(\cdot,0)\|_{0})|t-s|+|\nabla \sigma(s,Y_s)Y_s'|\cdot {\bb{X}}_{s,t}+O(|t-s|^{3\alpha})\\
&\leq (C\|Y\|_{0}+\|b(\cdot,0)\|_{0})|t-s|+[\nabla \sigma]_{0,0}\|Y'\|_{0}\cdot {\bb{X}}_{s,t}+O(|t-s|^{3\alpha})\,.
\end{align}

Notice that $T_0$ is independent of the initial value $\xi$, so we are able to construct a solution on $[kT_0, (k+1)T_0]$, for $k=1,...,k_0$ with $(k_0+1)T_0\leq T$, by iteratively applying the Schauder fixed point theorem on the  set
\begin{align}
	B_{[kT_0, (k+1)T_0]}\coloneqq\Big\{(\bar{Y}, \bar{Y}')\in &\scr{D}_{X}^{2\alpha'}([kT_0, (k+1)T_0];\RR^d): \bar{Y}_{kT_0}=Y_{kT_0}, \\&\bar{Y}_{kT_0}'=a_0Y_{kT_0}+a_1(kT_0),
	\|\bar{Y}, \bar{Y}'\|_{X,2\alpha'}\leq \frac{\lambda}{1-\lambda}\big(1+|Y_{kT_0}|\big)\Big\}.
\end{align} 
Hence, we can construct a solution on $[0,T]$ by gluing together the solutions on the sub-intervals $[kT_0, (k+1)T_0]$, $k=0,1,...,k_0$.

\vspace{2mm}
Secondly, we treat the case $\|\sigma\|_{1,4}<\infty$. The proof in this setting relies on  treating the RDE without drift first and then obtaining existence for the RDE with drift by using the so-called Doss--Sussmann transformation. To this end, for any $s\in [0,T]$, we first consider the following RDE without drift
\begin{equation}\label{RDEnodrift}
	 \d Y_t=\sigma(t, Y_t)\,\d \bfX_t\,,\;\;\;\;\;\;\;\;\;\;t\in [s,T] \, ,
	\end{equation}
The above RDE has a unique solution $Y \in \scr{D}^{2\alpha}_{X}([s,T];\R^d)$. The proof of this follows from a modification of the fixed-point argument in~\cite[Theorem 8.3]{FrizHairer.2014.251} to deal with the time-inhomogeneity of the coefficient. We have included the proof in the appendix (see \cref{well-RDE}) for the sake of completeness. 
We now define $\Phi(s,t,\xi)$ to be solution of~\eqref{RDEnodrift} started from time $s$ with $\Phi(s,s,\xi)=\xi$. One can show that $\Phi(s,t,\xi)$ is continuously differentiable in $\xi$  and its gradient $(\nabla \Phi)(s,t,\xi)$ is the unique solution of the matrix-valued  linearised RDE 
\begin{align}
	 \d \zeta_t=(\nabla \sigma)(t, \Phi(s,t,\xi))\zeta_t \,\d \bfX_t \, ,\;\;\;\;\;\;\;\;\;\;t\in [s,T]\, ,
\end{align}
started at time $s$ with initial data $(\nabla \Phi)(s,s,\xi)= \mathrm{Id}_{\R^d}$.   We also consider the following backward RDE
\begin{align}\label{backRDEnodrift}
	\d h_t=\sigma(t, h_t)\,\d \bfX_t\,,\,\,\,\,\,\,\,\,\,\,\,\,\,\,\,t\in[s,T]\,,
\end{align}
which also has a unique solution in the sense of \cref{def:backRDEsol}.  We denote by $\Psi(s,t,\delta)$ the solution of~\eqref{backRDEnodrift} with terminal condition $\Psi(t,t,\delta)=\delta$. Note that  $\Psi(s,t,\delta)$ is also continuously differentiable in $\delta$  and its gradient $(\nabla \Psi)(s,t,\delta)$ is the unique solution of the matrix-valued  linearised backward RDE
\begin{equation}
	\d \eta_t=(\nabla\sigma)(t, \Psi(s,t,\delta))\eta_t \,\d \bfX_t \, ,\;\;\;\;\;\;\;\;\;\;t\in [s,T] \, .\label{eq:linearisedbackwardRDE}
\end{equation}
with terminal condition $(\nabla \Psi)(t,t,\delta)=\mathrm{Id}_{\R^d}$.
 Then  we have 
\begin{align}
	\Phi(0, t,\Psi(0, t,\delta))=	\Psi(0, t,\Phi(0, t,\delta))&=\delta\,,\label{eq:compositionPsi}
    \end{align}
    and by the inverse function theorem, it holds that 
\begin{align}
	(\nabla \Phi)(0,t,\Psi(0, t,\delta)) (\nabla \Psi)(0, t,\delta)=
	(\nabla \Psi)(0,t,\Phi(0, t,\delta)) (\nabla \Phi)(0, t,\delta)&=\mathrm{Id}_{\R^d}\,.\label{eq:inversefunctionPsi}
\end{align}
To see ~\eqref{eq:compositionPsi},  we notice $F(u)\coloneqq\Phi(0, u,\Psi(0, t,\delta)), u\in [0, t],$ solves
\begin{align}
	\label{eq:intermediateRDE}F(u)&=\Psi(0, t,\delta)+\int_0^u \sigma(r, F(r))\,\d \bfX_r\\
	&=\delta-\int_0^t\sigma(r,\Psi(r, t,\delta))\,\d \bfX_r+\int_0^u \sigma(r, F(r))\,\d \bfX_r\,.
\end{align} 
Clearly, $u \mapsto F(u)$ is the unique controlled rough path solving the RDE~\eqref{eq:intermediateRDE}. Notice however that if we choose $F(u)=\Psi(u,t,\delta)$ then we also have a solution to~\eqref{eq:intermediateRDE}, thus implying~\eqref{eq:compositionPsi} (by repeating this argument with $\Psi$) and~\eqref{eq:inversefunctionPsi}.

Consider now the ODE 
\begin{align}\label{eq:DossODE}
	\dot{z}_t=(\nabla \Psi)(0, t,\Phi(0,t,z_t))\,b(t,\Phi(0,t,z_t))\eqqcolon F(t, z_t),\,\,\,\,\,\,\,\,t\in[0,T]\,.
\end{align}
One can check that $\nabla \Psi$ is globally Lipschitz continuous and bounded in the state variable,  and from~\eqref{RDEnodrift} to see that $\Phi$ is globally Lipschitz continuous  so the composition  $(\nabla \Psi)(0, t,\Phi(0,t,z_t))$ is globally Lipschitz continuous and bounded in $z_t$. 
Since $b$ is also assumed to be Lipschitz continuous with  linear growth,  it follows that $F(t, z_t)$ fulfills the local Lipschitz continuity with the linear growth. Hence, the ODE~\eqref{eq:DossODE} is well-posed on $[0,T]$. We claim that $Y_t\coloneqq\Phi(0,t,z_t)$, $t\in [0,T]$, is a solution to the RDE
\begin{align}
	\d Y_t=b(t, Y_t)\,\d t+\sigma(t, Y_t)\,\d\bfX_t \, .\end{align}
Indeed, it follows from the definition of $\Phi$ that
\begin{align}
	Y_{s, t}&=\Phi(0, t, z_t)-\Phi(0,t, z_s)+\Phi(0, t, z_s)-\Phi(0,s, z_s)\\
	&=\Big[\Phi(0, t, z_t)-\Phi(0,t, z_s)\Big]+\int_0^t  \sigma(u, \Phi(0,u,z_s))\,\d\bfX_u-\int_0^s  \sigma(u, \Phi(0,u,z_s))\,\d\bfX_u \\
	 &=\Big[\Phi(0, t, z_t)-\Phi(0,t, z_s)\Big]+\int_s^t  \sigma(u, \Phi(0,u,z_s))\,\d\bfX_u\,.\label{incrementY}
	\end{align}
For the first term of~\eqref{incrementY}, we can Taylor expand to obtain
\begin{align}
	&\big|\Phi(0,t, z_t)-\Phi(0,t, z_s)-(\nabla \Phi)(0, s,   z_s)z_{s, t}\big|\\
	\leq \,\,\,&\big|\Phi(0,t, z_t)-\Phi(0,t, z_s)-(\nabla \Phi)(0, t,   z_s)z_{s, t}\big|+\big|(\nabla \Phi)(0, t,   z_s)-(\nabla \Phi)(0, s,   z_s)\big|\cdot|z_{s, t}|\\
	\leq \,\,\,&C|z_{s,t}|^2+\Big|\int_s^t (\nabla \sigma)(u, \Phi(0,u,z_s))(\nabla \Phi )(0, u,   z_s) \,\d\bfX_u\Big|\cdot|z_{s, t}|\\
	\leq\,\,\, & C|t-s|^{1+\alpha}\,,\label{incrementY1}
\end{align}
where in the last inequality, we have used the fact that $\nabla \sigma (\cdot, \Phi(0,\cdot,z_s)), \nabla \Phi(0,u,z_s)$ are both controlled rough paths and that the product of two controlled rough paths is still a controlled rough path along with the estimate from \cite[Theorem 4.10, Equation (4.22)]{FrizHairer.2014.251} . Using ~\eqref{eq:inversefunctionPsi} and~\eqref{eq:DossODE} and Taylor expanding, we arrive at 
\begin{align}
	(\nabla \Phi)(0, s,   z_s)z_{s,t}
	&=(\nabla \Phi) (0, s,   z_s)\int_s^t F(u,  z_u)\, \d u\\
	&=(\nabla\Phi)(0, s,   z_s)F(s,  z_s)|t-s|+O(|t-s|^{1+\alpha})\\
	&=(\nabla\Phi)(0, s,   z_s)(\nabla  \Psi)(0, s,  \Phi(0,s, z_s))\,b(s,\Phi(0,s, z_s))|t-s|+O(|t-s|^{1+\alpha})\\
	&= b(s,Y_s)|t-s|+O(|t-s|^{1+\alpha})\,,\label{derivativePhi}
	\end{align}
	where in the second equality, we use the fact that $F(u,z_u) \in C^{\alpha}([0,T];\R^d)$. By the definition of rough integral,  the second term of~\eqref{incrementY} follows that 
\begin{align}\label{eq:increamentintegralsigma}
	\int_s^t  \sigma(u, \Phi(0,u,z_s))\,\d\bfX_u =& \,\sigma(s, \Phi(0,s,z_s))X_{s,t}\\&+(\nabla \sigma)(s, \Phi(0,s,z_s))[\Phi(0,s,z_s)]'{\bb{X}}_{s,t}\\
	&+O(|t-s|^{3\alpha})\,.
\end{align}
In combination of~\eqref{incrementY},~\eqref{incrementY1},~\eqref{derivativePhi} and~\eqref{eq:increamentintegralsigma}, it holds that
\begin{align}
	Y_{s, t}=b(s, Y_s)|t-s|&+\sigma(s, \Phi(0,s,z_s))X_{s,t}+(\nabla \sigma)(s, \Phi(0,s,z_s))[\Phi(0,s,z_s)]'{\bb{X}}_{s,t}\\
	&+O(|t-s|^{3\alpha})\,,
\end{align} 
i.e.\ $Y$ is controlled by $\bfX$ with the $Y_s'=(\nabla \sigma)(s, \Phi(0,s,z_s))[\Phi(0,s,z_s)]'=(\nabla \sigma)(s, Y_s)Y_s'$. Taking any partition $\Pi$ of $[0,T]$ for $|\Pi|=\max_{[s,t]\in \Pi}|t-s|\to 0$, we arrive at 
\begin{align}
	Y_{0, T}=Y_0+\sum_{[s,t]\in\Pi} Y_{s,t}=\xi+\int_0^T b(u, Y_u)\,\d u+\int_0^T \sigma(u, Y_u)\,\d \bfX_u\,,
\end{align} 
which proves the existence part.

\vspace{3mm}
\paragraph{Uniqueness.} Now, we proceed to the proof of the stability estimate which enables us to obtain the uniqueness of the solution to RDE~~\eqref{eq:RDE}. Let $(Y,Y')$ and $(\tilde{Y},\tilde{Y}')$ denote the solutions of the two RDEs in~\eqref{eq:twoRDEs} and define 
$$(\Xi,\Xi')\coloneqq(\sigma(\cdot,Y), \sigma(\cdot,Y)'), \;\;\;\;\;(\tilde{\Xi},\tilde{\Xi}')\coloneqq(\tilde{\sigma}(\cdot,\tilde{Y}), \tilde{\sigma}(\cdot,\tilde{Y})'),$$ 
along with
\begin{align}
	(Z,Z')&\coloneqq\left(\xi+\int_0^{\cdot}b(s,Y_s)\;\d s+\int_0^{\cdot}\Xi_s \;\d\bfX_s,\Xi\right),\\
	(\tilde{Z},\tilde{Z}')&\coloneqq\left(\tilde{\xi}+\int_0^{\cdot}\tilde{b}(s,\tilde{Y}_s)\;\d s+\int_0^{\cdot}\tilde{\Xi}_s\;\d{\tilde\bfX}_s,\tilde{\Xi}\right). 
\end{align}
We can now make use of well-known stability results for rough integrals  provided to us by \cite[Theorem 4.17]{FrizHairer.2014.251} which yield \footnote{The theorem cited here does not consider a drift in the RDE but the estimate below can be obtained by a straightforward modification of the same proof.}
\begin{align}\label{es:XtildeX}
&\|Y,Y'; \tilde{Y}, \tilde{Y}'\|_{X,\tilde{X},2\alpha}\\ 
= \, &\|Z,Z'; \tilde{Z}, \tilde{Z}'\|_{X,\tilde{X},2\alpha}\\
\leq & \, C_{M,\alpha} \left[\varrho_{\alpha}(\bfX,{\tilde\bfX})
+|\xi-\tilde{\xi}|+T^{\alpha} \|\Xi,\Xi';\tilde{\Xi},\tilde{\Xi}'\|_{X,\tilde{X}, 2\alpha}+\bigg\|\int_0^{\cdot}b(s,Y_s)-\tilde{b}(s,\tilde{Y}_s)\;\d s\bigg\|_{2\alpha}\right] \, .
\end{align}
Applying the triangle inequality, we see that 
\begin{align}\label{es:XitildeXi}
&\|\Xi,\Xi';\tilde{\Xi},\tilde{\Xi}'\|_{X,\tilde{X}, 2\alpha}\\
= \,\,\,&\|\sigma(\cdot,Y), \sigma(\cdot,Y)';\tilde{\sigma}(\cdot,\tilde{Y}), \tilde{\sigma}(\cdot,\tilde{Y})'\|_{X,\tilde{X}, 2\alpha}\\
\leq \,\,\,& \|\sigma(\cdot,Y), \sigma(\cdot,Y)';\sigma(\cdot,\tilde{Y}), \sigma(\cdot,\tilde{Y})'\|_{X,\tilde{X}, 2\alpha}+\|\sigma(\cdot,\tilde{Y}), \sigma(\cdot,\tilde{Y})';\tilde{\sigma}(\cdot,\tilde{Y}), \tilde{\sigma}(\cdot,\tilde{Y})'\|_{X,\tilde{X}, 2\alpha} \, . \label{eq:trianglefftilde}
\end{align}
 In order to control the first term on the right hand side of the above expression, we need to modify the proof of \cite[Theorem 7.6]{FrizHairer.2014.251} on the stability of compositions of controlled rough paths by regular functions, since the function we are composing with is now time-dependent. We will see that the first term on the right hand side of~\eqref{eq:trianglefftilde} is bounded by
 \begin{align}\label{es:stability1}
 C_M\left(\|X-\tilde{X}\|_{\alpha}+|Y_0-\tilde{Y}_0|+|Y_0'-\tilde{Y}_0'|+\|Y,Y';\tilde{Y},\tilde{Y}'\|_{X,\tilde{X}, 2\alpha}\right)\,,
 \end{align}
 and the second term on the right hand side of~\eqref{eq:trianglefftilde} is controlled by 
 $$C_{M}\big([\nabla(\sigma-\tilde{\sigma})]_{{\rm{Lip}},0}+[\sigma-\tilde{\sigma}]_{{\rm{Lip}},0}+[\nabla(\sigma-\tilde{\sigma})]_{0,0}+[D^2(\sigma-\tilde{\sigma})]_{0,0}\big)\,.$$
 We left the proof in \cref{es:compostionsigma}.
 Moreover, by the Lipschitz continuity of  $b$ with respect to the state variable, 
\begin{align}
&\bigg\|\int_0^{\cdot}b(s,Y_s)-\tilde{b}(s,\tilde{Y}_s)\,\d s\bigg\|_{2\alpha}\\
=\,\,\,&\sup\limits_{s\neq t\in[0,T]}\frac{\Big|\int_s^tb(u,Y_u)-b(u,\tilde{Y}_u)\,\d u\Big|+\Big|\int_s^tb(u,\tilde{Y}_u)-\tilde{b}(u,\tilde{Y}_u)\,\d u\Big|}{|t-s|^{2\alpha}}\\
\leq\,\,\,&\|Y-\tilde{Y}\|_{0}|t-s|^{1-2\alpha}+\sup\limits_{s\neq t\in[0,T]}\frac{\int_s^t|b(u,\tilde{Y}_u)-\tilde{b}(u,\tilde{Y}_u)|\,\d u}{|t-s|^{2\alpha}}\\
\leq\,\,\,&(|Y_0-\tilde{Y}_0|+T^{\alpha}\|Y-\tilde{Y}\|_{\alpha})|t-s|^{1-2\alpha}+\sup\limits_{s\neq t\in[0,T]}\frac{\int_s^t|b(u,\tilde{Y}_u)-\tilde{b}(u,\tilde{Y}_u)|\,\d u}{|t-s|^{2\alpha}}\\
\lesssim \,\,\,&\big(|\xi-\tilde{\xi}|+\varrho_{\alpha}(\bfX,\tilde{\bfX})+T^{\alpha}\|Y,Y';\tilde{Y},\tilde{Y}'\|_{X,\tilde{X}, 2\alpha}\big)+\sup\limits_{s\neq t\in[0,T]}\frac{\int_s^t|b(u,\tilde{Y}_u)-\tilde{b}(u,\tilde{Y}_u)|\,\d u}{|t-s|^{2\alpha}}\,,
\end{align}
where the last inequality is due to the estimate  \cite[(4.30)]{FrizHairer.2014.251}. Thus, combining the above estimates together and taking $T$ small enough, we arrive at
\begin{align}
\|Y,Y';\tilde{Y},\tilde{Y}'\|_{X,\tilde{X},2\alpha}\lesssim &\,\,|\xi-\tilde{\xi}|+\varrho_{\alpha}(\bfX,\tilde{\bfX})+[\sigma-\tilde{\sigma}]_{{\rm{Lip}},0}+[\nabla(\sigma-\tilde{\sigma})]_{{\rm{Lip}},0}+[\nabla(\sigma-\tilde{\sigma})]_{0,0}\\
&\,\,\,+[D^2(\sigma-\tilde{\sigma})]_{0,0}+\sup\limits_{s\neq t\in[0,T]}\frac{\int_s^t|b(u,\tilde{Y}_u)-\tilde{b}(u,\tilde{Y}_u)|\,\d u}{|t-s|^{2\alpha}}\,.
\end{align}
\end{proof}
\subsection{Proof of \texorpdfstring{\cref{exisrds}}{Theorem~\ref{exisrds}}}
Before we can present the proof of~\cref{exisrds}, we consider the following sequence of systems 
\begin{subequations}
\begin{align}
    \d Y_t^n &=  b(Y_t^n, \mu_t)\,\d t+\sigma(Y_t^n, \mu_t) \,\d {{\mathbf{W}^n_t}}-2^{-1}\big((\nabla\sigma)\sigma\big)(Y_t^n,\mu_t)\,\d t,\;\;\;\;Y_0^n=y\,,\label{ODE}\\
    \partial_t \mu_t &=-\nabla\cdot (b(\cdot,\mu_t)\mu_t)+\frac{1}{2}D^2:\big((\sigma\sigma^*)(\cdot,\mu_t)\mu_t\big),\;\;\;\;\;\;\;\;\;\;\mu_0=\mu\,,\label{PDEb}
\end{align}
\end{subequations}
where $(\mathbf{W}^n=(W^n,{\bb{W}^n}))_{n\geq 1}$ is the dyadic piecewise linear approximation of $\mathbf{W}$. Notice that term $2^{-1}(\nabla \sigma)\sigma$ comes from the It\^{o}-Stratonovich correction. 
We start by presenting the following preliminary result.

\vspace{3mm}
\begin{proposition}\label{meas}
 Let \cref{assumPDE,ass:RDE,assume3} be satisfied. Then, 
for any $n\geq 1$, system~\eqref{ODE},~\eqref{PDEb} is well-posed for any $y\in\RR^d$ and $\mu\in\ca{P}_{p}(\RR^d),\, p \geq \kappa$. Let $Y_t^n$ be the unique solution to~\eqref{ODE} and $\mu_t$ be the unique solution to~\eqref{PDEb} (in the sense of \cref{defofPDE}) with initial datum $(y,\mu)$. Then, the map $f:[0,\infty)\times \Omega \times \uE_p\to \uE_p$,
\begin{align}
f(t,\omega, (y,\mu))= & \,\,(f_1(t,\omega,(y,\mu)), f_2(t,\mu)) \coloneqq (Y_t^n,\mu_t)
\end{align}
is continuous. 
\end{proposition}
\vspace{0mm}
\begin{proof}
The function $f$ is continuous if and only if $f_1,f_2$ are continuous. It follows from~\cref{ThmWell-posednessPDE} (see~\eqref{eq:timeregmu} and~\eqref{PDE:Stability}) that $f_2$ is continuous in $t,\mu$. 

We now move on to $f_1$. Fix $n\geq 1$ and $\omega\in \Omega$. According to the definition of $W^n$, for $k=0,1,\ldots,2^{n}-1$, and $t\in[k\delta_n,(k+1)\delta_n)$,
 the ODE~\eqref{ODE} reduces to 
\begin{align}
\d Y_t^n=b(Y_t^n, \mu_t)\,\d t+\delta_n^{-1}\big(W_{(k+1)\delta_n}-W_{k\delta_n}\big)\sigma(Y_t^n, \mu_t)\,\d t-2^{-1}\big((\nabla\sigma)\sigma\big)(Y_t^n, \mu_t)\,\d t\,,
\end{align}
where $\delta_n=T/2^n$. Note that the coefficients $b$, $\sigma$ and $(\nabla \sigma)\sigma$ are all Lipschitz continuous in the state variable. It is clear that~\eqref{ODE} is well-posed and the solution $Y^n$ is continuous in time $t$.
Let $Y^n(\omega)$ and $\tilde{Y}^n(\tilde{\omega})$ be two solutions to (\ref{ODE}) 
started from initial values $(y,\mu_0)$, $(\tilde{y},\tilde{\mu}_0)$ and driven by two sample paths $\omega$ and $\tilde{\omega}$. For notational simplicity, we set $W^{n}(\omega)=\omega^n, W^{n}(\tilde{\omega})=\tilde{\omega}^n$. Then, it holds that 
\begin{align}
Y_t^n(\omega)-\tilde{Y}^n_t(\tilde{\omega})= \,\,\,&y-\tilde{y}+\int_0^t \Big[b(Y_s^n(\omega),\mu_s)-b(\tilde{Y}^n_s(\tilde{\omega}),\tilde{\mu}_s   )\Big] \;\d s\\
&+\int_0^t \Big[\sigma(Y_s^n(\omega),\mu_s)(\omega_s^n)'-\sigma(\tilde{Y}^n_s(\tilde{\omega}),\tilde{\mu}_s)(\tilde{\omega}_s^n)'\Big]\; \d s\\
&+2^{-1}\int_0^t \Big[\big((\nabla\sigma)\sigma\big)(Y_s^n(\omega),\mu_s)-\big((\nabla\sigma)\sigma\big)(\tilde{Y}^n_t(\tilde{\omega}),\tilde{\mu}_s)\Big]\;\d s\,,
\end{align}
where $(\mu_s)_{s\in [0,T]},(\tilde{\mu}_s)_{s\in [0,T]}$ are solutions to PDE (\ref{PDEb}) started from  $\mu_0, \tilde{\mu}_0 \in \ca{P}_p(\R^d)$ and the piecewise derivatives $(\omega_s^n)'=\delta_n^{-1}\big(\omega_{(k+1)\delta_n}-\omega_{k\delta_n}\big)$\,, $(\tilde{\omega}_s^n)'=\delta_n^{-1}\big(\tilde{\omega}_{(k+1)\delta_n}-\tilde{\omega}_{k\delta_n}\big)$ for any $k=0,1,\ldots,2^{n}-1$ and $s\in \big[k\delta_n, (k+1)\delta_n\big]$.  It is straightforward to see that for any $s\in [0,T]$, 
$(\omega_s^n)'\leq C_{n,T}\|\omega\|_0$, $(\tilde{\omega}_s^n)'\leq C_{n,T}\|\tilde{\omega}\|_0$ and $$|(\omega_s^n)'-(\tilde{\omega}_s^n)'|\leq C_{n,T}\|\omega-\tilde{\omega}\|_0,$$ which imply that 
\begin{align}
&|Y_t^n(\omega)-\tilde{Y}^n_t(\tilde{\omega})|\\
\leq\,\,&|y-\tilde{y}|+C_{n,T}\|\bar{\sigma}_{\tilde{\mu}}\|_{0,0}\|\omega-\tilde{\omega}\|_0\\[2mm]
&+K_1\int_0^t \left(1+|\tilde{Y}_s^n(\tilde{\omega})|\right)\Big[|Y_s^n(\omega)-\tilde{Y}^n_s(\tilde{\omega})|+\ca{D}_{p}(\mu_s,\tilde{\mu}_s) +G(\ca{D}_{p}(\mu_s,\tilde{\mu}_s))\Big]\,\d s\\[2mm]
\leq\,\,& |y-\tilde{y}|+C_{n,T}\|\bar{\sigma}_{\tilde{\mu}}\|_{0,0}\|\omega-\tilde{\omega}\|_0\\
&+K_1\int_0^t \left(1+|\tilde{Y}_s^n(\tilde{\omega})|\right)\Big[|Y_s^n(\omega)-\tilde{Y}^n_s(\tilde{\omega})|+\ca{D}_{p}(\mu_0,\tilde{\mu}_0) +G(C\ca{D}_{p}(\mu_0,\tilde{\mu}_0))\Big]\,\d s\\
\leq\,\,& K_2{\rm{e}}^{\int_0^t|\tilde{Y}_s(\tilde{\omega})|\d s}\Big[|y-\tilde{y}|+\|\omega-\tilde{\omega}\|_0+G(C\ca{D}_{p}(\mu_0,\tilde{\mu}_0))+\ca{D}_{p}(\mu_0,\nu_0)\Big]\,,
\end{align}
where $K_1$ depends on $n,\|\omega\|_0,M_{p, T}(\mu),M_{p,T}(\tilde{\mu})$, $K_2$ depends on $n, K_1, \|\bar{\sigma}_{\tilde{\mu}}\|_{0,0}, \sup_{t\in [0,T]}|\tilde{Y}_s(\tilde{\omega})|$, and we use  stability estimate~\eqref{PDE:Stability} in the second to last step. 
Since the distance $$d(\omega,\tilde{\omega})=\sum\limits_{i=1}^{\infty}\frac{1}{2^i}\frac{\|\omega-\tilde{\omega}\|_{[0,i]}}{1+\|\omega-\tilde{\omega}\|_{[0,i]}}\longrightarrow 0\,,$$ 
implies that $\|\omega-\tilde{\omega}\|_0\longrightarrow 0$, we obtain that  for fixed $n\geq 1$,
$$|Y_t^n(\omega)-\tilde{Y}_t^n(\tilde{\omega})|\longrightarrow 0\,,$$
provided that $|y-\tilde{y}|+\ca{D}_{p}(\mu_0, \tilde{\mu}_0)+d(\omega,\tilde{\omega})\longrightarrow 0$, which completes the proof.
\end{proof}
\vspace{3mm}
Now, we continue the proof of the main result ~\cref{exisrds}.
\begin{proof}[Proof of~\cref{exisrds}]
For any $0 \leq t_1\leq t_2 < \infty$,  $y\in\RR^d$, $\mu \in \mathcal{P}_{p}(\RR^d)$, denote by $S_{t_1,t_2} (\mu)$ the unique solution of PDE (\ref{drift-PDE}), in the sense of~\cref{defofPDE}, at time $t_2$ starting from time $t_1$ with  initial distribution $\mu$. We will write $S_{t_2}(\mu)$ whenever $t_1=0$.  Because the coefficients do not depend explicitly on time $t$, we have the identity 
$$S_{t_1,t_2}(\mu)=S_{t_2-t_1}(\mu)\,.$$
Let $Y_t$ be the unique solution  of RDE, in the sense of~\cref{def:RDEsol}, which exists due to~\cref{thm:MVRDE},
\begin{align}\label{intedrift-SDE}
Y_t=y+\int_0^tb(Y_r,S_{r}(\mu))\,\d r+\int_0^t\sigma(Y_r,S_{r}(\mu))\,\d\mathbf{W}_r\,,
\end{align}
which starts from $y \in \R^d$, driven by the path $\mathbf{W}(\omega)$. Define 
$$\phi_t(\omega,y,\mu)\coloneqq Y_t\,,$$
the unique solution to \eqref{intedrift-SDE}. We now consider at the joint solution map corresponding to the equations~\eqref{drift-SDE}  and  ~\eqref{drift-PDE}, which for any $t\geq0$, $\omega\in \Omega$, we define as 
\begin{align} 
&\varphi(t,\omega): \uE_p \longrightarrow \uE_p\,,\\
&\varphi(t,\omega):(y,\mu)
\mapsto  
       ( \phi_{t}(\omega,y,\mu), 
        S_{t}(\mu))\,.
    \label{defphi}
\end{align}
Clearly, $\varphi(0,\omega)=\mathrm{Id}_{\uE_p}$ for all $\omega \in \Omega$. We shall prove that $\varphi$ is a  RDS, i.e.\ $\varphi$ is $(\ca{B}([0,\infty))\otimes \ca{F} \otimes \ca{B}(\uE_p),\ca{B}(\uE_p))$-measurable, and  it satisfies the cocycle property, 
\begin{align}
\varphi(s+t,\omega)\ue=\varphi(t,\theta_s\omega)\circ\varphi(s,\omega)\ue\,,\;\;\;\;\;\omega\in \Omega,\, s,t\geq 0,\, \ue\in \uE_p\,.
\end{align}

The  $(\ca{B}([0,\infty))\otimes \ca{F} \otimes \ca{B}(\uE_p),\ca{B}(\uE_p))$-measurability of $\varphi$ follows from an approximation argument and~\cref{meas}. Indeed, the solution $(Y_{\cdot}^n,\mu_{\cdot})$ to~\eqref{ODE}  and  ~\eqref{PDEb} is $(\ca{B}([0,\infty))\otimes \ca{F} \otimes \ca{B}(\uE_p),\ca{B}(\uE_p))$-measurable. For any $t\geq 0$, due to the fact that 
 \begin{align}
|Y_{t}-Y_{t}^{n}|\leq CT^{\alpha}\big(\varrho_\alpha(\mathbf{W},\mathbf{W}^n)+\|Y,Y';Y^{n},(Y^{n})'\|_{W,W^{n},2\alpha}\big),
\end{align}
(see the estimate below \cite[(4.29)]{FrizHairer.2014.251}),   the fact that for all $\omega \in \Omega$, $\varrho_\alpha(\mathbf{W},\mathbf{W}^n) \to 0$ as $n \to \infty$, (see~\cref{setup}) and the stability estimate from~\eqref{sta-RDE}, it follows that for all $\omega \in \Omega$, $t\geq 0$, $y\in\R^d$ and $\mu \in \ca{P}_p(\R^d)$, 
\begin{align}
\lim_{n\to \infty}Y_t^n=\phi_{t}(\omega,y,\mu)
\end{align}
and thus, $\phi_t$ is $(\ca{F} \otimes \ca{B}(\uE_p),\ca{B}(\uE_p))$-measurable. Besides, since the mapping 
$$t\longmapsto \varphi(t)$$ is continuous, we conclude the $(\ca{B}([0,\infty))\otimes \ca{F} \otimes \ca{B}(\uE_p),\ca{B}(\uE_p))$-measurability of $\varphi$ by \cite[Lemma 4.51]{AliprantisBorder.2006.703}. We now move to prove the cocycle property of  $\varphi$, i.e.\ we need to verify
\begin{align}\label{cocycle}
\big( \phi_{s+t}(\omega,y,\mu),
S_{s+t}(\mu)\big)
=\big(\phi_{t}(\theta_s\omega, \phi_{s}(\omega,y,\mu),S_{s}(\mu)), 
S_{t}(S_{s}(\mu))\big)\,.
\end{align}
Applying the uniqueness of the PDE~\eqref{drift-PDE}, we can deduce  the semigroup property of $S$ in the sense that for any $s,t\geq 0$, $\mu\in\ca{P}_p(\R^d)$,
\begin{align}\label{cocycle2}
S_{s+t}(\mu)=S_{t}S_{s}(\mu) \, .
\end{align}
Due to the fact that
$\phi_{s+t}(y,\mu,\omega)$ is the unique solution of the RDE~\eqref{drift-SDE}, together with the definition of rough integral, we obtain that
\begin{align}\label{shift}
&\;\;\;\;\;\phi_{s+t}(\omega,y,\mu)\\
&=\,y+\int_0^{s+t}b(\phi_{r}(\omega,y,\mu),S_{r}(\mu))\, \d r+\int_0^{s+t}\sigma(\phi_{r}(\omega,y,\mu),S_{r}(\mu))\, \d \mathbf{W}_r(\omega)\\
&=\,y+\int_0^{s}b(\phi_{r}(\omega,y,\mu),S_{r}(\mu))\,\d r+\int_0^{s}\sigma(\phi_{r}(\omega,y,\mu),S_{r}(\mu))\,\d\mathbf{W}_r(\omega)\\
&\;\;\;\;\;\;\;\;+\int_s^{s+t}b(\phi_{r}(\omega,y,\mu),S_{r}(\mu))\, \d r+\int_s^{s+t}\sigma(\phi_{r}(\omega,y,\mu),S_{r}(\mu))\, \d\mathbf{W}_r(\omega)\\
&\,=\phi_{s}(\omega,y,\mu)+\int_s^{s+t}b(\phi_{r}(\omega,y,\mu),S_{r}(\mu))\, \d r\\
&\;\;\;\;\;\;\;\;\;\;\;\;\;\;\;\;\;\;\;\;+\lim_{|\ca{P}|\rightarrow0}\sum_{[a_i, a_{i+1}]\in \ca{P}}\sigma(\phi_{a_i}(\omega,y,\mu),S_{a_i}(\mu))(W_{a_{i+1}}-W_{a_i})(\omega)\\
&\;\;\;\;\;\;\;\;\;\;\;\;\;\;\;\;\;\;\;\,\,+\lim_{|\ca{P}|\rightarrow0}\sum_{[a_i, a_{i+1}]\in \ca{P}}(\nabla\sigma)(\phi_{a_i}(\omega,y,\mu),S_{a_i}(\mu))\phi_{a_i}'(\omega,y,\mu)\bb{W}_{a_i,a_{i+1}}(\omega)\\
&\,=\phi_{s}(\omega,y,\mu)+\int_0^{t}b(\phi_{s+r}(\omega,y,\mu),S_{s+r}(\mu))\, \d r\\
&\;\;\;\;\;\;\;\;\;\;\;\;\;\;\;\;\;\;\;\;+\lim_{|\tilde{\ca{P}}|\rightarrow0}\sum_{[\tilde{a}_i, \tilde{a}_{i+1}]\in \tilde{\ca{P}}}\sigma\big(\phi_{s+\tilde{a}_i}(\omega,y,\mu),S_{\tilde{a}_i}(S_s(\mu))\big)(W_{\tilde{a}_{i+1}+s}-W_{\tilde{a}_{i}+s})(\omega)\\
&\;\;\;\;\;\;\;\;\;\;\;\;\;\;\;\;\;\;\;\;\,+\lim_{|\tilde{\ca{P}}|\rightarrow0}\sum_{[\tilde{a}_i, \tilde{a}_{i+1}]\in \tilde{\ca{P}}}(\nabla\sigma)\big(\phi_{s+\tilde{a}_i}(\omega,y,\mu),S_{\tilde{a}_i}(S_s(\mu))\big)\phi_{s+\tilde{a}_i}'(\omega,y,\mu)\bb{W}_{\tilde{a}_i+s,\tilde{a}_{i+1}+s}(\omega)\\
&\,=\phi_{s}(\omega,y,\mu)+\int_0^{t}b(\phi_{s+r}(\omega,y,\mu),S_{s+r}(\mu))\, \d r\\
&\;\;\;\;\;\;\;\;\;\;\;\;\;\;\;\;\;\;\;+\lim_{|\tilde{\ca{P}}|\rightarrow0}\sum_{[\tilde{a}_i, \tilde{a}_{i+1}]\in \ca{P}}\sigma\big(\phi_{s+\tilde{a}_i}(\omega,y,\mu),S_{\tilde{a}_i}(S_s(\mu))\big)\big(W_{\tilde{a}_{i+1}}(\theta_s\omega)-W_{\tilde{a}_{i}}(\theta_s\omega)\big)\\
&\;\;\;\;\;\;\;\;\;\;\;\;\;\;\;\;\;\;\;+\lim_{|\tilde{\ca{P}}|\rightarrow0}\sum_{[\tilde{a}_i, \tilde{a}_{i+1}]\in \ca{P}}(\nabla\sigma)\big(\phi_{s+\tilde{a}_i}(\omega,y,\mu),S_{\tilde{a}_i}(S_s(\mu))\big)\phi_{s+\tilde{a}_i}'(\omega,y,\mu)\bb{W}_{\tilde{a}_i,\tilde{a}_{i+1}}(\theta_s\omega)\\
&=\phi_{s}(\omega,y,\mu)+\int_{0}^tb\big(\phi_{s+r}(\omega,y,\mu),S_{r}(S_{s}(\mu))\big)\,\d r\\
&\;\;\;\;\;\;\;\;\;\;\;\;\;\;\;\;\;\;+\int_0^{t}\sigma\big(\phi_{s+r}(\omega,y,\mu),S_{r}(S_{s}(\mu))\big)\,\d \mathbf{W}_r(\theta_{s}\,\omega)\,,
\end{align}
where $\ca{P}=\{a_0, a_1,...\}$ denotes a partition of $[s, s+t]$  and $\tilde{\ca{P}}=\{\tilde{a}_0, \tilde{a}_1,...\}=\{a_0-s, a_1-s,...\}$ is the corresponding shifted partition of $[0,t]$, $|\ca{P}|$ denotes the length of the largest element of $\ca{P}$. For each $i\in \{0,\ldots,n-1\}$, $\tilde{a}_i=a_i-s$, and we used the identity~\eqref{eq:roughshift} in the last to the second equality. Here, the rough integral is justified, since $(\phi_{s+{.}}(\omega,y,\mu),\phi'_{0,s+{.}}(\omega,y,\mu))\in \scr{D}_{\theta_sW}^{2\alpha}([0,t];\RR^d)$ provided that $(\phi_{{.}}(\omega,y,\mu),\phi'_{0,{.}}(\omega,y,\mu))\in \scr{D}_W^{2\alpha}([s,s+t];\RR^d)$.  Due to the path-by-path uniqueness of the RDE (\ref{drift-SDE}) (or~\eqref{intedrift-SDE}), we know that 
\begin{align}
\phi_{s+t}(\omega,y,\mu)=\phi_{t}(\theta_s\omega,\phi_{s}(\omega,y,\mu),S_{s}(\mu))\,.
\end{align}
Therefore, $\varphi$ defined in (\ref{defphi}) generates an RDS. The time continuity of the RDS follows from the continuity in time of the trajectory of the solution to~\eqref{drift-SDE} and the estimate~\eqref{eq:timeregmu}
  of the  solution to~\eqref{drift-PDE} along with
  \begin{align}
  	|\varphi(t, \omega)(y,\mu)-\varphi(s, \omega)(y,\mu)|\leq |Y_t(\omega)-Y_s(\omega)|+\ca{D}_{p}(\mu_t,\mu_s)\,.
  \end{align}

\end{proof}

\section{Applications}\label{sec-appli}
\subsection{Ensemble Kalman sampler}
In this section, we check that~\cref{assumPDE,ass:RDE,assume3} apply to the ensemble Kalman system~\eqref{eksdecoupled} discussed in~\cref{ssub:the_ensemble_kalman_sampler} and thus obtain the existence of an RDS. In this example, we define 
\begin{align}\label{eq:Cmu}
	 C(\mu)\coloneqq {\rm{Cov}}(\mu)=\int_{\RR^d}\big(y-m(\mu)\big)\otimes\big(y-m(\mu)\big)\;\d\mu\,,
\end{align}
where $m(\mu)$ is the mean as defined in~\eqref{eq:mean}. $C (\mu)$ can be expressed component-wise as follows
$$C_{ij}(\mu)=\int_{\RR^d}y_i y_j \, \mathrm{d}\mu-\left(\int_{\R^d}y_i \, \rm{d}\mu\right)\left(\int_{\R^d}y_j \, \rm{d}\mu\right).$$

We now make the following assumptions on the prior we wish to sample from.
\begin{assumption}\label{ass:V}
	$V: \RR^d\to \RR$ is three times continuously differentiable, and there exists a constant $C$ such that for all $x, y\in\RR^d$,
\begin{align}
		-\langle C(\mu)\nabla V(y), y\rangle &\leq \, C(1+|y|^2+M_2(\mu))\,,\\
		|\nabla V(0)|+|D^2 V(y)|&+|D^3 V(y)|\leq C\,.
	\end{align} 

where $C(\mu)$ is as in~\eqref{eq:Cmu}.
	\end{assumption}

\begin{proposition}\label{eks} Under \cref{ass:V},
there exists an RDS on the product space $\uE_p=\RR^d\times \ca{P}_p(\RR^d)$ for all $p\geq 2$ associated with the equations
\begin{align}
\d Y_t&=-C(\mu_t)\nabla V(Y_t)\,\d t+\sqrt{C(\mu_t)}\,\d\mathbf{W}_t\,,\;\;\;\;\;\;\;\;\;\;\;\;\;\;\;Y_0=y\,,\label{EnsembleSDE}\\
\partial_t \mu_t&=\nabla \cdot\big(C(\mu_t) \nabla V(\cdot)\mu_t\big)+\frac{1}{2}D^2:(C(\mu_t)\mu_t)\,,\;\;\;\;\;\;\;\mu_0=\mu\,,\label{EnsemblePDE}
\end{align}
where $(y,\mu)\in \uE_p$\,.
\end{proposition}

\begin{proof}
We just need to verify that~\cref{assumPDE,ass:RDE,assume3} are satisfied  with $\kappa=2$ under~\cref{ass:V}. We start by observing that for any $\mu\in\ca{P}_2(\RR^d)$, the  Hilbert--Schmidt norm  can be bounded as  follows
\begin{equation}\label{CmuHS}
	|C(\mu)|\leq 2d M_2(\mu)\,.
\end{equation}
 Indeed, by the definition of Hilbert--Schmidt norm, for any $\mu\in\ca{P}_2(\RR^d)$, 
\begin{align}
|C(\mu)|^2
&=\,\sum_{i,j=1}^{d}\bigg|\int_{\RR^d} y_iy_j\,\d\mu-\int_{\RR^d} y_i\,\d \mu\int_{\RR^d} y_j\;\d\mu\bigg|^2\\
&\leq \,2\sum_{i,j=1}^{d}\left(\int_{\RR^d}|y|^2\,\d\mu\right)^2+2\sum_{i,j=1}^{d}\left(\int_{\RR^d} |y|^2\,\d\mu\right)\left(\int_{\RR^d} |y|^2\,\d\mu\right)\\
&= \,4 d^2 M_2^2(\mu).
\end{align}
Moreover, for any $\mu,\nu\in\ca{P}_2(\RR^d)$,  
\begin{equation}
	|C(\mu)-C(\nu)|\leq d(1+M_1(\mu)+M_1(\nu) )\ca{D}_2(\mu,\nu)\, ,
\end{equation}
and 
\begin{equation}
	\big|\sqrt{C(\mu)}-\sqrt{C(\nu)}\,\big|\leq d\sqrt{(1+M_1(\mu)+M_1(\nu) )\ca{D}_2(\mu,\nu) }\, .
\end{equation}
To see this, for any $i, j\in \{1,\ldots,d\}$, we compute 
\begin{align}
&|C_{ij}(\mu)-C_{ij}(\nu)|\\
=\,\,&\bigg|\int_{\RR^d} y_iy_j\,\d(\mu-\nu)-\int_{\RR^d} y_i\,\d\mu\int_{\RR^d} y_j \,\d\mu+\int_{\RR^d} y_i\,\d\nu\int_{\RR^d} y_j\,\d\nu\bigg|\\
\leq\,\, &\ca{D}_2(\mu,\nu)+\bigg|\left(\int_{\RR^d} y_i\,\d(\mu-\nu)\right)\int_{\RR^d} y_j\,\d\mu\bigg|+\bigg|\int_{\RR^d} y_i\,\d \nu\left(\int_{\RR^d} y_j\,\d(\mu-\nu)\right)\bigg|\\
\leq \,\,& (1+M_1(\mu)+M_1(\nu) )\ca{D}_2(\mu,\nu)\,,
\end{align}
which leads that
\begin{align}
|C(\mu)-C(\nu)|^2=&\sum_{i,j=1}^{d}|C_{ij}(\mu)-C_{ij}(\nu)|^2
\leq d^2(1+M_1(\mu)+M_1(\nu))^2\ca{D}^2_2(\mu,\nu)\,.
\end{align}
Moreover, let $\ca{C}^{ij}(\mu,\nu)$ be a matrix for which 
\begin{align}
(\ca{C}^{ij}(\mu,\nu))_{kl}=
\begin{cases}
&C_{ij}(\mu)-C_{ij}(\nu)\;\;\text{if}\;\;(k, l)=(i,j)\,,\\ 
 &0, \;\;\text{otherwise}.
\end{cases}
\end{align}
Using \cite[Lemma 4.1]{PowersSto.1970.CMP1} and the sub-additivity of trace norm, we derive that
\begin{align}\label{sqrtmu}
\big|\sqrt{C(\mu)}-\sqrt{C(\nu)}\,\big|^2\leq\,& \,\|C(\mu)-C(\nu)\|_{{\rm{Tr}}}\\
\leq\,& \sum_{i,j=1}^d\|\ca{C}^{ij}(\mu,\nu)\|_{{\rm{Tr}}}\\
\leq\,& \sum_{i,j=1}^d {\rm{Tr}}\Big[\big(\ca{C}^{ij}(\mu,\nu)\big)^*\ca{C}^{ij}(\mu,\nu)\Big]^{\frac{1}{2}}\\
\leq\, &\,d^2 (1+M_1(\mu)+M_1(\nu) )\ca{D}_2(\mu,\nu)\,.
\end{align}
\noindent Hence, one is readily to check that
$$b(y,\mu)=-C(\mu)\nabla V(y)\,,\;\;\;\sigma(\mu)=\sqrt{C(\mu)}\,,$$
satisfy
\cref{assumPDE} and \cref{assume3} for the choice of $\kappa=2$. We are left to verify  regularity requirements of \cref{ass:RDE}. Because in the ensemble Kalman sampler model, the diffusion term  is state-independent, it suffices to show that for each $i,j=1,\ldots,d$, 
\begin{align}\label{Cijmu}
C_{ij}(\mu_t)=\int_{\RR^d}y_iy_j\,\d\mu_t-\int_{\RR^d}y_i\,\d\mu_t\int_{\RR^d}y_j\,\d\mu_t
\end{align}
is once differentiable with respect to time $t$. Recall that $\mu_t$ solves the evolution PDE
\begin{align}\label{recallensembPDE}
\partial_t \mu_t&=\nabla \cdot\big(C(\mu_t)\nabla V(\cdot)\mu_t\big)+\frac{1}{2}D^2:(C(\mu_t)\mu_t) \,,
\end{align}
in the sense of~\cref{defofPDE}. We choose the admissible test function $\phi(y)=y_iy_j\in\ca{K}_2$, for some fixed $i,j\in \{1,\ldots,d\}$, to obtain
for any $t_0\in [0,T]$,
\begin{align}\label{es:Lipoft}
\begin{split}
	&\int_{\RR^d}y_iy_j\,\d\mu_t-\int_{\RR^d}y_iy_j\,\d\mu_{t_0}\\
=\,\,\,&\int_{t_0}^t\int_{\RR^d}\langle \nabla\big(y_iy_j
\big), -C(\mu_s)\nabla V(y)\rangle\, \d\mu_s\d s+\int_{t_0}^t\int_{\RR^d}{\rm{Tr}}\big(C(\mu_s)D^2(y_iy_j) \big)\,\d\mu_s\d s\\
=\,\,\,&\int_{t_0}^t \int_{\RR^d}-\Big(y_j\cdot \sum_{k=1}^d C_{ik}(\mu_s)\partial_{y_k} V(y)+y_i\cdot \sum_{k=1}^d C_{jk}(\mu_s)\partial_{y_k} V(y)\Big)\, \d \mu_s\d s\\
&+\int_{t_0}^t \int_{\RR^d} \big[C_{ij}(\mu_s)+C_{ji}(\mu_s)\big]\,\d\mu_s\d s\\
\leq\,\,\, &\sum_{k=1}^d \int_{t_0}^t \int_{\RR^d}
\big[|C_{ik}(\mu_s)||y_j|(1+|y|)+|C_{jk}(\mu_s)||y_i|(1+|y|)\big]\,\d \mu_s\d s+2\int_{t_0}^t  C_{ij}(\mu_s)\,\d s\\
\leq\,\,\, &C_d\Big(1+\sup_{t\in [0,T]}M_2(\mu_t)\Big)^2|t-t_0|\,,
\end{split}
\end{align}
where we used the estimate~\eqref{CmuHS} in the last inequality. Therefore, for each $i,j\in \{1,2,\ldots,d\}$, $C_{ij}(\mu_{\cdot})$ is once differentiable in time $t$ and the derivative  is bounded by $(1+\sup_{t\in [0,T]}M_2^2(\mu_t))^2$ by the fundamental theorem of calculus.
\end{proof}
\subsection{Lagrangian formulation of the Landau equation}

In this section, we check that~\cref{assumPDE,ass:RDE,assume3} apply to the system~\eqref{eq:landaudecoupled}. We have the following result:
\begin{proposition}\label{landauprop}
	Consider the system~\eqref{eq:landaudecoupled} with $\gamma=0$. Then, there exists an RDS on the product space $\uE_p=\RR^d\times \ca{P}_p(\RR^d)$, $p\geq 2$, associated to~\eqref{eq:landaudecoupled}.
	\begin{proof}
		One can readily see that the distribution-dependent SDE~\eqref{eq:Landau1} satisfies \cref{assumPDE} and \cref{assume3}. It remains to justify that the coefficients of the RDE~\eqref{eq:Landau1} fulfill \cref{ass:RDE}, i.e.\ $\|m(\mu_t)\|_0$ is finite and $\sigma_0(\mu_t)$ is continuously differentiable, which is  straightforward using a similar argument as in~\eqref{es:Lipoft}. Therefore, it follows from \cref{exisrds} that an RDS exists associated with the system~\eqref{eq:landaudecoupled}.   
	\end{proof}
		\end{proposition}

\appendix
\section{Proof of~\texorpdfstring{\cref{convergenceequiva}}{Theorem~\ref{convergenceequiva}}}\label{sec:proofconvergenceequiva}

\begin{proof} (i) One can readily  see that $\ca{D}_p$ is symmetric and satisfies the triangle inequality. Let $\eta^{\varepsilon}$ be a standard mollifier. Suppose $\ca{D}_p(\mu,\nu)=0$,  we see that for any bounded continuous function $f$,
\begin{align}
\begin{split}\label{metric}
\int_{\RR^d} f \,\d (\mu-\nu)&=\int_{\RR^d} \big(f-f*\eta^{\varepsilon} \big)\,\d(\mu-\nu)+\int_{\RR^d}f*\eta^{\varepsilon} \,\d (\mu-\nu)\\
&=\int_{\RR^d} \big(f-f*\eta^{\varepsilon} \big)\,\d (\mu-\nu)\,,
\end{split}
\end{align}
where we used the fact that the integrand in second integral  lies in the class $\ca{K}_p$. Taking $\varepsilon\to 0$ in~\eqref{metric} and using the dominated convergence theorem, we derive that for any bounded continuous function $f$,
$$\int_{\RR^d} f \,\d (\mu-\nu)=0\,,$$
which implies that $\mu=\nu$. We show separability and completeness of $\ca{D}_p$ after we complete the proof of (ii).

\vspace{2mm}
(ii) $\Longrightarrow$ Let us assume that $\mu_n\xrightarrow{\ca{D}_p}\mu$ as $n\to \infty$. It is straightforward to check that the function $f(x)\coloneqq \frac{1}{p}|x|^p\in \ca{K}_p$ from which it follows that
\begin{equation}\label{eq:xsquareconvergenc}
	\int_{\RR^d} |x|^p\,\d \mu_n\longrightarrow\int_{\RR^d} |x|^p\,\d \mu\,,\,\,\,\,\,\,\,\,\,\,n\to \infty\,,
\end{equation} 
and for any $\varphi\in C_b^2(\R^d)\subseteq \ca{K}_p$ (up to a multiplicative constant) \footnote{Here $C_b^2(\R^d)$ is the family of twice continuously differentiable functions with bounded derivatives up to order 2.} we see that
\begin{equation}\label{eq:Cb2convergence}
	\int_{\RR^d} \varphi(x)\,\d \mu_n\longrightarrow\int_{\RR^d} \varphi(x)\,\d \mu\,,\,\,\,\,\,\,\,\,\,\,n\to \infty\,.
\end{equation} 
It follows from~\eqref{eq:xsquareconvergenc} that
$\sup_n M_p(\mu_n)<\infty$. Thus,  there exists
a subsequence $(n_k)_{k\geq 1}$ and $\mu^*\in\ca{P}(\RR^d)$ such that for any bounded and continuous $\varphi$,
\begin{equation}
	\int_{\RR^d} \varphi(x)\,\d \mu_{n_k}\longrightarrow\int_{\RR^d} \varphi(x)\,\d \mu^*\,,\,\,\,\,\,\,\,\,\,\,n\to \infty\,.
\end{equation}
Because of~\eqref{eq:Cb2convergence} and  the uniqueness of the limit, we must have that $\mu^*=\mu$. In combination with~\eqref{eq:xsquareconvergenc}, we have that $\mu_n\xrightarrow{d_p} \mu$ using the characterisation of $d_p$-convergence provided in  \cite[Definition 6.8 (iv), Theorem 6.9]{Villani.2009.973}.

$\Longleftarrow$  We now assume $\mu_n\xrightarrow{d_p}\mu$. Let us also define 
\begin{align}
\overline{\ca{D}}(\mu,\nu)=\sup_{
 |\varphi|\leq 1, |\nabla \varphi|\leq 1}\int_{\RR^d}\varphi \,\d (\mu-\nu)\,,
\end{align}
which is also a metric (using   a similar argument as above), in the space $\ca{P}_p(\RR^d)$. Using the Kantorovich--Rubinstein formulation of $d_1$ (see \cite[Remark 6.5]{Villani.2009.973}), we see that for any $\mu, \nu\in\ca{P}_p(\RR^d)$,
\begin{align}\label{W1comparison}
	\overline{\ca{D}}(\mu,\nu)\leq  \sup_{
\varphi:
 |\nabla\varphi|\leq 1} \int_{\RR^d}\varphi \,\d(\mu-\nu)=d_1(\mu,\nu)\,.
\end{align}
In particular, as $n\to \infty$, if $\mu_n\xrightarrow{d_p} \mu$, then $\mu_n\xrightarrow{\overline{\ca{D}}} \mu$. 
 Let $\overline{\ca{D}}(\mu_n,\mu)\to 0$ and $\psi_{\delta}$ be smooth, compactly supported on a ball of size $2/{\delta}$ and be equal to 1 on  the ball of size $1/{\delta}$. For any $\varphi\in\ca{K}_p$, we assume $\varphi(0)=0$ without loss of generality  
 and approximate it by $\varphi \psi_{\delta}$. Then 
\begin{align}
\ca{D}_p(\mu_n,\mu)\leq &\,\sup_{\varphi\in\ca{K}_p}\int_{\RR^d} \varphi(1-\psi_{\delta})\,\d(\mu_n-\mu)+\sup_{\varphi\in\ca{K}_p}\int_{\RR^d} \varphi\psi_{\delta}\,\d(\mu_n-\mu)\\
\leq &\;\sup_{\varphi\in\ca{K}_p}\int_{B^c_{1/\delta}} \varphi(1-\psi_{\delta})\,\d(\mu_n-\mu)+\sup_{\varphi\in\ca{K}_p} \max\{\|\varphi\psi_{\delta}\|_{0}, \|\nabla (\varphi\psi_{\delta})\|_{0}\}\overline{\ca{D}}(\mu_n,\mu)\\
\leq &\;\int_{B^c_{1/\delta}} 2\big(1+|y|^p\big)\,\d(\mu_n+\mu)\\
&\,+\sup_{\varphi\in\ca{K}_p}\max\{\|\varphi\psi_{\delta}\|_{0}, \|\nabla (\varphi\psi_{\delta})\|_{0}\} \overline{\ca{D}}(\mu_n,\mu)\\
\leq &\;\int_{B^c_{1/\delta}} 2\big(1+|y|^p\big)\,\d(\mu_n+\mu)\\
&\,+\max\Big\{2\Big(1+\frac{2^p}{\delta^p}\Big), \Big(1+\Big(\frac{2}{\delta}\Big)^{p-1}\Big)\|\psi_{\delta}\|_0+2\Big(1+\frac{2^p}{\delta^p}\Big)\|\nabla \psi_{\delta}\|_0\Big\} \overline{\ca{D}}(\mu_n,\mu)\,.\label{barD}
\end{align}

Using $\mu_n\xrightarrow{d_p}\mu$, we have that for any bounded and continuous $\varphi$,
 \begin{equation}
 	\int_{\RR^d}\varphi(x)|x|^p\,\d \mu_n\to \int_{\RR^d}\varphi(x)|x|^p\,\d \mu\,.
 \end{equation}
 Consider the sequence of positive measures $\d \nu_n\coloneqq |x|^p\,\d \mu_n$, which converge weakly to $\d \nu\coloneqq |x|^p\,\d \mu$. Using Prokhorov's theorem, we obtain that $\nu_n$ is tight and thus, for any $\eps>0$, we may pick $\delta$ small enough such that the first term on the right hand side of~\eqref{barD} is smaller than $\eps$ for all $n\geq 1$. We can then pick $n$ large enough so that the second term on the right hand side of~\eqref{barD} is smaller than $\eps$. Since $\eps>0$ is arbitrary it follows that $\mu_n\xrightarrow{\ca{D}_p} \mu$. Since $d_p$ is complete and separable so is $\ca{D}_p$.

  
 \end{proof}

\section{Moment estimates}\label{Moment estimates}
\label{uniformappendix}
\begin{lemma}
Let $Y^n$ be the unique solution to the approximating SDE~\eqref{approxiSDE}, then we have  
\begin{align}\label{uniformlymomentestimateappendix}
\sup_{n\geq 1}\mathbb{E}\bigg[\sup_{t\in[0,T]}|Y^{n}_t|^{p}\bigg]\leq C_{p,T}\,,
\end{align}
and for any 
 $s, t\in [0,T]$,
\begin{align}
\sup_{n\geq 1}\EE|Y_t^{n}-Y_s^{n}|^{p}&\leq C_{T,  \kappa,p} |t-s|^{\frac{p}{2}}\,.
\end{align}
\end{lemma}
\begin{proof}
Applying It\^{o}'s formula, we have
\begin{align}
|Y^{n}_t|^2=\;&|Y_0|^2+2 \int_0^t  \langle Y^{n}_s, b(Y^{n}_s,\mu^{n}_{\delta^n_s})\rangle \,\d s+2 \int_0^t \langle Y^{n}_s, \sigma(Y^{n}_s,\mu^{n}_{\delta^n_s}) \;\d W_s\rangle\\
&+\int_0^t |\sigma(Y^{n}_s,\mu^{n}_{\delta^n_s}) |^2\,\d s\,,
\end{align}
and using It\^{o}'s formula again, then
\begin{align}
|Y^{n}_t|^p=\;\;\;&|Y_0|^p+\frac{p(p-2)}{2}   \int_0^t|Y^{n}_s|^{p-4}|\sigma^*(Y^{n}_s,\mu^{n}_{\delta^n_s})Y^{n}_s|^2\,\d s\\
&+\frac{p}{2}\int_0^t|Y^{n}_s|^{p-2}\big(2\langle Y^{n}_s, b(Y^{n}_s,\mu^{n}_{\delta^n_s})\rangle
 +|\sigma(Y^{n}_s,\mu^{n}_{\delta^n_s})|^2\big)\,\d s
 \\
&+p\int_0^t|Y^{n}_s|^{p-2}\langle Y^{n}_s,\sigma(Y^{n}_s,\mu^{n}_{\delta^n_s})\,\d W_s\rangle
 \\
 \leq\;\;\;&|Y_0|^p+C_p\int_0^t |Y^{n}_s|^{p-2}\big(1+|Y^{n}_s|^2+(\mathbb{E}|Y^{n}_{\delta^n_s}|^{\kappa})^{\frac{2}{\kappa}}\big)\,\d s
  \\
 &+p\int_0^t|Y^{n}_s|^{p-2}\langle Y^{n}_s,\sigma(Y^{n}_s,\mu^{n}_{\delta^n_s})\,\d W_s\rangle
 \\
 \leq\;\;\;&|Y_0|^p+C_p\int_0^t (1+|Y^{n}_s|^p+\mathbb{E}|Y^{n}_{\delta^n_s}|^p)\,\d s
  \\
&+p\int_0^t|Y^{n}_s|^{p-2}\langle Y^{n}_s,\sigma(Y^{n}_s,\mu^{n}_{\delta^n_s})\,\d W_s\rangle\,,
\end{align}
where we use the assumptions (\ref{Weakcoercivity}) and (\ref{Growthsigma}) in the first inequality, Jessen's inequality and  Young's inequality in the last step. By Burkholder-Davis-Gundy's inequality, it holds that
\begin{align}
&\mathbb{E}\Big[\sup_{t\in[0,T]}|Y^{n}_t|^p\Big]\\
\leq \;\;\;&\mathbb{E}|Y_0|^p+C_p \int_0^T1+\mathbb{E}\Big[\sup_{u\in[0,s]}|Y^{n}_u|^p\Big]\,\d s+C_p\EE\Big(\int_0^T |Y^{n}_s|^{2p-2}\cdot |\sigma(Y^{n}_s,\mu^{n}_{\delta^n_s})|^2\,\d s\Big)^{\frac{1}{2}}\\
\leq \;\;\;&\mathbb{E}|Y_0|^p+C_p \int_0^T 1+\mathbb{E}\Big[\sup_{u\in[0,s]}|Y^{n}_u|^p\Big]\;\d s\\
&+C_p\EE\bigg[\sup_{s\in [0,T]} |Y^{n}_s|^p\cdot  \int_0^T |Y^{n}_s|^{p-2}  \big(1+|Y^{n}_s|^2+(\mathbb{E}|Y^{n}_{\delta^n_s}|^{\kappa})^{\frac{2}{\kappa}}\big) \; 
 \d s\bigg]^{\frac{1}{2}}\\
\leq \;\;\;&\mathbb{E}|Y_0|^p+C_p \int_0^T 1+\mathbb{E}\Big[\sup_{u\in[0,s]}|Y^{n}_u|^p\Big]\;\d s\\
&+\frac{1}{2}\mathbb{E}\bigg[\sup_{t\in[0,T]}|Y^{n}_t|^p\bigg]+C_p \EE \int_0^T (1+|Y^{n}_s|^p+\mathbb{E}|Y^{n}_{\delta^n_s})|^p)  
 \,\d s\\
 \leq \;\;\;&\mathbb{E}|Y_0|^p+C_p \int_0^T1+\mathbb{E}\Big[\sup_{u\in[0,s]}|Y^{n}_u|^p\Big]\,\d s+\frac{1}{2}\mathbb{E}\Big[\sup_{t\in[0,T]}|Y^{n}_t|^p\Big]\,,
\end{align}
which implies the uniform estimate~\eqref{uniformlymomentestimateappendix} by using Gr\"{o}nwall's lemma. For the sake of convenience,  we omit the procedure of defining the stopping times $\tau_k\coloneqq\inf\{t\geq 0, |Y_t^n|\geq k\}$ to ensure the finiteness of $\EE\Big[\sup_{t\in[0,T\wedge \tau_k]}|Y^{n}_t|^p\Big]$ before using Gr\"{o}nwall's lemma.

Moreover, because of  the growth conditions~\eqref{Growthsigma},~\eqref{Growthb} of $b, \sigma$, using H$\rm{\ddot{o}}$lder's inequality and Burkholder-Davis-Gundy's inequality, it holds that
\begin{align}
&\mathbb{E}|Y^{n}_t-Y^{n}_s|^{p}\\
\leq\,\,\, &|t-s|^{p-1}\mathbb{E}\int_{s}^{t} |b(Y^{n}_u,\mu^{n}_{\delta^n_u})|^{p} \,\d u+\mathbb{E}\Big|\int_{s}^{t}\sigma(Y^{n}_u,\mu^{n}_{\delta^n_u})\,\d W_u\Big|^{p}\\
\leq\,\,\, &|t-s|^{p-1}\mathbb{E}\int_{s}^{t} |b(Y^{n}_u,\mu^{n}_{\delta^n_u})|^{p} \,\d u+|t-s|^{\frac{p}{2}-1}\mathbb{E}\int_{s}^{t}\big|\sigma(Y^{n}_u,\mu^{n}_{\delta^n_u})\big|^p\,\d u\\
\leq\,\,\,  &C_{T,\kappa,p}|t-s|^{p-1}\int_{s}^{t}1+\bb{E}\Big[\sup_{u\in[0,T]}|Y_u^n|\Big]^p\,\d u
 \\
&+C|t-s|^{\frac{p}{2}-1}\int_{s}^{t}1+\mathbb{E}\Big[\sup_{u\in[0,T]}|Y^{n}_u|^{p}\Big]\,\d u\\
\leq ~&C_{T,\kappa, p}|t-s|^{\frac{p}{2}}\,.
\end{align}
\end{proof}

\section{Differentiability with respect to the initial condition}\label{Differentiability of backward PDE}
\begin{lemma}
Let $\xi_T^{t, y}$ be the unique solution to SDE~\eqref{FeynmanSDE} at time $T$, and $\varphi\in\ca{K}_{\kappa}$. Define
$$u(t, y)\coloneqq \EE\big[\varphi(\xi_T^{t, y})\big]\,.$$
Then $u$ is twice continuously differentiable with respect to $y$ with the derivatives satisfying
\begin{align}
(\partial_{y_i} u)(t,y)
&=\sum_{\ell=1}^{d}\EE \Big[(\nabla \varphi)_{\ell}(\xi_T^{t,y}) \eta_T^{\ell i}(t,y)\Big]\,,\label{firstderivaof}\\
(\partial_{y_iy_j} u)(t,y)&=\sum_{\ell,m=1}^{d} \EE\Big[(D^2 \varphi)_{i\ell}(\xi_T^{t,y}) \eta^{\ell m}_T(t,y) \eta^{mj}_T(t,y)\Big]+\sum_{n=1}^d\EE\big[(\nabla \varphi)_{n}(\xi_T^{t,y}) \gamma^{nij}_T(t,y)\big]\,,\label{secondderivaof}
\end{align}
where
\begin{equation}
    \eta_T^{ij}\coloneqq \partial_{y_j}(\xi_T^{t,y})_i\,\quad \gamma_T^{ijk}\coloneqq \partial_{y_k}(\eta_T^{ij}) \, .
\end{equation}
As a consequence, we have for any $t\in [0,T]$,
	\begin{align}
|\nabla u| (t, y)&\leq C_1(1+|y|^{\kappa-1})\,,\\
|D^2 u| (t, y)&\leq C_2(1+|y|^{\kappa-1})\,,
\end{align}
where the constants $C_1, C_2$ depend on $d, T, M_{\kappa, T}(\mu)$.
\end{lemma}
\begin{proof} 
It follows from \cite[V.7. Lemma 1, Remark 5]{Krylov.1995.271} that the derivatives are given by the expressions in~\eqref{firstderivaof} and~\eqref{secondderivaof}. Moreover, component-wise  $\eta$ and $\gamma$ solve the following SDEs:
\begin{align}
&\eta^{\ell i}_T(t, y)\\
=\,\,&\delta_{\ell i}+\sum_{n=1}^{d}\int_t^T \partial_{y_n} b_\ell(\xi_r^{t,y},\mu_r) \eta^{ni}_r(t,y)\,\d r+\sum_{m,n=1}^d\int_t^T\partial_{y_n} \sigma_{\ell m}(\xi_r^{t,y},\mu_r) \eta^{ni}_r(t, y) \,\d W_r^m\,,\label{first}\\
&\gamma^{\ell ij}_{T}(t, y)\\
=\,\,&\sum_{k,n=1}^d
\int_t^T \partial_{y_ky_n}b_\ell(\xi_r^{t,y},\mu_r) \eta^{kj}_r(t,y) \eta^{ni}_r(t, y)\,\d r+\sum_{n=1}^d
\int_t^T \partial_{y_n}b_\ell(\xi_r^{t,y},\mu_r) (\partial_{y_j} \eta^{ni}_r)(t, y)\,\d r\\
&+\sum_{k,m,n=1}^{d}\int_t^T \partial_{y_ky_n}\sigma_{\ell m}(\xi_r^{t, y},\mu_r)\eta^{kj}_r(t,y)     \eta^{ni}_r(t, y)\, \d W_r^m\\
&+\sum_{m,n=1}^{d}\int_t^T\partial_{y_n}\sigma_{\ell m}(\xi_r^{t,y},\mu_r) (\partial_{y_j}\eta^{ni}_r)(t, y)\,\d W_r^m.\label{second}
\end{align}
In what follows we will use $C_{d,T,M_{\kappa, T}(\mu)}$ to denote a generic constant whose actual value may change from line to line. Now, taking square and expectation on both sides of~\eqref{first}, by H$\rm{\ddot{o}}$lder's inequality and It\^{o}'s isometry, it holds that
\begin{align}\label{firstderivative}
&\EE|\eta^{\ell i}_T(t, y)|^2\\
\leq \,\,\,&3+3\,\EE\bigg| \sum_{n=1}^{d}\int_t^T \partial_{y_n} b_\ell(\xi_r^{t,y},\mu_r) \eta^{ni}_r(t, y)\,\d r \bigg|^2\\
&\,\,\,\,+3\,\EE\bigg|\sum_{m,n=1}^d\int_t^T\partial_{y_n}  \sigma_{\ell m}(\xi_r^{t,y},\mu_r) \eta^{ni}_r(t, y) \,\d W_r^m \bigg|^2\\
\leq \,\,\,& C_{d,T,M_{\kappa, T}(\mu)}\left(1+\sum_{n=1}^{d} \int_t^T \EE\big|  \eta^{ni}_r(t, y)\big|^2\,\d r+\sum_{m,n=1}^d\int_t^T\EE\big| \eta^{ni}_r(t, y)\big|^2\,\d r\right)\\
\leq \,\,\,& C_{d,T,M_{\kappa, T}(\mu)}\left(1+\int_t^T \EE|\eta_r(t, y)|^2\,\d r\right)\,,
\end{align}
where in the second to last inequality, we use the fact from \cref{assumPDE},~\eqref{Growth1}, that the first derivative of $b, \sigma$ with respect to $y$ only depend on the moment of $\mu$. Analogously, raising both sides of the equation to the fourth power, taking the expectation on both sides of~\eqref{first}, and applying the Burkholder--Davis--Gundy inequality, we obtain
\begin{align}\label{forthmoment}
\EE\big|\eta^{\ell i}_T(t,y)\big|^4\leq C_{d,T,M_{\kappa, T}(\mu)}\left(1+\int_t^T \EE|\eta_r(t, y)|^4\,\d r\right)\,.
\end{align}
Taking square and expectation on both sides of~\eqref{second}, by H$\rm{\ddot{o}}$lder's inequality and It\^{o}'s isometry, we have
\begin{equation}\label{secondderivative}
\EE \big|\gamma^{\ell ij}_T(t, y)\big|^2
 \leq C_{d,T,M_{\kappa, T}(\mu)} \left( \EE\int_t^T |\eta_r(t, y)|^4\,\d r+ \int_t^T \EE|\gamma_r(t, y)|^2\,\d r\right)\,,
\end{equation}
where $|\cdot|$ is understood as the sum of the square of every components.
Summing up in $\ell, i, j$ on both sides of~\eqref{firstderivative},~\eqref{forthmoment}, and ~\eqref{secondderivative}, using Gr\"{o}nwall's inequality, we have 
\begin{align}\label{etasquare}
\EE|\eta_T(t, y)|^2+\EE|\eta_T(t, y)|^4\leq C_{d,T,M_{\kappa, T}(\mu)}\,,
\end{align}
and then we have
\begin{align}\label{gammasquare}
\EE |\gamma_T(t,y)|^2\leq \,\,\,&C_{d,T,M_{\kappa, T}(\mu)}\EE\int_t^T |\eta_r(t, y)|^4\,\d r+C_{d,T,M_{\kappa, T}(\mu)} \int_t^T |\gamma_r(t, y)|^2\,\d r\\
\leq\,\,\,&  C_{d,T,M_{\kappa, T}(\mu)}\,.
\end{align}

 Plugging the above estimates into formulations ~\eqref{firstderivaof},~\eqref{secondderivaof}, H\"older's inequality gives rises that
\begin{align}
|(\partial_{y_i} u)| (t,y)
&\leq\,\,\, \sum_{\ell=1}^{d} \Big[\EE|(\nabla \varphi)_{\ell}(\xi_T^{t,y})|^2\Big]^{\frac{1}{2}}\cdot \Big[\EE|\eta_T^{\ell i}(t, y)|^2\Big]^{\frac{1}{2}}\\
&\leq\,\,\,C_{d,T,M_{\kappa, T}(\mu)}\big(1+|y|^{\kappa-1}\big)\,,
\end{align}
and 
\begin{align}
\big|(\partial_{y_iy_j} u)\big|(t, y)
\leq \,\,\,&\sum_{\ell,m=1}^{d}\Big[\EE|(D^2 \varphi)_{i\ell}(\xi_T^{t,y})|^2\Big]^{\frac{1}{2}} \Big[\EE|\eta_T(t,y)|^4\Big]^{\frac{1}{2}}\\
&+\sum_{n=1}^d \Big[\EE|(\nabla \varphi)_{n}(\xi_T^{t,y})|^2\Big]^{\frac{1}{2}}\cdot  \Big[\EE|\gamma_{T}(t, y)|^2\Big]^{\frac{1}{2}}\\
\leq \,\,\,& \;C_{d,T,M_{\kappa, T}(\mu)}(1+|y|^{\kappa-1}) \, .
\end{align}

\end{proof}
\section{Auxiliary results on  RDEs}\label{sec:calculations}
\begin{lemma}\label{es:compactnessofcontrolledroughpath}
Let $T>0$, $1/3<\alpha<\beta<1/2$ and $\bfX\in \scr{C}^{\alpha}([0,T];\RR^d)$. Then, the following canonical embedding:
\begin{equation}
	\scr{D}_{X}^{2\beta}([0,T];\RR^d)\to\scr{D}_{X}^{2\alpha}([0,T];\RR^d)
\end{equation}
is compact, i.e.\ the bounded sets in the  $\|\cdot\|_{X, 2\beta}$ norm are compact in $\|\cdot\|_{X, 2\alpha}$ norm.
\end{lemma}
\begin{proof}
	Let $(Y^n)$ be a bounded sequence in 
	$\scr{D}_{X}^{2\beta}([0,T];\RR^d)$, namely $\sup_{n\in \bb{N}}|Y_0^n|+|(Y_0^n)'|+\|(Y^n)'\|_{\beta}+\|R^{Y^n}\|_{2\beta}<\infty$, from which it follows $\sup_{n \in \bb{N} }\|Y^n\|_{\beta}< \infty$\,,
	\begin{align}
		&\Big|R^{Y^n}_{s_1, t_1}-R^{Y^n}_{s_2, t_2}\Big|\\
		=\,\,&|Y^n_{s_1,t_1}-(Y_{s_1}^n)'X_{s_1, t_1}-Y^n_{s_2,t_2}+(Y_{s_2}^n)'X_{s_2, t_2}|\\
		\leq \,\,&|t_2-t_1|^{\beta}+|s_2-s_1|^{\beta}+2|s_1-s_2|^{\beta}\|X\|_{0}+\|(Y^n)'\|_{0}(|t_2-t_1|^{\beta}+|s_2-s_1|^{\beta})\\
		\leq \,\,& (1+2\|X\|_0+\|(Y^n)'\|_{0})(|t_2-t_1|^{\beta}+|s_2-s_1|^{\beta})\,,
	\end{align}
i.e.\ for each $n\in \bb{N}$, $R^{Y^n}\in C^{\beta}([0,T]\times [0,T];\RR^d)$. Since the  embedding $C^{\beta}(\Omega)\to C^{\alpha}(\Omega)$ is compact if $\Omega$ is bounded, we are able to find a subsequence (which we continue to parameterise by $n$) along with $Y\in C^{\alpha}([0,T];\RR^d)$, $Y^*\in C^{\alpha}([0,T];\RR^d)$, and $Y^{**}\in C^{\alpha}([0,T]\times [0,T];\RR^d)$, such that as $n\to \infty$,
	\begin{align}
	&\lim_{n\to \infty} \lVert Y^n-Y \rVert_0+\|Y^n-Y\|_{\alpha}=\,0\,,\\
	&\lim_{n\to \infty}\|(Y^{n})'-Y^*\|_{0}+\|(Y^{n})'-Y^*\|_{\alpha}=\,0\,,\\
	&\lim_{n\to \infty} \|R^{Y^n}-Y^{**}\|_{0}+ \|R^{Y^n}-Y^{**}\|_{\alpha} =0\,.
	\end{align}
Here, the $\a$-H\"{o}lder norm for $R^{Y^n}-Y^{**}$ is taken over the product space $[0,T]\times [0,T]$. It suffices to check that $Y\in\scr{D}_{X}^{2\alpha}([0,T]; \RR^d)$. Since for any $s, t\in [0,T]$,
\begin{equation}
	Y_{s,t}^n=(Y_s^n)'X_{s,t}+R_{s,t}^{Y^n}\,,
\end{equation}
	passing to the limit $n\to \infty$, we see that 
	\begin{equation}
	Y_{s,t}=Y_s^*X_{s,t}+Y^{**}_{s,t}\,,
\end{equation}
with $Y^*\in C^{\alpha}([0,T];\RR^d)$ and $Y^{**}\in C^{\alpha}([0,T]\times [0,T];\RR^d)$. In addition,  $\|Y^{**}\|_{2\alpha}<\infty$ since 
\begin{align}
	|Y^{**}_{s,t}|=\lim_{n\to \infty}|R^{Y_n}_{s,t}|\leq C|t-s|^{2\alpha}\, .
\end{align}
	\end{proof}
 \begin{lemma}\label{well-RDE}
Let $T>0$, $\beta\in (1/3,1/2)$ and $\bfX\in\scr{C}^{\beta}([0,T];\RR^d)$.  Assume $\lVert\sigma\rVert_{1,3} <\infty\,.$ Then, the following hold:
\begin{enumerate}[(i)]
    \item 
    Given $\xi\in\RR^d$, there exists a unique solution $(Y,Y')\in\scr{D}_{X}^{2\beta}([0,T];\RR^d)$ to the RDE
    \begin{equation}
        Y_t=\xi+\int_0^t \sigma(s,Y_s)\,\d \bfX_s
        \label{eq:roughSDE}
    \end{equation}
    in the sense of \cref{def:RDEsol}.
    \item 
    The solution is uniformly continuous with respect to the initial value. More precisely, given initial data $\xi,\tilde{\xi}$ with the corresponding solutions $Y, \tilde{Y}$, there exists a constant $C=C(T,\beta,X,\sigma)$ such  that 
    \begin{align}
        \|\tilde{Y}-Y\|_0\leq C|\tilde{\xi}-\xi|\,.
    \end{align}
\end{enumerate}
\end{lemma}
\begin{proof} (i) We first demonstrate that $\sigma(\cdot, Y)\in \scr{D}_X^{2\alpha}([0,T];\RR^d)$ if $Y$ is controlled by $\bfX$ and $\|\sigma\|_{1,3}<\infty$. Indeed,  the Taylor expansion gives us 
\begin{align}
	\sigma(t, Y_t)-\sigma(s,Y_s)&=\,\,\partial_t \sigma(s, Y_s)(t-s)+(\nabla \sigma)(s,Y_s)(Y_t-Y_s)+(D^2 \sigma)(t_0, Y_{t_0})(Y_t-Y_s)^2\\
	&=\,\,\partial_t \sigma(s, Y_s)(t-s)+(\nabla \sigma)(s,Y_s)Y_s'X_{s,t}+(\nabla \sigma)(s,Y_s)R_{s,t}^Y\\
	&\,\,\,\,\,\,\,\,\,\,\,+(D^2 \sigma)(t_0, Y_{t_0})(Y_t-Y_s)^2\,,
\end{align}
where $t_0\in [s,t]$, it follows that the Gubinelli derivative $\sigma(\cdot, Y)'_s=\sigma(s,Y_s)Y_s'=(\nabla\sigma)(s,Y_s)\sigma(s,Y_s)$. 
We now define the map $\mathcal{F}: \scr{D}_X^{2\alpha} \to \scr{D}_X^{2\alpha}$, $  (Y,Y') \mapsto (F,F')$ as follows
\begin{equation}
    \mathcal{F}(Y,Y')=(F,F'):=\left(\xi+\int_0^{\cdot}\sigma(s,Y_s)\d\bfX_s, \sigma(\cdot,Y)\right)\in \scr{D}_X^{2\alpha}\, .
    \label{eq:fixedpointmap}
\end{equation}
 We are looking for a fixed point of this map in the controlled rough path space $\scr{D}_X^{2\alpha}$  which clearly will be a solution of our RDE in the sense of \cref{def:RDEsol}. Furthermore, since $\bfX\in \mathscr{C}^{\beta}$, it turns out that this solution is automatically an element in $\scr{D}_X^{2\beta}$ following similar arguments as in the proof of \cite[Theorem 8.3]{FrizHairer.2014.251}. 
 
 \vspace{2mm}
 \paragraph{Step (a)} We first check that $\mathcal{F}$ leaves the following set, $${B}_T:=\{(Y,Y')\in\scr{D}_X^{2\alpha}:Y_0=\xi,Y_0'=\sigma(0,\xi), \|(Y,Y')\|_{X,2\alpha}\leq 1\} \, ,$$
 invariant for $T \leq T_0$, where $T_0$ is sufficiently small. The invariance would follow if we could prove the following estimate for all $(Y,Y') \in B_T, T \leq 1$,
and
\begin{align}\label{inv2}
&\bigg\|\int_0^{\cdot}\Xi_s\d \bfX_s,\Xi\bigg\|_{X,2\alpha}
\lesssim \,\,T^{1-\alpha}+T^{\beta-\alpha}\,,
\end{align}
where $(\Xi,\Xi'):=(\sigma(\cdot,Y),\sigma(\cdot,Y)')\in \scr{D}_X^{2\alpha}$, and the implicit constant depends on $\alpha,\beta$ and the relevant norms of $\bfX, Y,Y', \sigma$.  Thus, we may choose $T_0$ sufficiently small such that $\mathcal{F}$ leaves $B_{T}$ invariant for all $T \leq T_0$.
We will first prove the estimate 
\begin{align}\label{inv1}
\|\Xi, \Xi'\|_{X,2\alpha}\lesssim\,\, (|Y_0'|+\|Y,Y'\|_{X,2\alpha}) \, ,
\end{align}
for all $(Y,Y') \in B_T, T \leq 1,$ and the implicit constant depends on relevant norms of $\sigma$.
To obtain~\eqref{inv1}, we notice that
\begin{align}
\|\Xi, \Xi'\|_{X,2\alpha}
=\,\,\,&\|(\sigma(\cdot,Y))'\|_{\alpha}+\|R^{\sigma}\|_{2\alpha}\\
=\,\,\,& \|(\nabla \sigma)(\cdot,Y)Y'\|_{\alpha}+\|R^{\sigma}\|_{2\alpha}\,,
\end{align}
where
\begin{align}
&\|(\nabla \sigma)(\cdot,Y)Y'\|_{\alpha} \\
\leq \,\,\,&\sup_{s\neq t\in[0,T]}\frac{|(\nabla \sigma)(t,Y_t)Y_t' -(\nabla \sigma)(s,Y_s)Y_s'|}{|t-s|^{\alpha}}\\
\leq\,\,\, &\|Y'\|_{0}[\nabla \sigma]_{\alpha,0}+\|Y'\|_{0}[\nabla \sigma(\cdot,Y)]_{0,1}\|Y\|_{\alpha}+[\nabla \sigma]_{0,0 }\|Y'\|_{\alpha}\\
\leq \,\,\,&{\|Y'\|_{0}[\nabla \sigma]_{\alpha,0}}+\|Y'\|_{0}[D^2 \sigma]_{0,0}\|Y\|_{\alpha}+[\nabla \sigma]_{0,0}\|Y'\|_{\alpha}\,,
\end{align}
and
\begin{align}
\|R^{\sigma}\|_{2\alpha}=&\,\,\sup_{s\neq t\in[0,T]}\frac{1}{|t-s|^{2\alpha}}|\sigma(t,Y_t)-\sigma(s,Y_s)-(\nabla \sigma(s,Y_s))Y_s'X_{s,t}|\\
=&\,\,\sup_{s\neq t\in[0,T]}\frac{1}{|t-s|^{2\alpha}}|\sigma(t,Y_t)-\sigma(s,Y_s)-(\nabla \sigma(s,Y_s))Y_{s,t}+(\nabla \sigma(s,Y_s))R_{s,t}^Y|\\
\leq&\,\, {[\sigma]_{2\alpha,0}}+\frac{1}{2}[D^2 \sigma]_{0, 0}\|Y\|_{\alpha}^2+[\nabla \sigma]_{0, 0}\|R^Y\|_{2\alpha} \, ,
\end{align}
where in the last inequality we have Taylor expanded $\sigma$.
\worknote{compared to the Lemma 7.3, here, we have  extra terms  as $\sigma$ is time dependent now.}
Hence,
\begin{align}
\|\Xi, \Xi'\|_{X,2\alpha}\leq\,\,\,& \|Y'\|_{0}[\nabla \sigma]_{\alpha, 0}+[\sigma]_{2\alpha,0} \\ 
&+([\nabla \sigma]_{0, 0}+[D^2 \sigma]_{0, 0})(\|Y'\|_{0} \|Y\|_{\alpha}+\|Y\|^2_{\alpha}+\|Y,Y'\|_{X,2\alpha})\,.
\end{align}
Considering that we have the following estimates
\begin{align}
&\|Y'\|_{\infty}\leq \,\,|Y_0'|+\|Y'\|_{\alpha}T^{\alpha}\,,\\
&\|Y\|_{\alpha}\leq\,\, (|Y_0'|+\|Y'\|_{\alpha}T^{\alpha})\|X\|_{\beta}T^{\beta-\alpha}+T^{\alpha}\|R^{Y}\|_{2\alpha}\,,\\
&[\nabla \sigma]_{\alpha, 0}\leq\,\, [\nabla \sigma]_{{\rm{Lip}}, 0}T^{1-\alpha}\,,\\
&[\sigma]_{2\alpha, 0}\leq \,\, [\sigma]_{{\rm{Lip}}, 0}T^{1-2\alpha}\,,
\end{align}
it follows from restricting $T\leq 1$ that
\begin{align}
&\|\Xi, \Xi'\|_{X,2\alpha}\\
\lesssim \,\,&\big([\sigma]_{{\rm{Lip}}, 0}+[\nabla \sigma]_{{\rm{Lip}}, 0}+[\nabla \sigma]_{0, 0}+[D^2 \sigma]_{0, 0}\big) (1+|Y_0'|+\|Y,Y'\|_{X,2\alpha})(|Y_0'|+\|Y,Y'\|_{X,2\alpha})\\
\lesssim \,\,&\big([\sigma]_{{\rm{Lip}}, 0}+[\nabla \sigma]_{{\rm{Lip}}, 0}+[\nabla \sigma]_{0, 0}+[D^2 \sigma]_{0, 0}\big) (|Y_0'|+\|Y,Y'\|_{X,2\alpha})\,.
\end{align}
We now move on to the proof of~\eqref{inv2}. Using Part b) of \cite[Theorem 4.10]{FrizHairer.2014.251}, we know that there exists a constant $C_\alpha$ such that
\begin{align}\label{rough}
\bigg\|\int_0^{\cdot}\Xi_s\d\bfX_s,\Xi\bigg\|_{X,2\alpha}&\leq \|\Xi\|_{\alpha}+C_\alpha(|\Xi_0'|+\|\Xi,\Xi'\|_{X,2\alpha})T^{\beta-\alpha}\,.
\end{align}
By  direct calculation, we have
\begin{align}
\|\Xi\|_{\alpha}
&\,\,\leq \sup_{s\neq t\in[0,T]}\frac{1}{|t-s|^{\alpha}}|\sigma(t,Y_t)-\sigma(s,Y_t)+\sigma(s,Y_t)-\sigma(s,Y_s)|\\
&\,\,\leq [\sigma]_{\alpha, 0}+[\nabla \sigma]_{0,0}\|Y\|_{\alpha}\\
&\,\,\leq [\sigma]_{{\rm{Lip}}, 0}T^{1-\alpha}+[\nabla \sigma]_{0, 0}(|Y_0'|+\|Y'\|_{\alpha}T^{\alpha})\|X\|_{\beta}T^{\beta-\alpha}+T^{\alpha}\|R^{Y}\|_{2\alpha}\\
&\,\,\leq     [\sigma]_{{\rm{Lip}}, 0}T^{1-\alpha}+[\nabla \sigma]_{0, 0}(|\sigma(0,\xi)|+1)T^{\beta-\alpha} \, ,
\end{align}
which completes the desired estimate~\eqref{inv2}.

\vspace{2mm}
\paragraph{Step (b)} We will now show that $\mathcal{F}$ is a contraction in ${B}_T$ for  $T \leq T_0$ chosen to be sufficiently small. We define  $$\Delta_s:=\sigma(s,Y_s)-\sigma(s,\tilde{Y}_s)\, .$$ We now use the same approach as that used for~\eqref{rough} to obtain
\begin{align}
&\|\mathcal{F}(Y,Y')-\mathcal{F}(\tilde{Y},\tilde{Y}')\|_{X,2\alpha}=\,\,\bigg\|\int_0^{\cdot}\Delta_s\d\bfX_s,\Delta\bigg\|_{X,2\alpha}\leq\,\, \|\Delta\|_{\alpha}+C_\alpha \|\Delta,\Delta'\|_{X,2\alpha} T^{\beta-\alpha}\,,
\end{align}
where  we have used the fact that  $\Delta_0=0$ and $\Delta_0'=0$.
The contraction property follows provided that
\begin{align}\label{contr1}
 \|\Delta\|_{\alpha}\leq C\|Y-\tilde{Y},Y'-\tilde{Y}'\|_{X,2\alpha}T^{\beta-\alpha},
\end{align}
and 
\begin{align}\label{contr2}
\|\Delta,\Delta'\|_{X,2\alpha}\leq \|Y-\tilde{Y},Y'-\tilde{Y}'\|_{X,2\alpha}.
\end{align}

To derive (\ref{contr1}), we see from fundamental theorem of calculus and chain rule that
\begin{align}
&\|\Delta\|_{\alpha}\\
=\,\,\,&\sup_{s\neq t\in[0,T]}\frac{|\sigma(t,Y_t)-\sigma(t,\tilde{Y}_t)-(\sigma(s,Y_s)-\sigma(s,\tilde{Y}_s))|}{|t-s|^{\alpha}}\\
\leq\,\,\, &[\sigma]_{{\rm{Lip}},0}\|Y-\tilde{Y}\|_{\alpha}+\|Y-\tilde{Y}\|_{0}[\nabla \sigma]_{{\rm{Lip}},0} T^{1-\alpha}\\
&+\|Y-\tilde{Y}\|_{0}[D^2 \sigma]_{0,0}\frac{|\int_0^1\tilde{Y}_t-\tilde{Y}_s+\theta(Y_t-\tilde{Y}_t-Y_s+\tilde{Y}_s)\d\theta|}{|t-s|^{\alpha}}\\
\leq\,\,\, &[\sigma]_{{\rm{Lip}},0}\|Y-\tilde{Y}\|_{\alpha}+\|Y-\tilde{Y}\|_{0}[\nabla \sigma]_{{\rm{Lip}},0}T^{1-\alpha}+\|Y-\tilde{Y}\|_{0}[D^2 \sigma]_{0,0}\\
&+\frac{|\frac{1}{2}\tilde{Y}_t-\frac{1}{2}\tilde{Y}_s+\frac{1}{2}Y_t-\frac{1}{2}Y_s|}{|t-s|^{\alpha}}\\
\leq \,\,\,& [\sigma]_{{\rm{Lip}},0}\|Y-\tilde{Y}\|_{\alpha}+\|Y-\tilde{Y}\|_{0}[\nabla \sigma]_{{\rm{Lip}},0}T^{1-\alpha}+\|Y-\tilde{Y}\|_{\alpha}^2[D^2 \sigma]_{0,0}\\
\leq\,\,\, & [\sigma]_{{\rm{Lip}},0}\|Y-\tilde{Y}\|_{\alpha}+\|Y-\tilde{Y}\|_{0}[\nabla \sigma]_{{\rm{Lip}},0}T^{1-\alpha}+(\|Y-\tilde{Y}, Y'-\tilde{Y}'\|_{X, 2\alpha})^2[D^2 \sigma]_{0,0}\\
\leq\,\,\, & [\sigma]_{{\rm{Lip}},0}\|Y-\tilde{Y}\|_{\alpha}+\|Y-\tilde{Y}\|_{0}[\nabla \sigma]_{{\rm{Lip}},0}T^{1-\alpha}+\|Y-\tilde{Y}, Y'-\tilde{Y}'\|_{X, 2\alpha}[D^2 \sigma]_{0,0}\,.
\end{align}\worknote{The last inequality is the same as the estimate below (8.10)}
We now go to (\ref{contr2}). Recall that
\begin{align}
\Delta_s&=\sigma(s,Y_s)-\sigma(s,\tilde{Y}_s)\\
&=\left(\int_0^1\nabla \sigma(s,\tilde{Y}_s+\theta(Y_s-\tilde{Y}_s))\d \theta\right)(Y_s-\tilde{Y}_s)\\
&=: G_s \times
H_s\,,
\end{align}
where $G_s=g(s,Y_s,\tilde{Y}_s)=\int_0^1\nabla \sigma(s,\tilde{Y}_s+\theta(Y_s-\tilde{Y}_s))\d\theta$, $H_s=Y_s-\tilde{Y}_s$. In contrast to the  calculation in homogeneous case, here $g$ is a triple variable function and we shall see that it doesn't affect the   computation that much. In view of a straightforward computation, we see that $(G,G')\in\scr{D}_{X}^{2\alpha}$ with $G'=(D_Yg)Y'+(D_{\tilde{Y}}g){\tilde{Y}}'$. In fact, by the fundamental theorem of calculus, one has
\begin{align}
&\int_0^1\nabla \sigma(t,\tilde{Y}_t+\theta(Y_t-\tilde{Y}_t))-\nabla \sigma(s,\tilde{Y}_s+\theta(Y_s-\tilde{Y}_s))\;\d\theta\\
\leq &\int_0^1   [\nabla \sigma]_{{\rm{Lip}},0}\cdot|t-s|\d\theta\\
&+\int_0^1D^2\sigma(s,\tilde{Y}_s+\theta(Y_s-\tilde{Y}_s))(\tilde{Y}_t-\tilde{Y}_s+\theta(Y_t-\tilde{Y}_t-Y_s+\tilde{Y}_s))\d\theta+O(|t-s|^{2\alpha})\\
=&\left(\int_0^1D^2\sigma(s,(1-\theta)\tilde{Y}_s+\theta Y_s)\d\theta\right)(1-\theta) (\tilde{Y}_s'X_{s,t})+\theta(Y_s'X_{s,t})+O(|t-s|^{2\alpha})\,,
\end{align}
which implies $G'=(D_Yg)Y'+(D_{\tilde{Y}}g){\tilde{Y}}'$. Moreover, 
\begin{align}
\|G,G'\|_{X,2\alpha}&=\|(D_Yg)Y'+(D_{\tilde{Y}}g){\tilde{Y}}'\|_{\alpha}+\|R^{G}\|_{2\alpha}\\
&=[D^2 \sigma]_{0,0}\cdot\|Y'\|_{\alpha}+[D^2 \sigma]_{0,0}\cdot\|\tilde{Y}'\|_{\alpha}+\|R^{G}\|_{2\alpha}\\
&\leq  [D^2 \sigma]_{0,0}\cdot\|Y'\|_{\alpha}+[D^2 \sigma]_{0,0}\cdot\|\tilde{Y}'\|_{\alpha}+[D^2 \sigma]_{0,0}\big((1-\theta)\|\tilde{Y}\|_{\alpha}^2+\theta \|Y\|_{\alpha}\big)\\
&\leq  [D^2 \sigma]_{0,0}\cdot\|Y'\|_{\alpha}+[D^2 \sigma]_{0,0}\cdot\|\tilde{Y}'\|_{\alpha}+[D^2 \sigma]_{0,0}(2+|\sigma(0,\xi)|^2)\,.
\end{align}
Since $(G,G')\in\scr{D}_{X}^{2\alpha}$, $(H,H')\in\scr{D}_{X}^{2\alpha}$, then $(GH,(GH)')\in\scr{D}_{X}^{2\alpha}$ with $(GH)'=G'H+GH'$, and the corresponding norm
\begin{align}
&\|(GH,(GH)')\|_{X,2\alpha}\\
=\,\,&\|G'H+GH'\|_{\alpha}+\|R^{GH}\|_{2\alpha}\\
\lesssim\,\,& (|G_0|+|G_0'|+\|G,G'\|_{X,2\alpha})(|H_0|+|H_0'|+\|H,H'\|_{X,2\alpha})\,.
\end{align}
\worknote{\todo{check $\|R^{GH}\|_{2\alpha}$? Is it equal to 0?}}
Since $H_0=Y_0-\tilde{Y}_0=\xi-\xi=0$, $H_0'=Y_0'-\tilde{Y}_0'=\sigma(0,\xi)-\sigma(0,\xi)=0$, we have for all $(Y,Y'), (\tilde{Y},\tilde{Y}')\in B_T$,
\begin{align}
\|\Delta,\Delta'\|_{X,2\alpha} \lesssim\,\,& (|G_0|+|G_0'|+\|G,G'\|_{X,2\alpha})\|H,H'\|_{X,2\alpha}\\
\lesssim\,\,& ([\nabla \sigma]_{0,0}+[D^2 \sigma]_{0,0}(|Y_0'|+|\tilde{Y}_0'|)+[D^2 \sigma]_{0,0})\|Y-\tilde{Y},Y'-\tilde{Y}'\|_{X,2\alpha}\\
\lesssim\,\,& ([\nabla \sigma]_{0,0}+[D^2 \sigma]_{0,0})(1+|\sigma(0,\xi)|)\|Y-\tilde{Y},Y'-\tilde{Y}'\|_{X,2\alpha}\,,
\end{align}
which completes the proof of (\ref{contr2}).
From the above estimates and taking $T_0$ sufficiently small, we see that $F_{T_0}(B_{T_0})\subseteq B_{T_0}$ and deduce that
$$\|F_T(Y,Y')-F_T(\tilde{Y},\tilde{Y}')\|_{X,2\alpha}\leq \frac{1}{2}\|Y-\tilde{Y},Y'-\tilde{Y}'\|_{X,2\alpha}\,.
$$

Because of the above contractivity of $F_{T_0}$, it follows from the fixed point theorem that $F_{T_0}$ has a unique fixed point $(Y,Y')\in B_{T_0}$, which is the unique solution of the RDE on the interval $[0,T_0]$. Here, we notice that $T_0$ can be chosen independent of the initial point and the global solution on $[0,T]$ is constructed by iterating the above well-posedness argument on $[T_0,2T_0], [2T_0,3T_0]$ and so on.

\vspace{3mm}
(ii) Recall that
\begin{align}
	\tilde{Y}-Y=\tilde{\xi}-\xi+\int_0^t \sigma(s,\tilde{Y}_s)-\sigma(s,Y_s)\,\d\bfX_s\,,
\end{align}
which is controlled by $\bfX$. By the estimate \cite[(4.30)]{FrizHairer.2014.251}, we obtain that 
\begin{align}
	\|\tilde{Y}-Y\|_0&\leq\,\,|\tilde{\xi}-\xi|+T^{\beta}\|\tilde{Y}-Y\|_{\beta}\\
	&\leq\,\,|\tilde{\xi}-\xi|+T^{\beta}\left(|Y_0'-\tilde{Y}_0'|+T^{\beta}\|Y,Y';\tilde{Y},\tilde{Y}'\|_{X,2\beta}\right)\,.
\end{align}
Using a similar argument as \cite[Theorem 8.5]{FrizHairer.2014.251} of proving the stability estimate, one is able to prove that $\|\tilde{Y}-Y\|_0\leq C_{T,\beta,X,\sigma}|\tilde{\xi}-\xi|$.
\end{proof}

\vspace{3mm}

\begin{lemma} \label{es:compostionsigma}
Let $\bfX, \tilde{\bfX}\in \scr{C}^{\alpha}([0,T];\RR^d)$, $(Y,Y')\in \scr{D}_{X}^{2\alpha}([0,T];\RR^d)$ and $(\tilde{Y},\tilde{Y}')\in \scr{D}_{\tilde{X}}^{2\alpha}([0,T];\RR^d)$. Assume $\sigma\colon [0, T]\times \RR^d\to \RR^{d\times d}$ satisfy the bounds in~\cref{ass:RDE}. Then 
\begin{align}
	&\|\sigma(\cdot,Y), \sigma(\cdot,Y)';\sigma(\cdot,\tilde{Y}), \sigma(\cdot,\tilde{Y})'\|_{X,\tilde{X}, 2\alpha}\\
	\leq\,\,\, &C_M\left(\|X-\tilde{X}\|_{\alpha}+|Y_0-\tilde{Y}_0|+|Y_0'-\tilde{Y}_0'|+\|Y,Y';\tilde{Y},\tilde{Y}'\|_{X,\tilde{X}, 2\alpha}\right)\,,\label{es:appendixstability1}
\end{align}
and \begin{align}
	&\|\sigma(\cdot,\tilde{Y}), \sigma(\cdot,\tilde{Y})';\tilde{\sigma}(\cdot,\tilde{Y}), \tilde{\sigma}(\cdot,\tilde{Y})'\|_{X,\tilde{X}, 2\alpha} \\[2mm]
	\leq \,\,\,&C_{M}\Big([\nabla(\sigma-\tilde{\sigma})]_{{\rm{Lip}},0}+[\sigma-\tilde{\sigma}]_{{\rm{Lip}},0}+[\nabla(\sigma-\tilde{\sigma})]_{0,0}+[D^2(\sigma-\tilde{\sigma})]_{0,0}\Big)\,.
\end{align}
\end{lemma}	
\begin{proof} Write $C_M(\varepsilon_X+\varepsilon_0+\varepsilon_0'+\varepsilon)$ for the right-hand side of~\eqref{es:appendixstability1}. By the definition of the controlled rough path norm, 
 \begin{align}
 	\|\sigma(\cdot,Y), \sigma(\cdot,Y)';\sigma(\cdot,\tilde{Y}), \sigma(\cdot,\tilde{Y})'\|_{X,\tilde{X}, 2\alpha}=\|\sigma(\cdot,Y)'- \sigma(\cdot,\tilde{Y})'\|_{\alpha}+\|R^{\sigma(\cdot,Y)}-R^{\sigma(\cdot,\tilde{Y})}\|_{2\alpha}\,,
 \end{align}
where 
 \begin{align}
 	&\,\,\,\,\,\,\,\,\,\,\,\|\sigma(\cdot,Y)'- \sigma(\cdot,\tilde{Y})'\|_{\alpha}\\
 	&=\,\,\,\|(\nabla \sigma)(\cdot, Y)Y'-(\nabla \sigma)(\cdot, \tilde{Y})\tilde{Y}'\|_{\alpha}\\
 &\leq \,\,\, \sup_{s\neq t}\left(\frac{1}{|t-s|^{\alpha}}|\nabla \sigma(t, Y_t)(Y_t'-\tilde{Y}'_t-Y_s'+\tilde{Y}'_s)+(\nabla \sigma (t, \tilde{Y}'_t)-\nabla \sigma (s, \tilde{Y}'_s))(Y_s'-\tilde{Y}'_s)|\right)\\
 &\,\,\,\,\,\,\,\,\,\,+\sup_{s\neq t}\left(\frac{1}{|t-s|^{\alpha}}|(\nabla \sigma(t, Y_t)-\nabla \sigma (t, \tilde{Y}_t))(Y_t'-\tilde{Y}'_s)|\right)\\
 &\,\,\,\,\,\,\,\,\,\,+\sup_{s\neq t}\left(\frac{1}{|t-s|^{\alpha}}|\nabla \sigma(t, Y_t)-\nabla \sigma (t, \tilde{Y}_t)-\nabla \sigma(s, Y_s)+\nabla \sigma (s, \tilde{Y}_s)|\cdot |\tilde{Y}'_s|\right)\\
 	&\leq \,\,\,\|\nabla \sigma\|_{0,0}\|Y'-\tilde{Y}'\|_{\alpha}+\big(\|\nabla \sigma\|_{{\rm{Lip}}, 0}+\|\nabla \sigma\|_{0,1}\|Y\|_{\alpha}\big)\|Y'-\tilde{Y}'\|_0\\[2mm]
 	&\,\,\,\,\,\,\,\,\,\,+\|D^2 \sigma\|_{0,0}\|Y-\tilde{Y}\|_{\alpha}\|\tilde{Y}'\|_0+\|D^2 \sigma\|_{{\rm{Lip}}, 0}\|Y-\tilde{Y}\|_0\|\tilde{Y}'\|_0\\
 	&\,\,\,\,\,\,\,\,\,\,+\sup_{s\neq t}\left(\frac{1}{|t-s|^{\alpha}}\|\tilde{Y}'\|_0\Big|\int_0^1 D^2 \sigma(t, \tilde{Y}_t+\theta(Y_t-\tilde{Y}_t))(Y_t-\tilde{Y}_t)\,\d \theta\right.\\
&\,\,\,\,\,\,\,\,\,\,\,\,\,\,\,\,\,\,\,\,\,\,\,\,\,\,\,\,\,\,\,\,\,\,\,\,\,\,\,\,\,\,\,\,\,\,\,\,\,\,\,\,\,\,\,\,\,\,\,\,\,\,\,\,\,\,\,\,\,\,-\int_0^1 D^2 \sigma(s, \tilde{Y}_s+\theta(Y_s-\tilde{Y}_s))(Y_s-\tilde{Y}_s)\,\d \theta\Big|\\
&\,\,\,\,\,\,\,\,\,\,\,\,\,\,\,\,\,\,\,\,\,\,\,\,\,\,\,\,\,\,\,\,\,\,\,\,\,\,\,\,\,\,\,\,\,\,\,\,\,\,\,\,\,\,\,\,\,\,\,\,\,\,\,\,\,\,\,\,\,\,\left.-\int_0^1 D^2 \sigma(s, \tilde{Y}_s+\theta(Y_s-\tilde{Y}_s))(Y_s-\tilde{Y}_s)\,\d \theta\Big|\right)\\
 	&\leq \,\,\,\|\nabla \sigma\|_{0,0}\|Y'-\tilde{Y}'\|_{\alpha}+\big(\|\nabla \sigma\|_{{\rm{Lip}}, 0}+\|\nabla \sigma\|_{0,1}\|Y\|_{\alpha}\big)\|Y'-\tilde{Y}'\|_0\\[2mm]
 	&\,\,\,\,\,\,\,\,\,\,+\|D^2 \sigma\|_{0,0}\|Y-\tilde{Y}\|_{\alpha}\|\tilde{Y}'\|_0+\|D^2 \sigma\|_{{\rm{Lip}}, 0}\|Y-\tilde{Y}\|_0\|\tilde{Y}'\|_0\\
 	&\,\,\,\,\,\,\,\,\,\,+\|D^3 \sigma\|_{0,0}\|Y-\tilde{Y}\|_0\|\tilde{Y}'\|_0
\frac{\big|\int_0^1\int_0^1 \big[\tilde{Y}_t+\theta(Y_t-\tilde{Y}_t)-\tilde{Y}_s-\theta(Y_s-\tilde{Y}_s) \big]\,\d \theta\big|}{|t-s|^{\alpha}}\\
 	&\,\,\,\,\,\,\,\,\,\,+\|D^2 \sigma\|_{0,0}\|Y-\tilde{Y}\|_{\alpha}\|\tilde{Y}'\|_0\\
 	&\leq \,\,\,\|\nabla \sigma\|_{0,0}\|Y'-\tilde{Y}'\|_{\alpha}+\big(\|\nabla \sigma\|_{{\rm{Lip}}, 0}+\|\nabla \sigma\|_{0,1}\|Y\|_{\alpha}\big)\|Y'-\tilde{Y}'\|_0\\
 	&\,\,\,\,\,\,\,\,\,\,+\|D^2 \sigma\|_{0,0}\|Y-\tilde{Y}\|_{\alpha}\|\tilde{Y}'\|_0+\|D^2 \sigma\|_{{\rm{Lip}}, 0}\|Y-\tilde{Y}\|_0\|\tilde{Y}'\|_0\\[2mm]
 	&\,\,\,\,\,\,\,\,\,\,+\|D^3 \sigma\|_{0,0}\|Y-\tilde{Y}\|_0\|\tilde{Y}'\|_0
(\|Y\|_{\alpha}+\|\tilde{Y}\|_{\alpha})+\|D^2 \sigma\|_{0,0}\|Y-\tilde{Y}\|_{\alpha}\|\tilde{Y}'\|_0\,,
 \end{align}
 and 
 \begin{align}
 	&R^{\sigma(\cdot,Y)}-R^{\sigma(\cdot,\tilde{Y})}\\
 	=\,\,\,&\sigma(t, Y_t)-\sigma(s, Y_s)-\nabla \sigma(s, Y_s)Y_s'X_{s,t}-\big[\sigma(t, \tilde{Y}_t)-\sigma(s, \tilde{Y}_s)-\nabla \sigma(s, \tilde{Y}_s)\tilde{Y}_s'X_{s,t}\big]\\
 	=\,\,\,&\sigma(t, Y_t)-\sigma(s, Y_s)-\nabla \sigma(s, Y_s)(Y_{s,t}-R_{s,t}^Y)-\big[\sigma(t, \tilde{Y}_t)-\sigma(s, \tilde{Y}_s)-\nabla \sigma(s, \tilde{Y}_s)(\tilde{Y}_{s,t}-R_{s,t}^{\tilde{Y}})\big]\\
 	\leq\,\,\,&\big(\partial_t \sigma(s, Y_s)-\partial_t\sigma(s, \tilde{Y}_s)\big)(t-s)\\
 	&+\int_0^1 \nabla \sigma (s, Y_s+\theta(Y_t-Y_s))(Y_t-Y_s)\,\d \theta-\int_0^1 \nabla \sigma (s, Y_s)Y_{s,t}\,\d \theta+\nabla \sigma(s, Y_s)R_{s,t}^{Y}
 \\
 & +\int_0^1 \nabla \sigma (s, \tilde{Y}_s+\theta(\tilde{Y}_t-\tilde{Y}_s))(\tilde{Y}_t-\tilde{Y}_s)\,\d \theta-\int_0^1 \nabla \sigma (s, \tilde{Y}_s)\tilde{Y}_{s,t}\,\d \theta+\nabla \sigma(s, \tilde{Y}_s)R_{s,t}^{\tilde{Y}}\\
 	\leq\,\,\,&\big(\partial_t \sigma(s, Y_s)-\partial_t\sigma(s, \tilde{Y}_s)\big)(t-s)+\int_0^1\int_0^1 (D^2 \sigma)(s, Y_s+\gamma\theta(Y_t-Y_s))(Y_{s,t}, Y_{s,t})\theta\,\d \gamma\d \theta\\
 	&-\int_0^1\int_0^1 (D^2 \sigma)(s, \tilde{Y}_s+\gamma\theta(\tilde{Y}_t-\tilde{Y}_s))(\tilde{Y}_{s,t}, \tilde{Y}_{s,t})\theta\,\d \gamma\d \theta+\nabla \sigma(s, Y_s)R_{s,t}^Y-\nabla \sigma(s, \tilde{Y}_s)R_{s,t}^{\tilde{Y}}\,.
 	\end{align}
 	Considering that $T\leq 1$, $\|Y-\tilde{Y}\|_{\alpha}\lesssim \varepsilon_X+\varepsilon_Y,\,\|Y-\tilde{Y}\|_{0}\lesssim \varepsilon_0+\varepsilon_Y,\|\tilde{Y}'\|_{0}+\|Y\|_{\alpha}\leq C_M,$
 	and 
 	$$\frac{\big|\partial_t \sigma(s, Y_s)-\partial_t\sigma(s, \tilde{Y}_s)\big|\cdot|t-s|}{|t-s|^{2\alpha}}\leq \|\partial_t \sigma\|_{{\rm{Lip}}, 0}\|Y-\tilde{Y}\|_0T^{1-2\alpha},$$
 	following the same techniques shown in \cite[Theorem 7.6]{FrizHairer.2014.251}, the estimate~\eqref{es:appendixstability1} holds.
 	
Moreover,  by the definition of the metric $\|\cdot\|_{X,\tilde{X}, 2\alpha}$, we have
\begin{align}
&\|\sigma(\cdot,\tilde{Y})'-\tilde{\sigma}(\cdot,\tilde{Y})'\|_{\alpha}\\
\leq \,\,\, &\sup_{t\neq s}\bigg(\frac{1}{|t-s|^{\alpha}}|\nabla(\sigma-\tilde{\sigma})(t,\tilde{Y}_t)\tilde{Y}_t'-\nabla(\sigma-\tilde{\sigma})(s,\tilde{Y}_t)\tilde{Y}_t'|\bigg)\\
&+\sup_{t\neq s}\bigg(\frac{1}{|t-s|^{\alpha}}|\nabla(\sigma-\tilde{\sigma})(s,\tilde{Y}_t)\tilde{Y}_t'-\nabla(\sigma-\tilde{\sigma})(s,\tilde{Y}_s)\tilde{Y}_t'|\bigg)\\
 &\,\,\, +\sup_{t\neq s}\bigg(\frac{1}{|t-s|^{\alpha}}|\nabla(\sigma-\tilde{\sigma})(s,\tilde{Y}_s)\tilde{Y}_t'-\nabla(\sigma-\tilde{\sigma})(s,\tilde{Y}_s)\tilde{Y}_s'|\bigg)\\
 \leq \,\,\,&[\nabla(\sigma-\tilde{\sigma})]_{{\rm{Lip}},0}T^{1-\a}\cdot \|\tilde{Y}'\|_{0}+[D^2(\sigma-\tilde{\sigma})]_{0, 0} \cdot \|\tilde{Y}\|_{\alpha} \cdot \|\tilde{Y}'\|_{0} +[\nabla(\sigma-\tilde{\sigma})]_{0, 0}\cdot \|\tilde{Y}'\|_{\alpha} \, .
\end{align}
For the norm of  the remainder term, we have the following estimate
\begin{align}
&\|R^{\sigma(\cdot,\tilde{Y})}-R^{\tilde{\sigma}(\cdot,\tilde{Y})}\|_{2\alpha}\\
\leq \,\,\, &\sup_{t \neq s} \left(\frac{1}{|t-s|^{2\a}}|(\sigma-\tilde{\sigma})(t,\tilde{Y}_t)-(\sigma-\tilde{\sigma})(s,\tilde{Y}_t)|\right)\\
&+\sup_{t \neq s}\left(\frac{1}{|t-s|^{2\a}}|(\sigma-\tilde{\sigma})(s,\tilde{Y}_t)-(\sigma-\tilde{\sigma})(s,\tilde{Y}_s)-\nabla(\sigma-\tilde{\sigma})(s,\tilde{Y}_s)\tilde{Y}_s^{'} X_{s,t}|\right)\\
\leq \,\,\,&[(\sigma-\tilde{\sigma})]_{{\rm{Lip}},0}T^{1-2\a}+[D^2(\sigma-\tilde{\sigma})]_{0,0}\cdot\|\tilde{Y}\|_{\alpha}^2+[\nabla(\sigma-\tilde{\sigma})]_{0, 0}\cdot\|R^{\tilde{Y}}\|_{2\alpha}\,.
\end{align}
  
\end{proof}
\endappendix
\bibliographystyle{amsalpha}
\bibliography{zong0625}{}
 \end{document}